\documentclass[10pt,a4paper,final]{article}
\usepackage[left=2.250cm, right=2.250cm, top=3.750cm, bottom=3.75cm]{geometry}

\usepackage{titlesec}

\titleformat*{\section}{\normalsize\bfseries}
\titleformat*{\subsection}{\normalsize\itshape} 
\titleformat*{\subsubsection}{\normalsize\itshape}  
\titleformat*{\paragraph}{\normalsize\itshape}
\titleformat*{\subparagraph}{\normalsize\itshape}
\usepackage[font={footnotesize},labelfont=bf]{caption}

\usepackage{amssymb}
\usepackage{amsmath}

\def \MyFigFolder{./}

\usepackage{bm}
\usepackage{tikz}
\usepackage{tikz-cd}
\usetikzlibrary{arrows,matrix,calc}
\usepackage{accents}
\usepackage{multirow}
\usepackage{tabularx}
\usepackage{adjustbox}
\usepackage{enumerate}

\usepackage[LGR,T1]{fontenc}

\DeclareMathAlphabet{\mathgtt}{LGR}{cmtt}{m}{n}

\usepackage{authblk}
\usepackage{fancyhdr}
\usepackage{changepage}
\newcommand{\address}[2]{\affil[#1]{#2}}
\newcommand{\corref}{\footnote{Corresponding author}}
\newcommand{\MPI}{$\dagger$}
\newcommand{\TUM}{$\ddagger$}

\newcommand\blfootnote[1]{%
  \begin{NoHyper}
  \renewcommand\thefootnote{}\footnote{#1}%
  \addtocounter{footnote}{-1}%
  \end{NoHyper}
}
 
\author[\MPI]{\underline{Francesco Fambri}\corref}
\author[\MPI\TUM]{Eric Sonnendr\"ucker}

\address{\MPI}{Max-Planck-Institut f\"ur Plasmaphysik, Boltzmannstra{\ss}e 2, D-85748 Garching, Germany} 
\address{\TUM}{Technische Universität München, Zentrum Mathematik, Boltzmannstra{\ss}e 3, D-85748 Garching, Germany} 

\date{}    

\newcommand{\keywords}[1]{\textbf{keywords}-- #1}

\allowdisplaybreaks

\usepackage{scalerel}
\usepackage{hyperref}
\usepackage{cleveref}
\usetikzlibrary{svg.path}
\definecolor{orcidlogocol}{HTML}{A6CE39}
\tikzset{
  orcidlogo/.pic={
    \fill[orcidlogocol] svg{M256,128c0,70.7-57.3,128-128,128C57.3,256,0,198.7,0,128C0,57.3,57.3,0,128,0C198.7,0,256,57.3,256,128z};
    \fill[white] svg{M86.3,186.2H70.9V79.1h15.4v48.4V186.2z}
                 svg{M108.9,79.1h41.6c39.6,0,57,28.3,57,53.6c0,27.5-21.5,53.6-56.8,53.6h-41.8V79.1z M124.3,172.4h24.5c34.9,0,42.9-26.5,42.9-39.7c0-21.5-13.7-39.7-43.7-39.7h-23.7V172.4z}
                 svg{M88.7,56.8c0,5.5-4.5,10.1-10.1,10.1c-5.6,0-10.1-4.6-10.1-10.1c0-5.6,4.5-10.1,10.1-10.1C84.2,46.7,88.7,51.3,88.7,56.8z};
  }
}
\newcommand\orcidicon[1]{\href{https://orcid.org/#1}{\mbox{\scalerel*{
\begin{tikzpicture}[yscale=-1,transform shape]
\pic{orcidlogo};
\end{tikzpicture}
}{|}}}}

\patchcmd{\subequations}
  {\theparentequation\alph{equation}}
  {\theparentequation.\alph{equation}}
  {}{}

 \title{Structure preserving hybrid Finite Volume Finite Element method for compressible MHD}

\begin{document}

\maketitle 




\begin{abstract}
In this manuscript we present a novel and efficient numerical method for the compressible viscous and resistive MHD equations for all Mach number regimes. The time-integration strategy is a semi-implicit splitting, combined with a hybrid finite-volume  and finite-element (FE) discretization in space. The non-linear convection is solved by a robust explicit FV scheme,  while the magneto-acoustic terms are treated implicitly in time. The resulting CFL stability condition depends only on the fluid velocity, and not on the Alfv\'enic and acoustic modes.
The magneto-acoustic terms are discretized by compatible FE  based on a continuous and a discrete de Rham complexes designed using Finite Element Exterior Calculus (FEEC). Thanks to the use of FEEC, energy stability, magnetic-helicity conservation and the divergence-free conditions can be preserved also at the discrete level.
A very efficient splitting approach is used to separate the acoustic and the Alfv\'enic modes in such a fashion that the original symmetries of the PDE governing equations are preserved. In this way, the algorithm relies on the solution of linear, symmetric and positive-definite algebraic systems, that are very efficiently handled by the simple matrix-free conjugate-gradient method.
The resulting algorithm showed to be robust and accurate in low and high Mach regimes even at large Courant numbers. Non-trivial tests are solved  in one-, two- and three- space dimensions to confirm the robustness, accuracy, and the low-dissipative and conserving properties of the final algorithm.
While the formulation of the method is very general, numerical results for a second-order accurate FV-FE scheme will be presented.
\end{abstract}






\paragraph{Orcid numbers.}  
F. Fambri \orcidicon{0000-0002-6070-8372},
E. Sonnendr\"ucker \orcidicon{0000-0002-8340-7230}\!\!\!\blfootnote{
\underline{\textit{francesco.fambri@ipp.mpg.de}} (F. Fambri),  
\textit{eric.sonnendruecker@ipp.mpg.de} (E. Sonnendr\"ucker)}.

\newcommand{\id}{{i^*}}
\newcommand{\jd}{{j^*}}
\newcommand{\kd}{{k^*}}
\newcommand{\imain}{i}
\newcommand{\jmain}{j}
\newcommand{\kmain}{k}
\newcommand{\dd}{\,{\rm d} }
\newcommand{\fracp}[2]{\frac{\partial #1}{\partial #2}}
\newcommand{\fract}[2]{\frac{\dd #1}{\dd #2}}
\newcommand{\vecu}{\bm{\mathsf{U}}_e}
\newcommand{\vecv}{\bm{\mathsf{V}}_e}
\newcommand{\vecp}{\bm{\mathsf{P}}}
\newcommand{\vecq}{\bm{\mathsf{Q}}}
\newcommand{\vecm}{\bm{\mathsf{M}}_e}
\newcommand{\vece}{\bm{\mathsf{E}}_f}
\newcommand{\vecb}{\bm{\mathsf{B}}_f} 
\newcommand{\vecbd}{\tilde{\bm{\mathsf{B}}}_e} 
\newcommand{\vechn}{\bm{\mathsf{H}}_{\texttt{n}}}
\newcommand{\veche}{\bm{\mathsf{H}}_{\texttt{e}}}
\newcommand{\vecheb}{\tilde{\bm{\mathsf{H}}}_{\texttt{e}}}
\newcommand{\vech}{\bm{\mathsf{H}}_{\texttt{n}}}
\newcommand{\vecc}{\bm{\mathsf{C}}_f}
 
\newcommand{\vecrho}{\bm{\mathsf{\rho}}}
\newcommand{\matgrad}{\mathbb{G}}
\newcommand{\matgradd}{\tilde{\mathbb{G}}}
\newcommand{\matcurl}{\mathbb{C}}
\newcommand{\matcurld}{\tilde{\mathbb{C}}}
\newcommand{\matdiv}{\mathbb{D}}
\newcommand{\matdivd}{\tilde{\mathbb{D}}}
\newcommand{\matmass}{\mathbb{M}}
\newcommand{\Hmat}{\mathbb{H}}
\newcommand{\mati}{\mathbb{I}}
\newcommand{\matproj}{\mathbb{P}}

\newcommand{\half}{\frac{1}{2}}
 
\newcommand{\Q}{{\mathbf{Q}}}
\newcommand{\F}{\mathbf{F}}
\newcommand{\Fv}{\mathbf{F}_{v}}
\newcommand{\Fp}{\mathbf{F}_{p}}
\newcommand{\Fb}{\mathbf{F}_{b}}
\newcommand{\Fc}{\mathbf{F}_{c}}
\newcommand{\Fd}{\mathbf{F}_{d}}
\newcommand{\Fdv}{\mathbf{F}_{d}}
\newcommand{\Feta}{\mathbf{F}_{\eta}} 
\newcommand{\B}{{\mathbf{B}}}
\newcommand{\C}{{\mathbf{C}}}
\newcommand{\D}{{\mathbf{D}}}
\renewcommand{\v}{\mathbf{v}} 
\renewcommand{\u}{\mathbf{u}} 
\newcommand{\x}{\mathbf{x}} 
\newcommand{\w}{\mathbf{w}} 
\newcommand{\m}{{\mathbf{m}}}
\newcommand{\vecF}{\mathcal{F}^h}

\newcommand{\lessspace}{\!\!\!\!\!\!}
\newcommand{\Avex}[5]{\left( #1 \right)_{#2,#3,#4}^{#5}}
\newcommand{\Avey}[5]{\left( #1 \right)_{#2,#3,#4}^{#5}}
\newcommand{\Avez}[5]{\left( #1 \right)_{#2,#3,#4}^{#5}}
\newcommand{\Avexy}[5]{\left( #1 \right)_{#2,#3,#4}^{#5}}
\newcommand{\Avexz}[5]{\left( #1 \right)_{#2,#3,#4}^{#5}}
\newcommand{\Aveyx}[5]{\left( #1 \right)_{#2,#3,#4}^{#5}}
\newcommand{\Aveyz}[5]{\left( #1 \right)_{#2,#3,#4}^{#5}}
\newcommand{\Avezx}[5]{\left( #1 \right)_{#2,#3,#4}^{#5}} 
\newcommand{\Avezy}[5]{\left( #1 \right)_{#2,#3,#4}^{#5}} 
\newcommand{\f}{{\mathbf{f}}} 

\newcommand{\levi}{\mathgtt{e}} 
\newcommand{\CFL}{\texttt{CFL}}  
\newcommand{\nvar}{n_{\text{var}}}

\newcommand{\blue}[1]{{\color{blue}#1}}

\newcommand{\dtilde}[1]{\accentset{\approx}{#1}}  

\newcommand{\cell}{T}
\newcommand{\MHD}{\text{MHD}}


\section{Introduction}
\label{sec:intro}

MHD is a single fluid model that describes the macroscopic dynamics of electrically conducting fluids. Main research branches are related to the study of equilibrium configurations and related instabilities (statics) and the description of low-frequency disruptive processes (dynamics) \cite{Biskamp1993,GoedbloedPoedts,Freidberg:book}.
 The design of effective inertial or magnetic confinement fusion reactors for civilian (clean) energy production is probably one of the most noble considered applications of plasma physics. 
Other physical and industrial applications  may concern the plasma flow driven by strong magnetic field, e.g. in the magnetosphere of neutron stars, the physics of solar flares, accretion disks and relativistic jets in astrophysics, as well as the simulation of plasma thrusters for propulsion engineering in the design spacecrafts and satellites.

A spectral analysis of the ideal MHD equations identifies up to eight different characteristic wave speeds, expressed in terms of the convective, the Alfv\'en, the \emph{slow-} and \emph{fast-} magneto-sonic wave speeds. 
For a hyperbolic PDE solver, \emph{explicit} methods are constraint by the Courant-Friedrich-Levy condition, which is a {necessary} condition on the maximum time-step size to get stability and convergence. 
In terms of computational efficiency, explicit methods are unbeatable for one single time-update.  However, depending on the physical regime, large characteristic speeds may impose extremely small time-steps. In these situations, the simulation becomes computationally very expensive with eventually a loss of accuracy due to excessive numerical stabilization (diffusion) or accumulation of truncation error.
On the contrary, \emph{implicit}-integration methods imposes a \emph{global} algebraic equation involving all the  degrees of freedom without  CFL stability restrictions. 
This is actually the reason why implicit time integration is nowadays considered an essential tool for the numerical simulation of stabilized plasma in magnetic confinement fusion devices 
 \cite{AYDEMIR1985,HARNED1986,NIMROD1999,SOVINEC2004,Chacon2008}.
On the other hand, the draw-back of fully implicit methods is that they generate  a \emph{fully-coupled} \emph{global} algebraic system for all the degrees of freedom, that may become very inefficient to be solved, especially for nonlinear problems.
In all physical situations where fast modes are mainly associated to slowly varying flows, 
semi-implicit methods may guarantee stable and accurate solutions at the minimal computational costs even for large Courant numbers, e.g. for low plasma-beta, low Mach numbers or (nearly) incompressible flows. 
For decades, semi-implicit methods have been explored to design efficient schemes for  plasmas physics, see e.g. \cite{SCHNACK1987,HarnedSchnack1986,Lerbinger1991,LIONELLO1999}. In the works by  Chacon (2008) \cite{Chacon2008} and L\"utjens and Luciani (2010) \cite{LUTJENS2010}, efficient non-linear implicit Newton  methods have been designed. Consider the well-known review-paper by Jardin (2012) \cite{Jardin2012} about implicit methods for magnetically confined plasma.
  
The study of semi-implicit methods  had wide application also in many other areas of research and engineering.
Pressure based semi-implicit schemes on staggered date back to the \lq\lq marker and cell\rq\rq\, method of Harlow and Welch (1965) \cite{markerandcell}, and became soon widely used, see  e.g. \cite{chorin1,chorin2, 
patankarspalding,BellColellaGlaz,
vanKan,HirtNichols}.
A remarkable example is given by the very efficient class of pressure-based semi-implicit methods that have been designed and theoretically analyzed by Casulli et al. in the field of gravity-driven viscous free-surface and sub-surface flows, see \cite{CasulliCompressible, Casulli1990,CasulliCheng,Casulli1999,CasulliCattani, BrugnanoCasulli,BrugnanoCasulli2,CasulliZanolli2012}, guaranteeing conservation properties, positivity and unconditional stability at large Courant numbers.
These methods have been recently extended to other PDE systems, ranging from the simulation of blood flow in compliant arterial vessels \cite{CasulliDumbserToro,Blood3D2014}, non-Newtonian continuum mechanics \cite{SIGPR,Peshkov2021}, incompressible and compressible MHD \cite{SIMHD,hybridhexa}, eventually to the high-order discontinuous Galerkin finite-element methods on structured \cite{DumbserCasulli2013,FambriDumbser}, unstructured \cite{TavelliDumbser2014b,TavelliDumbser2015,TavelliDumbser2016} and adaptive grids \cite{AMRDGSI,Fambri2020}.
Semi-implicit method are very popular  in ocean-  as well as in  climate-modeling  and numerical weather prediction    \cite{Giraldo2005,TumoloBonaventura} where slow and fast, small and large scales need to be correctly approximated in one single simulation.
 Whenever using flux-splitting methods for  multi-scale problems, consistency with the analytical system should be preserved in all the asymptotic limits, 
 see \cite{KLEIN1995,Jin1999,Klein2001} introducing the concept of asymptotic-preserving (AP) schemes and multi-scale modeling. 
Then, the first pressure-methods were made suitable for higher-Mach number and nonlinear shock-dominated flows, see e.g. \cite{ParkMunz2005,Kwatra2009,Smolarkiewicz2009,Cordier2012,DumbserCasulli2016,RussoAllMach}.
The design of asymptotic preserving schemes  and multi-scale modeling are still nowadays hot-topic of research,  see \cite{Vater2018,Benacchio2019,Schmid2021}  for atmospheric modeling or  \cite{BoscarinoRusso2009,BoscarinoPareschiRusso,Kucera2022,DimarcoPareschi2017,BoscarinoPareschiRusso2017} for applications to fluid and kinetic models.

It is a well-known fact that preserving some physical invariants in a closed discrete system can improve the stability properties of a numerical method. Since the first intuitions of Arakawa (1966) \cite{arakawa1966computational}, only in the last few decades preservation of differential structures  and vector calculus identities  \cite{Bossavit.1998.ap,hiptmair2002finite,hirani2003discrete,pavlov2011structure,gawlik2011geometric} have become an essential tool to guarantee conservation of the invariants of the original PDE system.
 Among the years, mimetic finite-differences (MFD) \cite{Lipnikov2014,Brezzi2005}, mixed finite-elements \cite{Brezzi1985,BBF2013}, discrete exterior calculus (DEC) \cite{Hiptmair2001}, and finally the finite element exterior calculus (FEEC) \cite{arnold_falk_winther_2006_anum,arnold2018finite}, providing a practical framework for high-order discretizations of the de Rham complex, were deeply investigated and became increasingly popular in many research fields, e.g. see \cite{alonso1996error,hiptmair2002finite,cotter2012mixed,kraus2017gempic,palha_mass_2017,brokenFEECNS}. Relevant advances in the numerical analysis with FEEC with applications to incompressible (Hall) MHD are \cite{HU2021,GAWLIK2022,Laakmann2023}, while we refer to \cite{Holderied2020,HOLDERIED2021,Li2024} for application to MHD-kinetic models.

In this paper we introduce a novel structure-preserving semi-implicit method for the viscous and resistive magnetohydrodynamics (MHD) equations at all-Mach regimes. 
 The final algorithm is shown to be sufficiently stable to preserve a nonlinear MHD equilibrium in the long-time with low dissipation, as well as sufficiently robust and accurate to properly solve nonlinear shock dominated MHD flows in multiple space dimensions.
A semi-implicit splitting method is built with a \emph{hybrid} finite-volume/finite-element approach: the conservation  and robustness properties of explicit finite-volume (FV) schemes  are exploited for the nonlinear advection, while relying on an implicit time discretization  with compatible finite-elements of the magneto-acoustic terms.
Thanks to the use of FEEC, some geometrical structures and physical invariants are preserved: energy stability, conservation of magnetic-helicity and the divergence-free condition of the magnetic field hold at the discrete level. 
The final scheme will be subject to a Courant-Friedrich-Levy stability condition which depend only on the fluid velocity, and not on the Alfv\'enic and acoustic modes.   
Moreover, the adopted compatible finite-elements for the magneto-acoustics are based on discrete and continuous de Rham complexes in such a fashion that the original symmetries of the governing PDE  are reflected to the resulting discrete algebraic systems.  
\smallskip

\noindent
The rest of the paper is organized as follows:
\begin{enumerate}[(i)]
\item first, the governing PDE system is presented in Section \ref{sec:model}; 
\item  the numerical method is described in Section \ref{sec:numerics}: more specifically, the flux-splitting (\ref{sec:split}-\ref{eq:char}); the  finite-volume and chosen  compatible finite-elements spaces and the associated discrete  operators (\ref{sec:FE}-\ref{sec:diff}); the chosen semi-implicit time integration method (\ref{sec:time}); 
\item  the numerical scheme is validated against a non-trivial set of MHD problems in one-, two- and three-space dimensions in different Mach and plasma-beta regimes, see  Section \ref{sec:numerics}; 
\item  and finally, Section \ref{sec:conclusions} concludes the paper with a summary and some final comments and remarks.
\end{enumerate}

\section{The compressible   MHD model}
\label{sec:model}

We consider here the following compressible viscous and resistive MHD model in primitive variables 
\begin{align}
  \fracp{\rho}{t} + \nabla\cdot (\rho \mathbf{u}) &=0 \label{eq:MHD-rho} \\
  \fracp{ \rho\mathbf{u}}{t} + \nabla\cdot \left(\rho\mathbf{u}\otimes \mathbf{u}  \right) + \nabla p -(\nabla\times \mathbf{B})\times \mathbf{B} &= \nabla \cdot  \mu \left( \nabla \mathbf{u} + \nabla \mathbf{u}^T - \frac{2}{3} \left( \nabla \cdot \mathbf{u} \right) \mathbf{I} \right)   \\
  \fracp{ \mathbf{B} }{t} + \nabla \times \mathbf{E} &=0 \\
  \fracp{p}{t} + \nabla\cdot(p \mathbf{u}) + (\gamma -1 ) p \nabla\cdot \mathbf{u} &=\nabla \cdot ( \kappa \nabla T ) \label{eq:MHD-p}   
\end{align}
where the electric field is $ \mathbf{E}= - \mathbf{u}\times \mathbf{B} + \eta\nabla \times \mathbf{B}$, $\eta$ being a constant resistivity.  Here, $\mathbf{I}$ is the identity matrix, $T$ is the temperature, 
$\mu$ and $\kappa$ are the kinematic viscosity and thermal conductivity. 
The temperature $T$ is defined by the thermal equation of state $T=T(p,\rho)$, see Sec. \ref{sec:EOS}.
The equation system (\ref{eq:MHD-rho}-\ref{eq:MHD-p}) can be also written in the conservative form
\begin{equation}
\partial_t \Q + \nabla \cdot \left( \Fv+\Fd + \Fp + \Fb +  \Feta\right) = 0,
\end{equation}
where the array of conservative variables and the tensor of conservative fluxes are defined as

\begin{equation}\begin{array}{rl}
& 	\Q:= \left( \begin{array}{c} \rho \\ \rho \mathbf{u} \\ \rho E \\ \mathbf{B} \end{array}  \right); \quad 
	 \F_v := \begin{pmatrix} \rho \mathbf{u} \\ 
				\rho \mathbf{u} \otimes \mathbf{u}  \\	
	\frac 12 \mathbf{u}  \rho \mathbf{u}^2    \\
	0  \end{pmatrix} 
	\quad 	\Fp:= \begin{pmatrix}   0 \\    p \mathbf{I} 
\\    \mathbf{u} \rho  h   \\  0 \end{pmatrix} 
			 \end{array}  
\end{equation}

\begin{equation}\begin{array}{l}
	 \Fd := -\left( \begin{array}{c} 
0 \\ 
 \mu \left( \nabla \mathbf{u} + \nabla \mathbf{u}^T - \frac{2}{3} \left( \nabla \cdot \mathbf{u} \right) \mathbf{I} \right) \\
\mu \mathbf{u}  \left( \nabla \mathbf{u} + \nabla \mathbf{u}^T - \frac{2}{3} \left( \nabla \cdot \mathbf{u} \right) \mathbf{I} \right) + \kappa \nabla T  \\
0    
\end{array} \right) 
 \end{array}  \label{eq:MHDFd}.
\end{equation} 

\begin{equation}\begin{array}{rl}
& \Fb:= \left( \begin{array}{c} 0 \\   \frac{\mathbf{B}^2}{2}  \mathbf{I} - \mathbf{B} \otimes \mathbf{B} 
	\\  \mathbf{u}^T \mathbf{B}^2   - \mathbf{B}^T \left( \mathbf{u} \cdot \mathbf{B} \right)   \\ \mathbf{B} \otimes \mathbf{u} - \mathbf{u} \otimes \mathbf{B} \end{array}  \right) \end{array}, \quad \begin{array}{l}
	 \Feta := -\left( \begin{array}{c} 
0 \\ 
0 \\ 
 \frac{\eta}{4 \pi} \mathbf{B}^T \left( \nabla \mathbf{B} - \nabla \mathbf{B}^T \right) \\
\eta \left( \nabla \mathbf{B} - \nabla \mathbf{B}^T \right)    
\end{array} \right) \end{array}  \label{eq:MHD-Fb}
\end{equation} 
where $\rho E = {p}/({\gamma-1}) +  \mathbf{u}\cdot \mathbf{m}/2 $ is the total energy density, $h$ is the specific enthalpy  $ h = e + \frac{p}{\rho}$, $e$ the specific internal energy.
In this work we may indicate the three spatial components of the tensor fluxes as $\F_\alpha:=\begin{pmatrix} \mathbf{f}_\alpha, \mathbf{g}_\alpha, \mathbf{h}_\alpha   \end{pmatrix}$, where $\alpha$ indicates one of the splitting components.
\subsection{Equation of state of ideal gases}
\label{sec:EOS}
In this work we used equation of state for ideal gases. This simply means that we choose a thermal and caloric equations of state 
\begin{equation} 
\label{eq:thermcal.ideal} 
 \frac{p}{\rho} = (c_p-c_v) T, \qquad \textnormal{ and } \qquad e = c_v T,  
\end{equation}  
where $c_v$ and $c_p$  are  the heat capacities   at constant volume and at constant pressure, respectively. Combining these equations
\eqref{eq:thermcal.ideal} one may derive the corresponding equation of state for the internal energy
\begin{equation}
 e = e(p,\rho) = \frac{p}{(\gamma-1) \rho},
\end{equation} 
with $\gamma = c_p/c_v$ is the ratio of specific heats, so that the specific enthalpy reads $h=(p / \rho)\gamma/(\gamma-1)$. 
Note that the equations of sate for ideal gases is linear in the pressure $p$. 
%

\section{Numerical method}
\label{sec:numerics} 
\subsection{Splitting} 
\label{sec:split}
The numerical method is based on a splitting of the model into convective, acoustic and Alfv\'enic parts:
\begin{itemize}
  \item The convection-diffusion part reads:
  \begin{align}
    \fracp{\rho}{t} + \nabla\cdot (\rho \mathbf{u}) &=0  \label{eq:FVini}\\
    \fracp{ \rho\mathbf{u}}{t} + \nabla\cdot (\rho\mathbf{u}\otimes \mathbf{u} )  &= \nabla \cdot  \mu \left( \nabla \mathbf{u} + \nabla \mathbf{u}^T - \frac{2}{3} \left( \nabla \cdot \mathbf{u} \right) \mathbf{I} \right)  \\
    \fracp{ \mathbf{B} }{t}  &=0 \\ 
    \fracp{p}{t} &=\nabla \cdot ( \kappa \nabla T )   \label{eq:FVend}
  \end{align}
  \item the acoustic part reads:
  \begin{align}
    \fracp{\rho}{t} &=0 \label{eq:splitac-1}\\
    \fracp{\rho\mathbf{u}}{t}+ \nabla p  &=0 \label{eq:splitac-2a}\\
    \fracp{ \mathbf{B} }{t} &=0 \label{eq:splitac-3}\\
     \fracp{p}{t} +   \nabla\cdot(p \mathbf{u}) + (\gamma -1 ) p \nabla\cdot \mathbf{u}  &= 0\label{eq:splitac-4a}
  \end{align}
  We can rewrite this system by taking the dot product of \eqref{eq:splitac-2a} with $\mathbf{u}$ and adding it to \eqref{eq:splitac-4a}. 
  Thus in this step $\rho$ and $ \mathbf{B}$ are constant and the momentum $\mathbf{m}=\rho \mathbf{u}$ and $p$ are the solutions of the following \emph{conservative} system
  \begin{align}
    \fracp{\mathbf{m}}{t}+ \nabla p  &=0 \label{eq:splitac-2}\\
    \frac{\gamma-1}{\gamma}\fracp{}{t} \left( \frac{p}{\gamma-1} + \frac {1}{2} \mathbf{u}\cdot \mathbf{m} \right) + \nabla\cdot (\frac{p}{\rho}\mathbf{m}) &=0 \label{eq:splitac-4}
  \end{align}
	where we notice the definition of the \emph{total energy density} $\rho E= {p}/({\gamma-1}) +  \mathbf{u}\cdot \mathbf{m}/2 $;
  \item the Alfv\'enic part is here formulated in \emph{non-conservative} form and it reads:
  \begin{align}
    \fracp{\rho}{t}  &=0 \label{eq:splitalf-1}\\
    \fracp{ \mathbf{m}}{t} - (\nabla\times \mathbf{B})\times \mathbf{B} &=0 \label{eq:splitalf-2}\\
    \fracp{ \mathbf{B} }{t} +  \nabla\times ( -\mathbf{u}\times \mathbf{B} + \eta\nabla \times\mathbf{B}) &=0 \label{eq:splitalf-3}\\
    \fracp{p}{t} &=0 \label{eq:splitalf-4}
  \end{align}
  \item and the resistive part reads:
  \begin{align}
    \fracp{\rho}{t}  &=0 \label{eq:splitres-1}\\
    \fracp{ \mathbf{m}}{t}  &=0 \label{eq:splitres-2}\\
    \fracp{ \mathbf{B} }{t} +  \nabla\times ( \eta\nabla \times\mathbf{B}) &=0 \label{eq:splitres-3}\\
    \fracp{p}{t} &=0 \label{eq:splitres-4}
  \end{align}
\end{itemize}

\subsection{Characteristics} \label{eq:char}
In the design of robust numerical schemes for hyperbolic PDE, it is important to consider the spectrum of the characteristic wave speed. In a simplified one-dimensional model of idea plasma flow with constant $B_x$, one may compute the eight eigenvalues explicitly as 
\begin{equation}
\lambda_{1,8}^{\MHD} = u \mp c_f, \quad \lambda_{2,7}^{\MHD} = u \mp c_a, \quad \lambda_{3,6}^{\MHD} = u \mp c_s, \quad \lambda_4^{\MHD} = u, \quad \lambda_5^{\MHD} = 0, 
\label{eq:eval.full} 
\end{equation} 
where the speeds of the Alfv\'en ($c_a$), the slow ($c_s$) and the fast ($c_f$) magnetosonic waves are defined by
\begin{equation}\begin{array}{rl}
 &c_a = B_x / \sqrt{4 \pi \rho}, \\
 &c_s^2 = \half \left( b^2+c^2 - \sqrt{(b^2+c^2)^2-4 c_a^2 c^2} \right), \\
 &c_f^2 = \half \left( b^2+c^2 + \sqrt{(b^2+c^2)^2-4 c_a^2 c^2} \right),  \end{array}
\label{eq:mhd.wavespeeds} 
\end{equation} 
respectively. We introduced $b^2 = \mathbf{B}^2/\rho$, while  $c$ is the adiabatic sound speed and, by specifying an equation of state for ideal gases, it is defined by the relation $c^2 = \gamma p/ \rho$.
We may cite that the presented method can be straightforward extended to a more general class of equations of state in the form $p=p(e,\rho)$, see \cite{DumbserCasulli2016,SIMHD,Fambri20}. In that case, we have $c^2 = \partial p / \partial \rho + p / \rho^2 \partial p / \partial e$.
As stated above, we are following the PDE splitting proposed by \cite{Fambri20}, trying to divide the eigenspectrum between  \emph{convection modes}
\begin{equation}\label{eq:lambdac}
\begin{aligned}
\lambda^{v,x}_{1,2,3,4} =  0, \quad  	 \lambda^{v,x}_{5,6,7,8} = v_x, 
\end{aligned}
\end{equation}
 \emph{acoustic modes} 
\begin{equation}\label{eq:lambdap}
\begin{aligned}
\lambda^{p,x}_{1,8} = \frac{1}{2} \left( v_x \mp \sqrt{v_x^2 + 4 c^2 } \right), \quad  	\lambda^{p,x}_{2,3,4,5,6,7}= 0, 
\end{aligned}
\end{equation}
and \emph{Alfv\'enic modes} 
\begin{equation}\label{eq:lambdab}
\begin{aligned}
 &\lambda^{b,x}_{1,8} = \frac{1}{2} \left( v_x \mp \sqrt{v_x^2 + 4 \left | \B \right |^2/\rho } \right),  \quad  	 \lambda^{b,x}_{2,7} = \frac{1}{2} \left( v_x \mp \sqrt{v_x^2 + 4 B_x^2/\rho} \right), \\
& \lambda^{b,x}_{3,4,5,6} =  0.  
\end{aligned}
\end{equation}

\subsection{The Finite Element spaces and operators} \label{sec:FE}
\subsubsection{The  de Rham Complex}
\label{sec:DeRham}
The convective part will be discretized with a classical conservative Finite Volume scheme, where the discrete unknowns $\rho$, $\mathbf{u}$, $\mathbf{m}$, $ \mathbf{B}$ and $p$ are characterized by their cell average. 
On the other hand, the last two parts will be discretized with compatible Finite Elements coming from Finite Element Exterior Calculus (FEEC). These are based on the following commuting diagram involving a continuous and a discrete de Rham complex as well as commuting projectors:
\begin{equation}\label{CGL}
  \begin{tikzpicture}[baseline=(current  bounding  box.center)]
    \matrix (m) [matrix of math nodes,row sep=3em,column sep=4em,minimum width=2em] {
             H^1(\Omega) & H(\mathrm{curl},\Omega) & H(\mathrm{div},\Omega) & L^2(\Omega)
      \\
             V_0 & V_1 & V_2 & V_3
      \\
      };
    \path[-stealth]
      (m-1-1) edge node [left]  {$\Pi_0$} (m-2-1)
      (m-1-2) edge node [left]  {$\Pi_1$} (m-2-2)
      (m-1-3) edge node [left]  {$\Pi_2$} (m-2-3)
      (m-1-4) edge node [left]  {$\Pi_3$} (m-2-4)
      (m-1-1) edge node [above] {$\mathrm{grad}$} (m-1-2)
      (m-1-2) edge node [above] {$\mathrm{curl}$} (m-1-3)
      (m-1-3) edge node [above] {$\mathrm{div}$}  (m-1-4)
      (m-2-1) edge node [above] {$\mathrm{grad}$} (m-2-2)
      (m-2-2) edge node [above] {$\mathrm{curl}$} (m-2-3)
      (m-2-3) edge node [above] {$\mathrm{div}$}  (m-2-4)
      ;
      \path[-stealth]
      ($(m-2-2.south)+(-0.35,+0.15)$) edge node [below] {$\mathrm{div}_w$}  ($(m-2-1.south)+(+0.35,+0.15)$)
      ($(m-2-3.south)+(-0.35,+0.15)$) edge node [below] {$\mathrm{curl}_w$}  ($(m-2-2.south)+(+0.35,+0.15)$)
      ($(m-2-4.south)+(-0.35,+0.15)$) edge node [below] {$\mathrm{grad}_w$}  ($(m-2-3.south)+(+0.35,+0.15)$)
      ;
  \end{tikzpicture}
\end{equation}
where the weak discrete differential operators are defined by 
\begin{align}
  \int \mathrm{grad}\, p_h \cdot \mathbf{v}_h \dd \mathbf{x} &= - \int \mathrm{div}_w\,\mathbf{v}_h \cdot p_h \dd \mathbf{x}, \quad p_h\in V_0, \mathbf{v}_h\in V_1, \label{eq:defgradw} \\
  \int \mathrm{curl}\, \mathbf{u}_h \cdot \mathbf{B}_h \dd \mathbf{x} &= \int \mathrm{curl}_w\, \mathbf{B}_h  \cdot  \mathbf{u}_h  \dd \mathbf{x}, \quad  \mathbf{B}_h \in V_2, \mathbf{u}_h\in V_1, \label{eq:defcurlw} \\
  \int q_h \mathrm{div}\, \mathbf{B}_h  \dd \mathbf{x} &= - \int \mathrm{grad}_w\, q_h \cdot \mathbf{B}_h  \dd \mathbf{x}, \quad q_h\in V_3, \mathbf{B}_h\in V_2, \label{eq:defdivw} 
\end{align}

  We use as Finite Element spaces tensor products space with polynomials of degree $k$ in each variable for the first space.
  This corresponds to the Lagrange Finite Element $V_0= \mathbb{Q}_p$. 
  Then $V_1$ is the space of edge elements,   $V_2$ is the space of face elements and $V_3$ the space of discontinuous elements of degree $p-1$. For higher degree the degrees of freedom of $V_0$ are the values at the tensor product Gauss-Lobatto mesh. In this paper we will consider only $p=1$, in which case the degrees of freedom in $V_0$ are the vertex values, the degrees of freedom in $V_1$ are the edge integrals, the degrees of freedom in $V_0$ are the face integrals and the degrees of freedom in $V_3$ are the cell averages.

  We will keep the scalar valued quantities $\rho_n, p_n, E_n \in V_0$, which corresponds to their \emph{nodal} value and location at the \emph{nodes} (or vertices) of the cells.
  
\begin{figure}\centering
	\begin{tabular}{cc} 
	\includegraphics[width=0.25\textwidth]{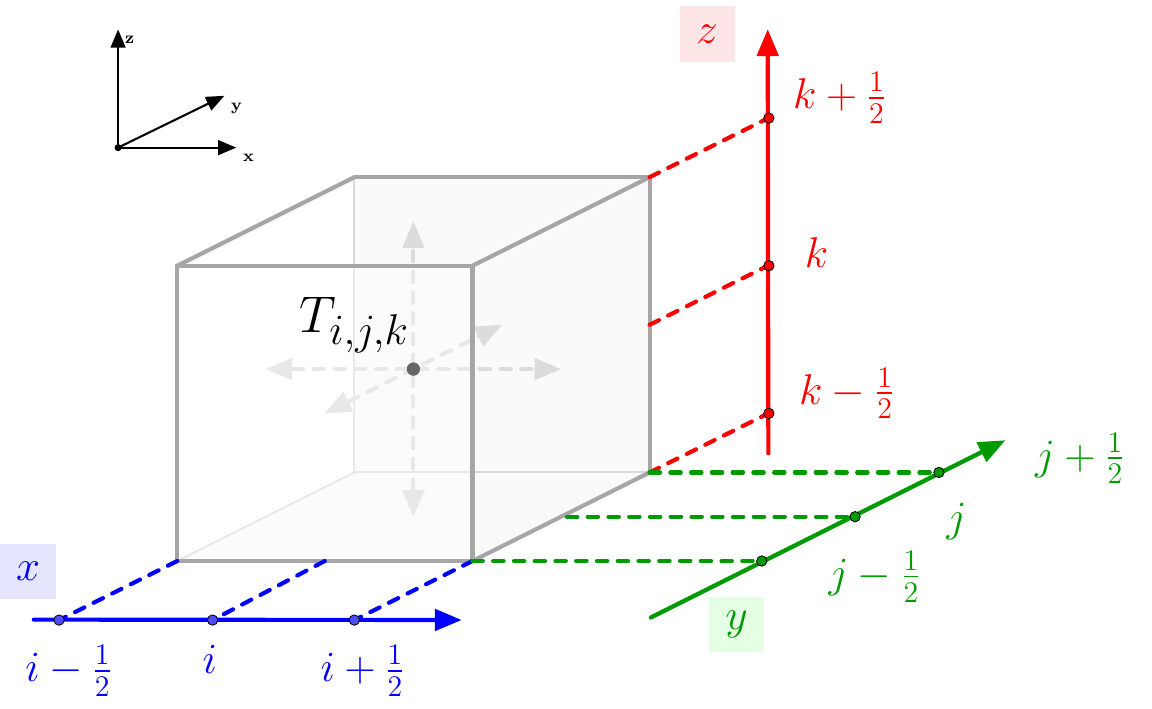}&
	\includegraphics[width=0.21\textwidth]{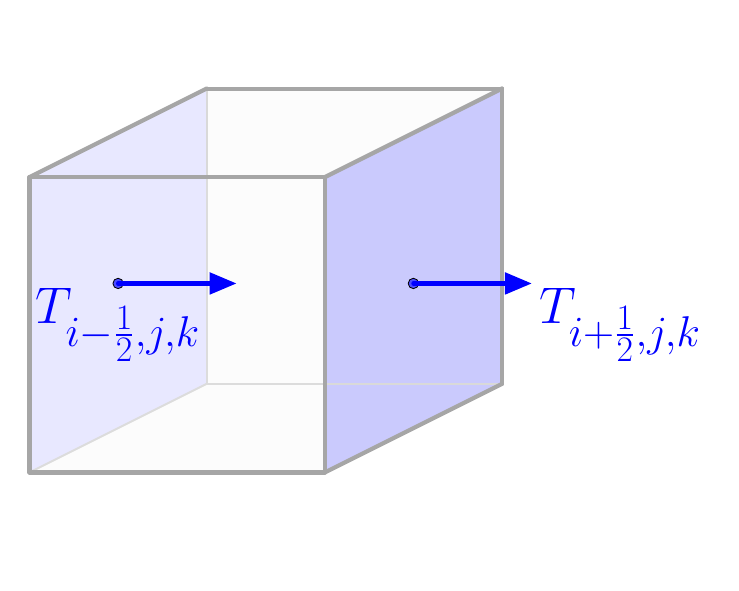} 
	\includegraphics[width=0.21\textwidth]{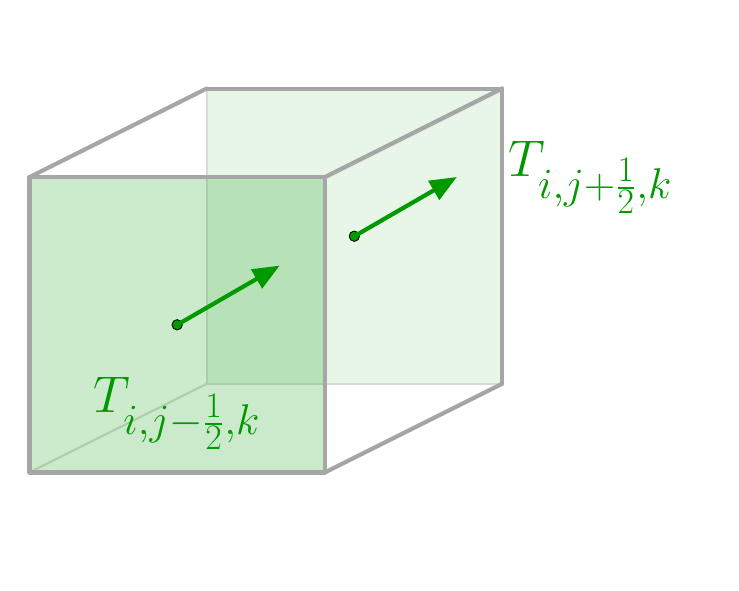} 
	\includegraphics[width=0.21\textwidth]{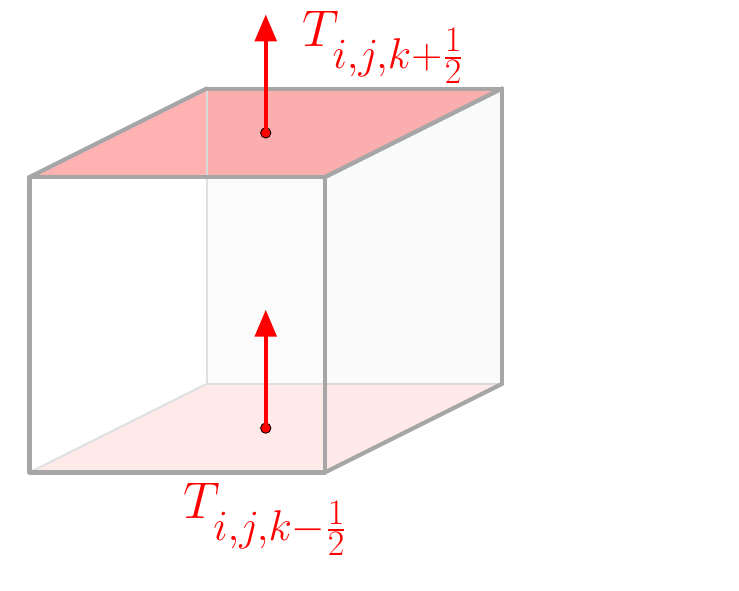} \\
	\includegraphics[width=0.25\textwidth,trim={0cm 0cm 2cm 0cm},clip]{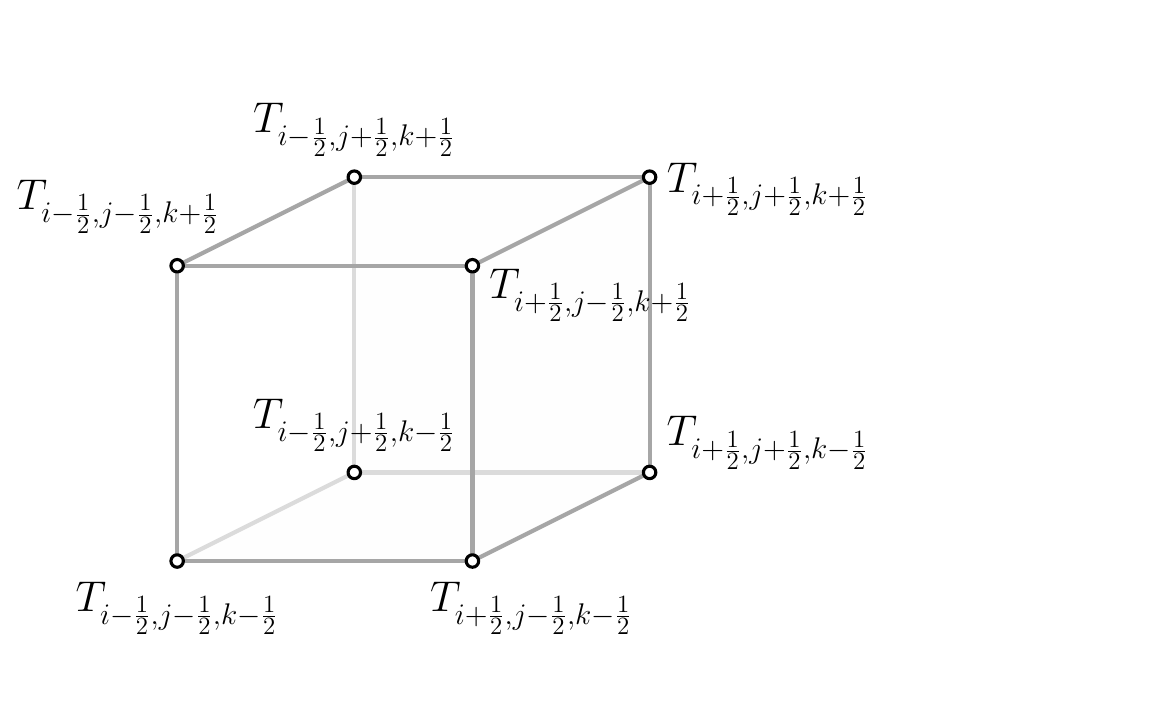}  &
	\includegraphics[width=0.21\textwidth]{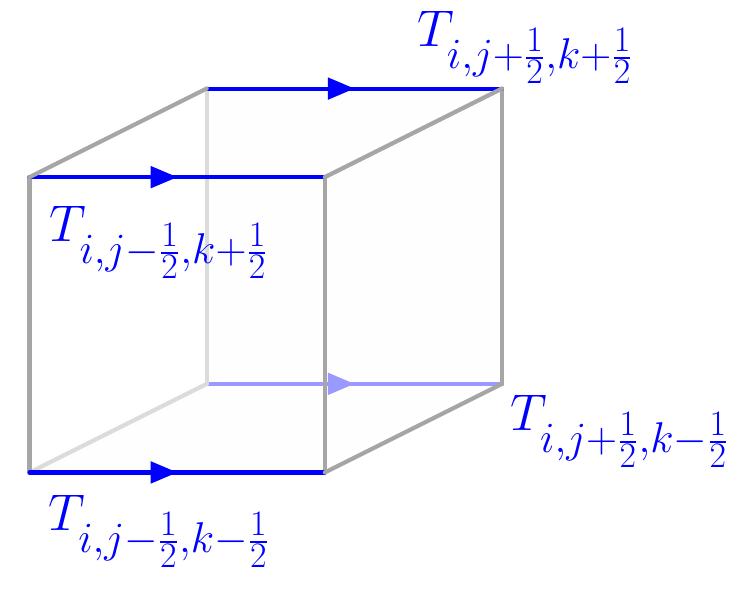}
	\includegraphics[width=0.21\textwidth]{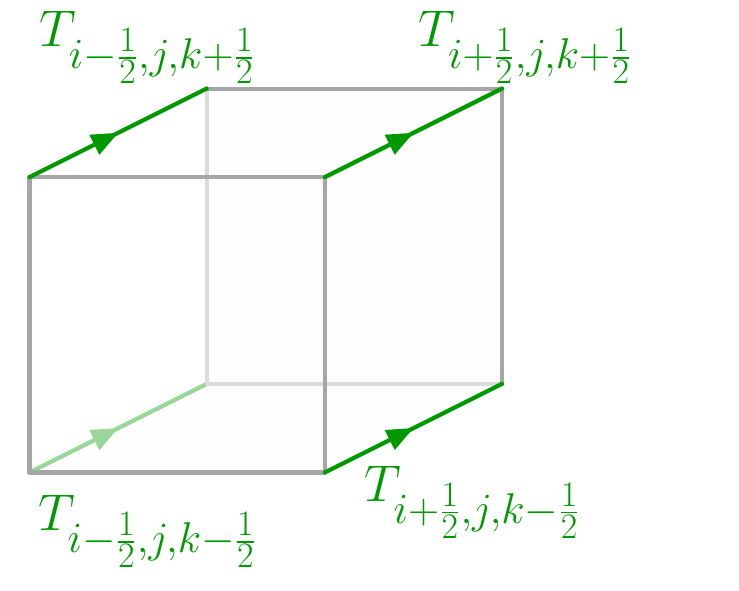} 
	\includegraphics[width=0.21\textwidth]{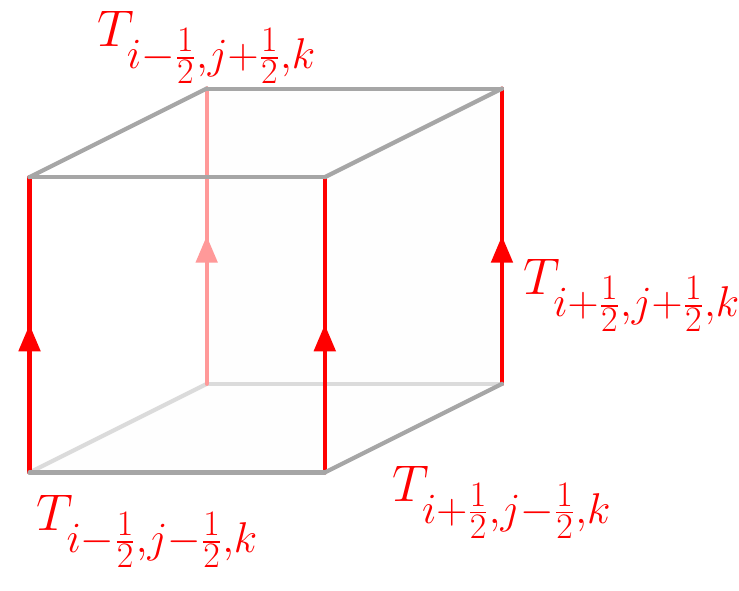}
	\end{tabular}
	\caption{Sketch of the location of the FEEC degrees of freedom: $V^0$ scalar variables are described by \emph{nodal} valued functions (bottom-left); $V^1$ vector variables by \emph{edge}-valued functions (bottom-right figures for the three spatial components, respectively); $V^2$ vector variables by \emph{face}-valued functions (top-right figures for the three components, respectively); $V^3$ scalar variables by \emph{barycentric}-valued functions (top-left figure). Note and remember that the \emph{nodal}-point corresponds to the \emph{dual} barycenters of the corresponding \emph{dual}-grid.}
	\label{fig:dualmeshes}  
	\end{figure}

 For the vector valued quantities we make the following choice: $\mathbf{u}_e,\,\mathbf{m}_e\in V_1$ (edge based) and $\mathbf{B}_f\in V_2$ which is the face based Finite Element space. Other choices are possible. In Fig. \ref{fig:dualmeshes}, a sketch of the location of the FEEC degrees of freedom in the adopted low-order approximation

\begin{figure} 
	\centering
\includegraphics[width=0.7\textwidth]{\MyFigFolder/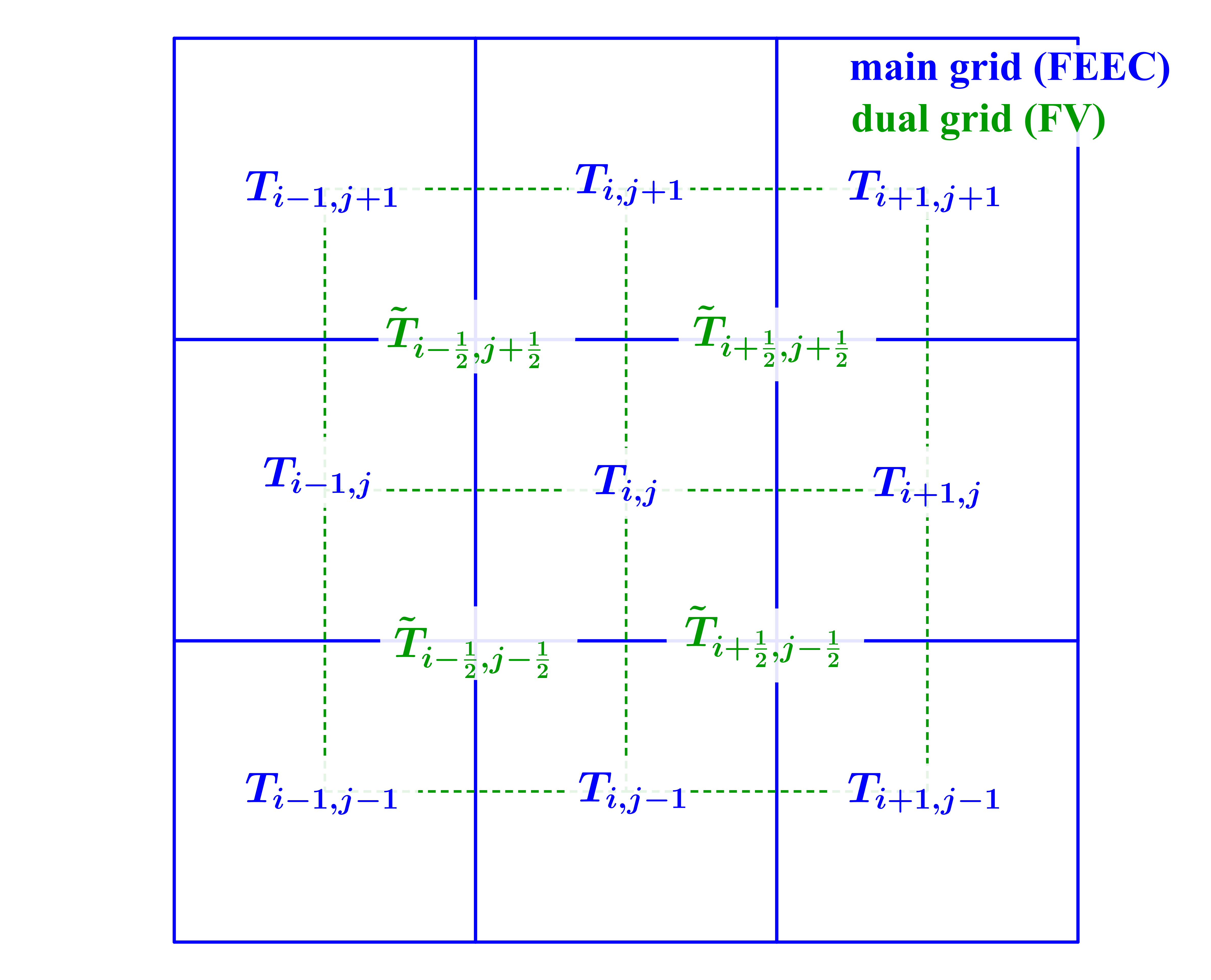}
\caption{Sketch of the main (blue) and the barycentric dual (green) grids.}
\label{fig:dualgrid}
\end{figure}

 We will thus have two different representations of all discrete quantities, their Finite Element representation characterized by the appropriate degrees of freedom and their Finite Volume representation characterized by their cell averages. In this work, we introduce a barycentric dual grid (see Fig.\ref{fig:dualgrid}) so that momentum $\mathbf{m}_e$  and velocities $\mathbf{u}_e$, whose degrees of freedom are edge-based, have continuous normal component unique value for flux at dual-cell interface. For Cartesian uniform meshes, this corresponds to a node-base staggered mesh. As an example, now the horizontal velocity has a unique value on the x-edges of the primal grid, as well as on the x-faces of the dual grid. 
In this way, the dual-barycenters are located in the nodes of the primal grid, that is the natural location for elements of $V_0$.
 
 Hence, the degrees of freedom for $\rho$ and $E$ are the average values on each dual-cell as in the Finite Volume case:
  $\rho_{i+\half,j+\half,k+\half}$, $p_{i+\half,j+\half,k+\half}$ which correspond to the cell averages on the dual-grid, which are given by the Finite Volume step. For convenience, we may assume the following notation for the indexes of the dual grid (later called \lq\lq dual indexes\rq\rq\,): $i^*:=i+\half$.
  
  The degrees of freedom for variables in $V_1$, e.g. $\mathbf{u}_e$,  are the edge averages on each edge of the mesh
 \begin{align}
   (u_{e,x})_{i,j+\frac 12,k+\frac 12} &= \frac{1}{\Delta x}\int_{x_{i-\frac 12}}^{x_{i+\frac 12}} u_e(x,y_{j+\frac 12},z_{k+\frac 12}) \dd x, \\
   (u_{e,y})_{i+\frac 12,j,k+\frac 12} &= \frac{1}{\Delta y}\int_{y_{j-\frac 12}}^{y_{j+\frac 12}} u_e(x_{i+\frac 12},y,z_{k+\frac 12}) \dd y, \\
   (u_{e,z})_{i+\frac 12,j+\frac 12,k} &= \frac{1}{\Delta z}\int_{z_{k-\frac 12}}^{z_{k+\frac 12}} u_e(x_{i+\frac 12},y_{j+\frac 12},z) \dd z,
 \end{align}
 and finally the degrees of freedom for variables in $V_2$, e.g. $ \mathbf{B}_f$,  are the face averages \textit{e.g.}
 \begin{align}
  (B_{f,x})_{i+\frac 12,j,k} &= \frac{1}{\Delta y\Delta z} \int_{y_{j-\frac 12}}^{y_{j+\frac 12}}\int_{z_{k-\frac 12}}^{z_{k+\frac 12}} B_f(x_{i+\frac 12},y,z) \dd y \dd z, \\
  (B_{f,y})_{i,j+\frac 12,k} &= \frac{1}{\Delta x\Delta z}\int_{x_{i-\frac 12}}^{x_{i+\frac 12}}\int_{z_{k-\frac 12}}^{z_{k+\frac 12}} B_f(x,y_{j+\frac 12},z) \dd x \dd z, \\
  (B_{f,z})_{i,j,k+\frac 12} &= \frac{1}{\Delta x\Delta y}\int_{x_{i-\frac 12}}^{x_{i+\frac 12}}\int_{y_{j-\frac 12}}^{y_{j+\frac 12}} B_f(x,y,z_{k+\frac 12}) \dd x \dd y.
\end{align}

We then define vectors containing all the degrees of freedom stacked in some order of the indices: 
for the elements of $V_0$  we define the node associated to cell $i,j,k$ to be the node at $(x_{i-\frac 12}, y_{j-\frac 12}, z_{k-\frac 12})$. 
On our periodic mesh, there is exactly one such node for each cell and so the node based quantities are characterized, e.g., by the degrees of freedom   $\vecp=(p_{i-\half,j-\half,k-\half})_{i,j,k}$,   $\vecrho=(\rho_{i-\half,j-\half,k-\half})_{i,j,k}$ for pressure and density.
For the vector based unknowns each $(i,j,k)$ index contains a small vector with three components:
\begin{equation}\label{eq:V2dof}
(\vecb)_{i,j,k} = \begin{pmatrix}
  (B_{f,x})_{i-\frac 12,j,k} \\ (B_{f,y})_{i,j-\frac 12,k} \\ (B_{f,z})_{i,j,k-\frac 12}
\end{pmatrix} ~~~~
(\vecu)_{i,j,k} = \begin{pmatrix}
  (u_{e,x})_{i,j-\frac 12,k-\frac 12}  \\ (u_{e,y})_{i-\frac 12,j,k-\frac 12} \\ (u_{e,z})_{i-\frac 12,j-\frac 12,k}
\end{pmatrix}
\end{equation}
 and finally for the barycentric degrees of freedom of variables in $V_3$, we define the node associated to cell $i,j,k$ to be $\vech = (h_{i,j,k})_{i,j,k}$.

We denote by $\Sigma_0, \Sigma_1, \Sigma_2, \Sigma_3$ the corresponding set of degrees of freedom, so that
$\vecp,\vecrho\in \Sigma_0$, $\vecu\in \Sigma_1$, $\vecb\in \Sigma_3$, $\vech\in \Sigma_3$.

\subsection{Discrete gradient, divergence and curl matrices} \label{sec:diff}

Let us now define the discrete gradient matrix associated to our degrees of freedom, denoted by $\mathbb{G}$, that applied to a vector $\vecp\in \Sigma_0$, describing a scalar variable with $N$ components, to 
a vector $\vecu$ describing a vector variable with altogether $3N$ components. We have
$$  
\begin{pmatrix}
  (u_{e,x})_{i,j-\frac 12,k-\frac 12}  \\ (u_{e,y})_{i-\frac 12,j,k-\frac 12} \\ (u_{e,z})_{i-\frac 12,j-\frac 12,k}
\end{pmatrix} =: (\vecu)_{i,j,k}=(\matgrad\vecp)_{i,j,k} = 
\begin{pmatrix}
  \frac{p_{i+\frac 12,j-\frac 12,k-\frac 12}- p_{i-\frac 12,j-\frac 12,k-\frac 12}}{\Delta x} \\ \frac{p_{i-\frac 12,j+\frac 12,k-\frac 12}- p_{i-\frac 12,j-\frac 12,k-\frac 12}}{\Delta y} \\ \frac{p_{i-\frac 12,j-\frac 12,k+\frac 12}- p_{i-\frac 12,j-\frac 12,k-\frac 12}}{\Delta z} 
\end{pmatrix}.
$$
Then we can define the discrete curl matrix for an edge based quantity $\vecu$ as
\begin{equation}
  (\matcurl\vecu)_{i,j,k} = 
  \begin{pmatrix}
    \frac{(u_{e,z})_{i-\frac 12,j+\frac 12,k}- (u_{e,z})_{i-\frac 12,j-\frac 12,k}}{\Delta y}
    - \frac{(u_{e,y})_{i-\frac 12,j,k+\frac 12}- (u_{e,y})_{i-\frac 12,j,k-\frac 12}}{\Delta z} \\
    \frac{(u_{e,x})_{i,j-\frac 12,k+\frac 12}- (u_{e,x})_{i,j-\frac 12,k-\frac 12}}{\Delta z}
    - \frac{(u_{e,z})_{i+\frac 12,j-\frac 12,k}- (u_{e,z})_{i-\frac 12,j-\frac 12,k}}{\Delta x}\\
    \frac{(u_{e,y})_{i+\frac 12,j,k-\frac 12}- (u_{e,y})_{i-\frac 12,j,k-\frac 12}}{\Delta x}
    - \frac{(u_{e,x})_{i,j+\frac 12,k-\frac 12}- (u_{e,x})_{i,j-\frac 12,k-\frac 12}}{\Delta y}
  \end{pmatrix}
\end{equation}
and the discrete divergence matrix for a face based quantity $\vecb$
\begin{multline} \label{eq:divdof}
  (\matdiv\vecb)_{i,j,k} = \left( \frac{(B_{f,x})_{i+\frac 12,j,k}- (B_{f,x})_{i-\frac 12,j,k}}{\Delta x}
+ \frac{(B_{f,y})_{i,j+\frac 12,k}- (B_{f,y})_{i,j-\frac 12,k}}{\Delta y} 
\right. \\ \left.
+ \frac{(B_{f,z})_{i,j,k+\frac 12}- (B_{f,z})_{i,j,k-\frac 12}}{\Delta z}
  \right) 
\end{multline}

Let us now compute the matrices corresponding to the weak differential operators as defined by \eqref{eq:defgradw}, \eqref{eq:defdivw}, denoted respectively by $\matgradd$, $\matcurld$ and $\matdivd$.
Denoting the mass matrix of $V_i$ by $\matmass_i$ for $i=0,1,2,3$. The definition relations immediately yield
\begin{equation}\label{eq:dualderivatives}
  \matmass_2\matgradd = - \matdiv^\top  \matmass_3, ~~
  \matmass_1\matcurld = \matcurl^\top  \matmass_2, ~~
  \matmass_0\matdivd = - \matgrad^\top  \matmass_1.
\end{equation}

We can summarize the action of these discrete differential operators on the degrees of freedom in the following de Rham complexes
\begin{equation}\label{discDiag}
  \begin{tikzpicture}[baseline=(current  bounding  box.center)]
    \matrix (m) [matrix of math nodes,row sep=3em,column sep=4em,minimum width=2em] {
             \Sigma_0 & \Sigma_1 & \Sigma_2 & \Sigma_3
      \\
      };
    \path[-stealth]
      (m-1-1) edge node [above] {$\matgrad$} (m-1-2)
      ($(m-1-2.south)+(-0.35,+0.15)$) edge node [below] {$\matdivd$}  ($(m-1-1.south)+(+0.35,+0.15)$)
      (m-1-2) edge node [above] {$\matcurl$} (m-1-3)
      ($(m-1-3.south)+(-0.35,+0.15)$) edge node [below] {$\matcurld$}  ($(m-1-2.south)+(+0.35,+0.15)$)
      (m-1-3) edge node [above] {$\matdiv$}  (m-1-4)
      ($(m-1-4.south)+(-0.35,+0.15)$) edge node [below] {$\matgradd$}  ($(m-1-3.south)+(+0.35,+0.15)$)
      ;
  \end{tikzpicture}
\end{equation}
where we have $\matcurl\matgrad=\matdiv\matcurl=0$ and also on the dual sequence
$\matcurld\matgradd=\matdivd\matcurld=0$.

\subsubsection{The discrete mass matrices}

We consider in this paper low order Finite Elements with $V_0 = \mathbb{Q}_1$ the tensor product piecewise linear element and $V_3 = \mathbb{Q}_0$, piecewise constant.

Let us denote by $\varphi_{i-\frac 12}(x)$ the linear hat functions associated to the node $x_{i-\frac 12}$ in 1D
\begin{equation}\label{eq:phi}
\varphi_{i-\frac 12}(x) = \left|\begin{matrix}
  \frac{x-x_{i-\frac 32}}{\Delta x} & x_{i-\frac 32}\leq x \leq x_{i-\frac 12} \\
  -\frac{x-x_{i+\frac 12}}{\Delta x} & x_{i-\frac 12}\leq x \leq x_{i+\frac 12} \\
  0 & \mbox{else}
\end{matrix}
  \right.
\end{equation}
and by $\chi_i$ the piecewise constant basis functions in 1D
\begin{equation}\label{eq:chi}
\chi_i(x) = \left|\begin{matrix}
  1 & x_{i-\frac 12}\leq x \leq x_{i+\frac 12} \\
  0 & \mbox{else}
\end{matrix}
  \right.
\end{equation}
Then functions in $V_0,V_1, V_2,V_3$ can be expressed in the following tensor product bases.
For any $p_h\in V_0$
\begin{equation}
  p_h(x,y,z)= \sum_{i,j,k} p_{i-\frac 12,j-\frac 12,k-\frac 12} \varphi_{i-\frac 12}(x)\varphi_{j-\frac 12}(y)\varphi_{k-\frac 12}(z),
\end{equation}
for any $ \mathbf{u}_h \in V_1$
\begin{align*}
  u_{h,x}(x,y,z) &= \sum_{i,j,k} (u_{e,x})_{i,j-\frac12,k-\frac12} \chi_i(x)\varphi_{j-\frac 12}(y)\varphi_{k-\frac 12}(z), \\
  u_{h,y}(x,y,z) &= \sum_{i,j,k} (u_{e,y})_{i-\frac12,j,k-\frac12}  \varphi_{i-\frac 12}(x)\chi_j(y)\varphi_{k-\frac 12}(z),\\
  u_{h,z}(x,y,z) &= \sum_{i,j,k} (u_{e,z})_{i-\frac12,j-\frac12,k} \varphi_{i-\frac 12}(x)\varphi_{j-\frac 12}(y)\chi_k(z),
\end{align*}
for any $ \mathbf{B}_h \in V_2$
\begin{align*}
  B_{h,x}(x,y,z) &= \sum_{i,j,k} (B_{f,x})_{i-\frac12,j,k} \varphi_{i-\frac 12}(y)\chi_j(x)\chi_k(z), \\
  B_{h,y}(x,y,z) &= \sum_{i,j,k} (B_{f,y})_{i,j-\frac12,k} \chi_i(x)\varphi_{j-\frac 12}(y)\chi_k(z), \\
  B_{h,z}(x,y,z) &= \sum_{i,j,k} (B_{f,z})_{i,j,k-\frac12} \chi_i(x)\chi_j(y)\varphi_{k-\frac 12}(z),
\end{align*}
and for any $q_h\in V_3$
\begin{equation}
  q_h(x,y,z)= \sum_{i,j,k} q_{i,j,k} \chi_i(x)\chi_j(y)\chi_k(z).
\end{equation}
We observe that due to the chosen basis functions, the coefficients of the expansion are exactly the value of the function at the corresponding point:
$$p_h(x_{i-\frac12},y_{j-\frac12},z_{k-\frac12}) =  p_{i-\frac12,j-\frac12,k-\frac12}, ~~ B_{h,x}(x_{i-\frac12},y_j,z_k)=(B_{f,x})_{i-\frac12,j,k}, 
$$
$$u_{h,x}(x_i,y_{j-\frac12},z_{k-\frac12}) = (u_{e,x})_{i,j-\frac12,k-\frac12}, ~~~q_h(x_{i},y_{j},z_{k}) = q_{i,j,k}$$ and similarly for the other components.

We will frequently need to express $ \mathbf{m}_h =  \rho_h\mathbf{u}_h$. Using the expressions with the basis functions, we have 
\begin{equation}\begin{aligned}
  m_{h,x}(x_{i},y_{j+\half},z_{k+\half}) & =  \rho_h(x_{i},y_{j+\half},z_{k+\half}) u_{h,x}(x_{i},y_{j+\half},z_{k+\half}) \\
& =\frac12(\rho_{i-\half,j+\half,k+\half}+ \rho_{i+\half,j+\half,k+\half}) (u_{e,x})_{i,j+\half,k+\half}
\end{aligned}
\end{equation}

and similarly for the other components so that 
\begin{align}
  (m_{e,x})_{i,j+\frac12,k+\frac12 } &= \frac12(\rho_{i-\frac12,j+\frac12,k+\frac12}+ \rho_{i+\frac12,j+\frac12,k+\frac12}) (u_{e,x})_{i,j+\frac12,k+\frac12} \label{eq:rhoux}\\
  (m_{e,y})_{i+\frac12,  j,k+\frac12} &= \frac12(\rho_{i+\frac12,j-\frac12,k+\frac12}+ \rho_{i+\frac12,j+\frac12,k+\frac12})(u_{e,y})_{i+\frac12,j,k+\frac12} \label{eq:rhouy}\\
  (u_{e,z})_{i+\frac12,j+\frac12,k }  &= \frac12(\rho_{i+\frac12,j+\frac12,k-\frac12}+ \rho_{i+\frac12,j+\frac12,k+\frac12})(u_{e,z})_{i+\frac12,j+\frac12,k } \label{eq:rhouz}
\end{align}
As averages of this type will be needed several times, we introduce the following notation on the vectors having these components:
$\vecm = \bar{\vecrho}\vecu$. Such products will also appear in weak formulations of the form $\vecm \cdot \vecm = \bar{\vecrho}\vecu\cdot \vecm$. Then a change of index in the sum will yield  
$$\vecm \cdot \vecm = \bar{\vecrho}\vecu\cdot \vecm = \vecrho \cdot \overline{\vecu\vecm},$$
where the $N$ components of $\overline{\vecu\vecm}$ are the averages over the components
\begin{align*}
  &(\overline{\vecu\vecm})_{i+\frac12,j+\frac12,k+\frac12} = \\ &  \frac 12 \left[ (u_{e,x})_{i,j+\frac12,k+\frac12 }(m_{e,x})_{i,j+\frac12,k+\frac12 } + (u_{e,x})_{i+1,j+\frac12,k+\frac12 }(m_{e,x})_{i+1,j+\frac12,k+\frac12 }
  \right.  \\ 
	& \left.  	+(u_{e,y})_{i+\frac12,j,k+\frac12}(m_{e,y})_{i+\frac12,j,k+\frac12}+(u_{e,y})_{i+\frac12,j+1,k+\frac12}(m_{e,y})_{i+\frac12,j+1,k+\frac12}
\right.   \\
& \left.   + (u_{e,z})_{i+\frac12,j+\frac12,k}(m_{e,z})_{i+\frac12,j+\frac12,k} +(u_{e,z})_{i+\frac12,j+\frac12,k+1}(m_{e,z})_{i+\frac12,j+\frac12,k+1}
  \right] 
\end{align*}
We will need the weighted mass matrix for tensor product scalar weights $\alpha(x)\beta(y)\gamma(z)$.
This maintains the tensor product structure so that only 1D integrals will need to be computed.
We will compute them by approximating the matrices involving products of $\varphi_{i-\frac 12}$ by the two-points Gauss-Lobatto quadrature (trapezoidal rule) on each element and the matrices involving products of $\chi_i$ by the one-point Gauss-Legendre quadrature (midpoint rule). Both are second order so that no loss of order is induced by the quadrature and all these so-called lumped mass matrix are diagonal.
We have 
\begin{equation}\label{eq:lamped}
\begin{aligned}
  \int \alpha(x)\varphi_{i-\frac 12}(x)\varphi_{l-\frac 12}(x) \dd x & \approx \alpha(x_{i-\frac 12}) \Delta x ~ \delta_{il}, \\
  \int \alpha(x)\chi_i(x)\chi_l(x) \dd x & \approx \alpha(x_{i})\Delta x ~ \delta_{il}.
	\end{aligned}
\end{equation} 
The mass matrices $\matmass_i$ of the spaces $V_i$ consists of product of those on the diagonals.
In particular when the weight is one, we have for $N$ grid points, denoting by $\mati_N$ the identity matrix of size $N$
\begin{equation} \label{eq:massmat}
\matmass_0 = \matmass_3 = \Delta x \Delta y\Delta z\mati_N, ~~~
\matmass_1 = \matmass_2 = \Delta x \Delta y\Delta z\mati_{3N}.
\end{equation} 
This implies that the relation between primal and dual differential operators \eqref{eq:dualderivatives} becomes in this case
\begin{equation}\label{eq:dualderivativesQ1}
  \matgradd = - \matdiv^\top, ~~
  \matcurld = \matcurl^\top, ~~
  \matdivd = - \matgrad^\top .
\end{equation}

Lastly, in order to treat nonlinear terms arising from products, projections to the discrete finite element spaces will be needed.
In order to treat in structure preserving manner the mixed products of the type $\int b\cdot (a\times c) \dd x$, which vanishes when $b=c$ for example, we will follow the idea introduced in \cite{HU2021} by applying the projections both on $a\times c$ and $c$ itself even though the later is not needed for our framework.

As this feature will be needed in the $V_1$ space for helicity conservation, we will introduce it only for $P_1$ the projection onto the $V_1$ space. However, the same could be done for other spaces. A natural projection operator $P_1$ is the $V_1$-orthogonal projection.
However in order to make the code more efficient, the mass matrix which appears in the discrete version can be lumped. In order to possibly account for this, 
we will  define the $V_1$ projection $P_1$ using the discrete scalar product introduced by the mass lumping.
Let us denote by $\langle a,b \rangle \approx \int a\cdot b \dd x$ this discrete scalar product. Note that when no mass lumping is used we have an equality.

Then we define the projection $P_1$ by the following relation $\langle P_1(d),a\rangle = \int d\cdot a \dd x$ for all $a$ in $V_1$.

It immediately follows that, in particular for $a\in V_1$ and $b\in V_2$, we have
\begin{equation}
	\int b\cdot P_1(a\times P_1(b)) = \langle P_1(b), P_1(a\times P_1(b)) \rangle = \int P_1(b)\cdot a\times P_1(b) = 0,
\end{equation}
as $P_1(b)$ appears twice in the mixed product. This relation will be used to prove helicity conservation.
We define the orthogonality preserving cross product by 
\begin{equation}\label{eq:crossp0}
	\int a\cdot P_1(b\times P_1(c)).
\end{equation}

\subsection{Discretization} \label{eq:disc}

In this work, we employed a 3+1-splitting semi-implicit discretization: a 3-splitting for ideal MHD, combined with a symmetric Strang splitting of the resistivity terms.
Next we present a simplified sketch of the semi-implicit operator-splitting time integration that we used to run the test-cases. After writing the \lq\lq $\alpha$ \emph{splitting-step}\rq\rq\, using the following notation
\begin{equation} \label{eq:PDE1}
 \partial_t Q =  \mathcal{L}^\alpha ( Q),
\end{equation}
 where $\alpha$ labels one between the resistive ($\eta$), convection and diffusion ($v$), acoustic ($p$) or Alfv\'enic ($b$) steps. 
In the following, the numerical methods to solve the four different sub-systems are described in separate and independent sections. 

We would like to anticipate the fact that, in this work, \textbf{all implicit} \textbf{time-discretizations} required the (recursive) solution of \textbf{linear, symmetric and positive definite systems} that are solved by means of the very efficient \textbf{matrix-free} \textbf{conjugate-gradient method}, even without using any kind of preconditioning. Symmetry and positivity should be considered an important result of the proposed compatible discretization, combined with the adopted splitting-step formulation of the governing equations.
In particular, the solution of the original Alfv\'enic step is reduced to a decoupled algebraic system for only the three spatial components of the velocity, while the acoustic step is reduced to a decoupled algebraic system for the (scalar) pressure. For some tests, we also report the time evolution of the average number of iterations per time-step. Independently on the computational costs per time-step, the implicit schemes are shown to require less stabilization (i.e. less diffusion) with respect to explicit counterpart. The design and application of proper preconditioners could definitely be considered in the future works to further improve the computational efficiency.

\subsubsection{Convection-diffusion} \label{sec:convdiff}
$$
 \partial_t Q =  \mathcal{L}^v ( Q)
$$
The convection-diffusion step is dis\-cre\-tized with an explicit upwind Finite Volume scheme, where all quantities are defined at the barycenters of a \emph{dual}  Cartesian mesh, see Fig. \ref{fig:dualgrid}. In this way, the FE degrees of freedom of quantities discretized in $V_0$ may be used also as the dual FV degrees of freedom without introducing other projections or averages, as well as eventual reconstruction procedure at the end of the Finite-Volume step. We make clear that this is only one of many possible options, and it is the choice of the authors. A possible alternative is to evaluate the exact cell-average from the Finite-Element piecewise polynomial solution.
In particular the equation system (\ref{eq:FVini}-\ref{eq:FVend}) can be cast in the conservative form
\begin{equation}
\partial_t \Q + \nabla \cdot \left( \Fv + \Fdv \right) = 0,
\end{equation}
where the array of conservative variables and the tensor of conservative fluxes are defined in (\ref{sec:model}).

A standard FV scheme is derived after integrating the conservative equations in every \emph{dual} cell $\cell_{i^*,j^*,k^*}:=[x_{i},x_{i+1}]\times[y_{j},y_{j+1}]\times[z_{k},z_{k+1}]$, resulting in the following discrete equations for $\Q \in \tilde V_3$
\begin{equation} 
\label{eq:expfv3d} 
   {\Q}_{\id,\jd,\kd}^{n+1} = \Q_{\id,\jd,\kd}^n   -  (\matdiv\vecF)_{\id,\jd,\kd}, \qquad
(\vecF_{\tilde f})_{\id,\jd,\kd}  = \begin{pmatrix}
  (\mathbf{f}_{\tilde f,x})_{\id-\half,\jd,\kd} \\ (\mathbf{g}_{\tilde f,y})_{\id,\jd-\half,\kd} \\(\mathbf{h}_{\tilde f,z})_{\id,\jd,\kd-\half} 
\end{pmatrix} ~
\end{equation}
where we used the notation introduced in (\ref{eq:V2dof}) and (\ref{eq:divdof}), for the degrees of freedom of $\tilde{V_2}$ variables, used for the numerical fluxes.
It's important to note that the solution $\Q^{n+1}$ in eq.(\ref{eq:expfv3d}) is not the real final solution of the global algorithm, but it's just a temporary solution of the current splitting-step. 

Since the convective system is only \emph{weakly hyperbolic}, i.e. there is no complete set of linearly independent eigenvectors,  for simplicity we used the \emph{Rusanov} (or \emph{local Lax-Friedrichs}) numerical fluxes
\begin{multline}
 \mathbf{f}_{\id+\half,\jd,\kd} :=   \half \left( \mathbf{f}_v( \w_{\id+\half,\jd,\kd}^-) + \mathbf{f}_v( \w_{\id+\half,\jd,\kd}^+) \right)+   \\
	 -  \half s_{\alpha}^x \left(\w_{\id+\half,\jd,\kd}^+  -  \w_{\id+\half,\jd,\kd}^- \right)   -   \Aveyz{ \mathbf{f}_d(\mathbf{Q}_h, \nabla \mathbf{Q}_h)}{\id+\half}{\jd}{\kd}{} 
		\label{eq:Rusanovx} %
\end{multline}
\begin{multline}
 \mathbf{g}_{\id,\jd+\half,\kd} :=   \half \left( \mathbf{g}_v( \w_{\id,\jd+\half,\kd}^-) + \mathbf{g}_v( \w_{\id,\jd+\half,\kd}^+) \right)+  \\
			  -  \half s_{\alpha}^y \left(\w_{\id,\jd+\half,\kd}^+  -  \w_{\id,\jd+\half,\kd}^- \right)   -   \Aveyz{\mathbf{g}_d(\mathbf{Q}_h, \nabla \mathbf{Q}_h)}{\id}{\jd+\half}{\kd}{}  
		 \label{eq:Rusanovy} %
\end{multline}
\begin{multline}
 \mathbf{h}_{\id,\jd,\kd+\half} :=   \half \left( \mathbf{h}_v( \w_{\id,\jd,\kd+\half}^-) + \mathbf{h}_v( \w_{\id,\jd,\kd+\half}^+) \right) +  \\
			 -  \half s_{\alpha}^z \left(\w_{\id,\jd,\kd+\half}^+  -  \w_{\id,\jd,\kd+\half}^- \right)  -   \Aveyz{\mathbf{h}_d(\mathbf{Q}_h, \nabla \mathbf{Q}_h)}{\id}{\jd}{\kd+\half}{} 
			\label{eq:Rusanovz} %
\end{multline}
Following \cite{ToroVazquez}, a valid alternative numerical-flux for the \emph{convection step} (\ref{eq:FVini}-\ref{eq:FVend}) is an \emph{upwind} numerical flux of the type
\begin{multline}\label{eq:Upwindx}
 \mathbf{f}_{\id+\half,\jd,\kd} :=   \half \left( \mathbf{f}_v( \w_{\id+\half,\jd,\kd}^-) + \mathbf{f}_v( \w_{\id+\half,\jd,\kd}^+) \right)+   \\
	 -  \half  \left( \omega^+ \mathbf{f}_v( \w_{\id+\half,\jd,\kd}^+)  -  \omega^- \mathbf{f}_v( \w_{\id+\half,\jd,\kd}^-)  \right)     -   \Aveyz{ \mathbf{f}_d(\mathbf{Q}_h, \nabla \mathbf{Q}_h)}{\id+\half}{\jd}{\kd}{}   
\end{multline}
where the positive weight $\omega$ is  defined by
$$
\omega =  \omega(u_{\id+\half,\jd,\kd}) = \frac{u_{\id+\half,\jd,\kd}}{\sqrt{\epsilon_\omega + (u_{\id+\half,\jd,\kd})^2}},
$$
and the corresponding formulas for the y- and z- fluxes, $\mathbf{g}_{\id,\jd+\half,\kd}$and $\mathbf{h}_{\id,\jd,\kd+\half}$.
Here $\epsilon_\omega$ is a sufficiently small non-negative number. In this work, we just tested this numerical flux against a high-Mach shock dominated problem, setting $\epsilon_\omega=10^{-14}$.

Here $\w^{\pm}$ must be interpreted as the space-time reconstructed value of the conserved variables at a given (dual) cell-interface from the right ($+$) and from the left ($-$) at the middle time $t=t^n + \Delta t/2$, i.e.
 the numerical fluxes are evaluated according to 
$$\w_{\id+\half}^{-} = \w_{\id}(x_{\id+\half},t^n+\Delta t/2) \qquad \text{and} \qquad \w_{\id+\half}^{+} = \w_{\id+1}(x_{\id+\half},t^n+\Delta t/2).$$
In the Rusanov flux, the stabilization terms are defined as jump-penalizations proportional to the maximal directional eigenvalue $s_{\alpha}^i = \max(\{ \lambda^{\alpha,i}_{1-8}\})$, $i=x,y,z$. In this case, only the convective terms are solved with an explicit-in-time scheme, and then we used $s_v =s  ( \lambda^{v} )  $, referring to (\ref{eq:lambdac}). For those tests where a fully explicit scheme is mentioned, we just set $s_{\text{max}} =s_{\text{max}} ( \lambda^{\MHD} )  $.

\paragraph{A polynomial reconstruction} 
Several choices can be done to build an accurate polynomial reconstruction $\w$ to evaluate the numerical fluxes (\ref{eq:Rusanovx}-\ref{eq:Rusanovz}). Accuracy and robustness are key features for a good reconstruction operator. In this work, a natural choice is the original FE polynomials, i.e.
\begin{equation}\label{eq:FEECrec}
\w^{\texttt{FEEC}}(x,t) = \left( \rho_h, \rho\u_h, E_h, \B_h \right)  \in     V^0\times  V^1 \times V^0 \times V^2,
\end{equation}
which preserves the order of accuracy of the FE splitting steps, and, moreover, these polynomials are already available in the algorithm flow. This choice is found to be efficient and accurate for smooth solutions like the Alfv\'en wave test, see Sec. \ref{sec:results}, but not sufficiently robust against strong shocks.
For this reason, we also introduce  
a  second-order TVD MUSCL-Hancock Finite-Volume reconstruction, see \cite{toro-book}. This can be efficiently done by introducing a space-time piece-wise linear reconstruction that, in the one dimensional case, reads as 
\begin{align} \label{eq:MUSCLrec}
	\w^{\texttt{FV}}_{\id}(x,t)  &= 	\Q_{\id}^n + \nabla \w_{\id}^n \cdot (x-x_{\id})+ \partial_t \w_{\id}^n   (t-t^n)\\
\partial_t \w_{\id }^n   &= - \frac{\mathbf{f}(\w_{\id }(x_{\id+\half},t^n))-\mathbf{f}(\w_{\id }(x_{\id-\half},t^n))}{\Delta x_{\id}} \\ \label{eq:minmod}
\nabla \w_{\id }^n   &=   \frac{\Delta \Q_{\id}^n}{\Delta x_\id}  := \text{\textbf{minmod}}\left(   \frac{\Q_{\id}^n - \Q_{\id-1}^n}{\Delta x_{\id-\half}},  \frac{ \Q_{\id+1}^n -  \Q_{\id}^n}{\Delta x_{\id+\half}} \right).
\end{align}

Using a \emph{non-limited} reconstruction $\w$, e.g. using centered averages instead of the \texttt{minmod} as slope estimator (\ref{eq:minmod}), would cause the scheme to lose the TVD property and then being similar to $\w^{\texttt{FEEC}}$ in terms of robustness and accuracy. 

\medskip
Finally, we want to anticipate that the FE momentum equation in the Alfv\'enic step (see Section \ref{sec:split-alf}) has been formulated in non-conservative form, see also eq.(\ref{eq:splitalf-2}). To restore the conservation properties, this FV scheme  will be also used to update the dual-barycentric variables in the following FE splitting steps. This will guarantee the final algorithm to be \emph{exactly conservative for the total energy density and momentum}.

\subsubsection{The acoustic part}\label{sec:acoustic}
$$
 \partial_t Q =  \mathcal{L}^p ( Q)
$$
Let us first perform the semi-discretization in time in the Finite Element formulation.
First as $ p_h\in V_0$, we have  $ \nabla p_h\in V_1$. So as $ \mathbf{m}_h\in V_1$, an implicit Euler discretization in time of equation \eqref{eq:splitac-2} yields 
\begin{equation}\label{eq:FEm}
 \mathbf{m}_h^{n+1}= \mathbf{m}_h^n - \Delta t \nabla p_h^{n+\theta_p}.
\end{equation}
On the other hand, a weak form of \eqref{eq:splitac-4} reads: \textrm{find} $p\in V_0$ \textrm{such that }
	   \begin{equation} \label{eq:p}
		 \fract{}{t} \int \left( \frac{p_h}{\gamma-1} + \frac {1}{2} \mathbf{u}_h\cdot \mathbf{m}_h \right) q_h - \int  \frac{\gamma}{\gamma-1}\frac{p_h}{\rho}\mathbf{m}_h \cdot \nabla q_h = 0 \quad\forall q_h \in V_0
	   \end{equation}

	   Discretizing in time, linearizing, plugging in the expression for $\mathbf{m}_h^{n+1}$ and introducing an iteration index $s_p$, noting that $\rho$ remains constant in this split step, yields
	   \begin{align} \label{eq:psystem}
		 &  \int \frac{p^{n+1,s_p+1}_h}{\gamma-1} q_h \dd \mathbf{x} + \theta_p \Delta t^2\int \frac{\gamma}{\gamma-1}\frac{p^{n+\theta_p,s_p}_h}{\rho^{n+1}} \nabla p^{n+\theta_p,s_p+1}_h\cdot\nabla q_h \dd\mathbf{x} 
		 = F^{n,s_p} \quad\forall q_h \in V_0 
		\end{align}
which is a \emph{symmetric positive definite linear} system \emph{at each Picard iteration}, assuming the positivity of the computed density $\rho^{n+1}$ and the predictor pressure $p^{n+\theta_p,s_p}$.
		The known-term in the right-hand-side is
		\begin{equation} \label{eq:prhs}
	\begin{aligned}
		 F^{n,s_p}  = & \int  \frac{ p^{n}_h}{\gamma-1} q_h \dd \mathbf{x} -  \int \left(  \frac {1}{2} \mathbf{u}_h^{n+1,s_p}\cdot \mathbf{m}_h^{n+1,s_p} -  \frac {1}{2} \mathbf{u}_h^{n}\cdot \mathbf{m}_h^{n} \right) q_h +  \\
	& +   \Delta t \int \frac{\gamma}{\gamma-1} \frac{p^{n+\theta_p,s_p}_h}{\rho^{n+1}} \mathbf{m}^n_h \cdot \nabla q_h . 
		 \end{aligned}
		\end{equation} 
Finally, we can express the non-linear solver for the acoustic step as a recursive implicit solve of linear system for  $p \in V_0$, $\mathbf{m}_f\in V_1$:
\begin{align}
	&(\matmass_0 + \theta_p^2 \Delta t^2 \matgrad^T \matmass_1 \matgrad)  \vecp^{n+1,s_p+1} = \vechn^{n,s_p} \label{eq:psyst-1} \\ 
	&\vecm^{n+1,s_p+1} = \vecm^m -\Delta t \matgrad\vecp^{n+\theta_p,s_p+1} \label{eq:psyst-2}
  \end{align}
	where $\matmass_0$ is the  mass matrix of $V_0$
	$$
	 (\matmass_0)_{mn}  = \int  q_m q_n \dd\mathbf{x}, \qquad q_m, q_n \in V_0,
	$$ the $\matmass_1$ is the positive weighted mass matrix of $V_1$
	$$
	 (\matmass_1)_{mn}  = \int \gamma \frac{p^{n+\theta_p,s_p}_h}{\rho^{n+1}} \C_m \C_n \dd\mathbf{x}, \qquad \C_m,\C_n \in V_1,
	$$
	while the  right-hand side
	\begin{equation}
\begin{aligned}
	 \vechn^{n,s_p} = &\matmass_0 \left(   \vecp^{n} - (\gamma-1)\left(\frac 12 \overline{\vecu \vecm}^{n+1,s_p}  
  + \frac 12 \overline{\vecu \vecm}^{n} \right)\right) + \Delta t \matgrad^\top \matmass_1 \vecm^{n}   +\\
  &-  (1-\theta_p )\Delta t^2 \matgrad^T \matmass_1 \matgrad  \vecp^{n}.
\end{aligned} 
	\end{equation}
\paragraph{Energy and momentum conservation} 
	The acoustic step is completed once the (dual) barycentric values of the total energy density $\rho E$ have been also updated according to eq.(\ref{eq:splitac-4}), with the aforementioned \emph{conservative} FV scheme (\ref{eq:expfv3d}) but using the analytical fluxes evaluated directly via interpolation of the FE variables  $p^{n+1}\in V_0$, $\u_e^{n+1}\in V_1$. 
	This guarantee exact conservation and accuracy also for the total-energy density.
	 Moreover, also the linear momentum is exactly conserved because the discrete formula (\ref{eq:FEm}) is conservative.
	Extension of the presented method to a more general nonlinear equation of state is clearly possible, following the work of \cite{DumbserCasulli2016} for compressible Euler, see also \cite{SIMHD,Fambri20} for MHD.
	\paragraph{Numerical stabilization and convergence}
	In this work we consider the following symmetric positive semi-definite term 
	   \begin{align} \label{eq:pstab}
		 &  \frac{\gamma}{\gamma-1}\int s_p \epsilon \nabla p^{n+\theta_p,s_p+1}_h\cdot\nabla q_h \dd\mathbf{x}  
		\end{align}
	which can be interpreted as a Rusanov inspired numerical stabilization terms and has to be added to the pressure equation (\ref{eq:p}). Here, $s_p$ is chosen to be the maximum eigenvalue of the acoustic system in absolute value, while $\epsilon$ is a non-negative scalar-valued function. Indeed, one can verify that the corresponding \emph{modified} discrete pressure system reads like
	   \begin{align}\label{eq:p2}
		 &  \frac{1}{\gamma-1}\int {p_h }q_h \dd \mathbf{x} + \Delta t^2 \theta_p^2 \int \tilde{h} \nabla p_h\cdot\nabla q_h \dd\mathbf{x} 
		 = \tilde{F} \quad\forall q_h \in V_0 
		\end{align}
		where, assuming for simplicity $\theta_p>0$, $\tilde{h}$ is an effective enthalpy defined as
\begin{equation}\label{eq:c}
\tilde{h} = h_h + \frac{s_p  \epsilon}{\theta_p \Delta t}  \geq h_h,
\end{equation}
where $h_h:= (p_h/\rho) \gamma/(\gamma-1)$. First, observe that the coefficient $(\gamma-1)\theta_p \Delta t \tilde{h}_h$ plays the role of an effective viscosity. 
 In this sense, the stabilization term (\ref{eq:pstab}) effectively introduces an additional artificial viscosity, proportional to 
$ {s_p  \epsilon} \geq 0.$
It seems to be reasonable to define $\epsilon = C O(\Delta x^N)$, where $N$  is the order of accuracy of the method, $C$ is a non-negative real constant.  Then, after introducing the $\lambda_\alpha$-based Courant number $\CFL_{\alpha} = s_\alpha {\Delta t}/{\Delta x}$,
the modified enthalpy can be finally written as 
\begin{equation}\label{eq:c2}
\tilde{h} = h_h + |C| \frac{s_p  s_\alpha}{\theta_p \CFL_\alpha} O(\Delta x^{N-1})  \geq h_h,
\end{equation}
and the solution is expected to converge to the solution of the original  discrete problem (\ref{eq:p}) with the order of the scheme.

While this choice may become useful to stabilize the scheme for higher Mach regimes, one should also consider the remote possibility of stabilizing the scheme against negative solutions. As it was mentioned above, a simple matrix-free conjugate gradient method may efficiently convergence to the unique solution of system (\ref{eq:psystem}) (or \ref{eq:psyst-1}), only after assuming the positivity of $h_h$.
Now, the stability and convergence is guaranteed by only requiring the \emph{positivity of the effective enthalpy},  i.e. wherever $h_h<0$ holds,  the positivity of the modified enthalpy (\ref{eq:c}) needs a lower bound for $\epsilon$, i.e. from(\ref{eq:c})
 \begin{align} \epsilon \geq -  \Delta t \theta_p \frac{ h_h}{ s_p}  
= - \Delta x  \CFL_p \theta_p \frac{ h_h}{ s_p^2} \qquad  \Longleftrightarrow \qquad  \tilde{h} \geq 0  . 
\end{align}
In particular, for negative enthalpies  a sufficient lower bound of the artificial viscosity is obtained by setting
\begin{align}
 \epsilon \left(h_h<0\right) > \Delta x (\gamma-1) \theta_p  \CFL_p, \qquad \Longrightarrow  \qquad \tilde{h}_h \geq 0 
\end{align}
which comes from (\ref{eq:c}), after using the definition of $h_h$ and the inequality  $$| \gamma p_h/\rho | / s_p \leq 1.$$

These estimates are important \emph{only to be aware} of the amount of \emph{local} artificial diffusion we are introducing by requiring the positivity of the effective enthalpy. In particular,  observe that, for negative solutions, the artificial viscosity scales \emph{locally} with the mesh resolution $\Delta x$. As a direct consequence, the presented method is robust even against negative solutions and, only when this happens, the method reduced to be \emph{locally} first order accurate in space.
However, to build an \emph{unconditionally stable solver}, in the practice it is easier to define the physical enthalpy and the error estimator as 
\begin{align}\label{eq:hdef}
h_h \longrightarrow \max\left( h_h, 0 \right), \qquad  \epsilon := c_{\tilde{h}} \frac{\Delta x}{2}
\end{align}
This special choice  yields to an effective enthalpy that is provably non-negative independently of the positivity of the numerical solution. This is a \emph{sufficient condition to guarantee the convergence to the unique solution} of  system (\ref{eq:p2}) by means of a very efficient matrix-free conjugate gradient method \cite{cgmethod}.
Observe that the modified pressure equation is still \emph{energy conserving}, because the additional viscous term is conservative.

 The matrix of the acoustic algebraic system is 
$$
\Hmat := \matmass_0 + \theta_p^2 \Delta t^2 \matgrad^T \matmass_1 \matgrad.
$$
Observe that the first term is the $V_0$ mass matrix, which is purely diagonal and strictly positive definite, while in the second term $ \matgrad^T \matmass_1 \matgrad$ is the induced discrete Laplacian, which is a tridiagonal matrix (1d), pentadiagonal (2d) or eventually eptadiagonal (3d), and it is also \emph{strictly diagonally dominant} with positive, i.e.
$$
\Hmat_{ii}>0, \qquad \Hmat_{ij}<0 \quad\text{for $i\neq j$}, \qquad \Hmat_{ii} > \sum_{j\neq i} |\Hmat_{ij}|
$$
 This is the result of the  \emph{mass lumping} approximation (\ref{eq:lamped}), yielding  
\begin{multline}
 \mathbb{L}_0 := \matgrad^T \matmass_1 \matgrad   =   \frac{1}{\Delta x^2}  \begin{pmatrix} 
-1 & 0 &  \cdots& 0 & +1 \\
+1 & -1 & 0 & \cdots & 0 \\
0 & +1 & -1  & \cdots &  \vdots \\
 \vdots  &   & \ddots   &  \ddots  &  0 \\
0 & \cdots & 0 & +1 & -1  
\end{pmatrix}   \matmass_1
 \begin{pmatrix} 
-1 & +1 & 0 & \cdots & 0 \\
0 & -1 & +1 &  			& \vdots\\
 0 &     & \ddots  &  \ddots & 0 \\
 0 &  \cdots  & 0  &   -1 & +1 \\
+1  &  0 & \cdots &   & -1  
\end{pmatrix} \\
 = \frac{1}{\Delta x^2}\begin{pmatrix}
 m^1_1 + m^1_{N^1_{tot}}   	& - m^1_1 		& 0 				& \cdots 						&	 				&   -  m^1_{N^1_{tot}} \\
- m^1_1   & m^1_1 + m^1_2 & -m^1_2 			& 0 								&	 \cdots	& 0 \\
0   		& - m^1_1     & m^1_1 + m^1_2 & -m^1_2 							& 0 			& \vdots \\
 \vdots	& \ddots 		&\ddots  		& \ddots 						& \ddots  & 0\\
-m^1_{N^1_{tot}} 	&  \cdots 		&     			& 	0				&  -  m^1_{N^1_{tot}-1} &	 m^1_{N^1_{tot}-1} +m^1_{N^1_{tot}} 
\end{pmatrix} \nonumber
\end{multline}
for the one dimensional case, where we defined $\matmass_1 = diag(\{m^1_i\})$, $i=1,\ldots,N^1_{tot}$, and $N^0_{tot}$, $N^1_{tot}$ being the total number of degrees of freedom of $V^0$ and $V^1$, respectively.
Observe that, thanks to the use of the compatible FE, these properties remain valid also in the higher order case, since the discrete gradient is the same and the $V^1$ mass matrix is still diagonal by using the corresponding higher-order accurate histopolant or interpolant polynomial basis function. Moreover, whenever Dirichlet pressure boundary conditions are applied, the induced discrete Laplacian $\mathbb{L}_0$ becomes non-singular and invertible, with $\mathbb{L}_0^{-1}>0$.
	
Note also that, by choosing $\theta_p=1/2$, the off diagonal terms of $\mathbb{H}$ are scaled by a factor $1/4$ and the dominance of the diagonal is further enforced, with a consequent faster convergence of the solver.

\paragraph{Remark about pressure positivity}

Given the assumptions above to guarantee the positivity of the numerical enthalpy, by using artificial stabilization terms (\ref{eq:hdef}), the resulting algebraic system (\ref{eq:psyst-1})	forms an \emph{M-matrix} and \emph{the unique solution $p^{n+1}$ is always positive provided that the right-hand side is also positive}.
Numerical strategies to ensure positivity of the right-hand side (\ref{eq:prhs}) are  not yet investigated by the authors. This may eventually require other constraints, see e.g. Greenspan and Casulli \cite{GreenspanCasulli}, and the use of upwinding, Lagrangian or semi-Lagrangian techniques.

\subsubsection{The Alfv\'enic part} \label{sec:split-alf}
$$
 \partial_t Q =  \mathcal{L}^b ( Q)
$$
We look for $\mathbf{B}_h\in V_2$ and $ \mathbf{m}_h\in V_1$, so that $\nabla\cdot\mathbf{B}_h\in V_3$ is defined \emph{strongly}. However,  $\mathbf{u}_h \times \mathbf{B}_h$ is not in $V_1$. So we will project it with the orthogonal projection in $V_1$ that we denote by $P_1$.
The Galerkin approximation of \eqref{eq:splitalf-2}-\eqref{eq:splitalf-3} then reads
\begin{align}
		&  \fract{}{t} \int\rho\mathbf{u}_h\cdot \mathbf{C}_h +  \int \nabla\times  P_1( \mathbf{C}_h \times P_1(\mathbf{B}_h)) \cdot  \mathbf{B}_h   = 0~ \forall{\C_h\in V_1}, \label{eq:splitalfdisc-2}  \\
		 &  \fracp{ \B_h}{t}  +   \nabla\times P_1 \left( -  \u_h  \times P_1 (\mathbf{B}_h) \right)   =0.  \label{eq:splitalfdisc-3}
	   \end{align}

The semi-discretization in time yields respectively
	   \begin{align}\label{eq:semidisc-strongB}
		 \B_h^{n+1}= \B_h^n +\Delta t \nabla\times P_1 ( \u_h^{n+\theta_b}  \times P_1(\mathbf{B}^{n+\theta_b}_h) )&\\
		\label{eq:semidisc-weakm0}
		  \int \rho_h \u_h^{n+1}\cdot \mathbf{w}_h +  \Delta t \int \nabla\times  P_1(\mathbf{w}_h \times  P_1(\mathbf{B}^{n+\theta_b}_h)) \cdot \B^{n+\theta_b}_h   = \int \rho_h \u_h^{n}\cdot \mathbf{w}_h& \\
			\nonumber ~\forall \mathbf{w}_h \in V_1&
	   \end{align}
Let us first observe that $\nabla \cdot \mathbf{B}_h$ stays zero  if it is zero at the initial time.

Plugging \eqref{eq:semidisc-strongB} into \eqref{eq:semidisc-weakm0} and introducing nonlinear iterations yields
	\begin{multline}\label{eq:semidisc-weakM}
		 \int \rho_h\u_h^{n+1,s_b+1}\cdot \mathbf{w}_h +  \theta_b \Delta t^2 \int \nabla\times  P_1(\mathbf{w}_h \times  P_1(\mathbf{B}^{n+\theta_b}_h)) \cdot \nabla\times P_1 ( \u_h^{n+\theta_b,s_b+1} \times P_1(\mathbf{B}^{n+\theta_b}_h) ) \\ = -\Delta t \int \nabla\times  P_1(\mathbf{w}_h \times  P_1(\mathbf{B}^{n+\theta_b}_h)) \cdot \B^{n}_h  +\int \rho_h \u_h^{n}\cdot \mathbf{w}_h ~\forall \mathbf{w}_h \in V_1
	\end{multline}  
We notice that the left-hand-side is a \emph{symmetric positive definite bilinear form} \emph{at each nonlinear iteration}.
This can be solved for $ \u_h^{n+1}$ by Picard iterations, since the nonlinear system for $\vecv$ is decoupled from  $ \vecb$. Finally, $ \mathbf{B}^{n+1}$ is given by \eqref{eq:semidisc-strongB}.

Again, we can express the non-linear solver for the Alfv\'enic step as a recursive implicit iterative solve of  nonlinear system for $ \vecv \in V_1$:
  \begin{equation}\label{eq:disc-weakm0}
	(\matmass_1^\rho + \theta_b^2\Delta t^2 \matproj_{B^s_b}^T \matcurl^T\matmass_2 \matcurl \matproj_{B^s_b} )  \vecv^{n+1,s_b+1} = \veche^{n,s_b} 
  \end{equation}
	where $ \matproj_{B^s_b}\vecv^{n+1,s_b+1}$ is associated to $   P_1(\v_h^{n+1,s_b+1} \times  \B^{n+1,s_b}_h)$,
	  $\matmass_2$ is the  mass matrix of $V_2$
	$$
	 (\matmass_2)_{mn}  = \int  \v_m \cdot \v_n \dd\mathbf{x}, \qquad \v_m, \v_n \in V_2,
	$$   $\matmass_1^\rho$ is the positive weighted mass matrix of $V_1$
	$$
	 (\matmass_1^\rho)_{mn}  = \int  \rho^{n+1} \w_m \cdot\w_n \dd\mathbf{x}, \qquad \w_m,\w_n \in V_1,
	$$
	while the  right-hand side
	\begin{equation}
\begin{aligned}
	 \veche^{n,s_b} = &\matmass_1^\rho \vecv^{n} -  \theta_b (1-\theta_b) \Delta t^2 \matproj_{B^s_b}^T \matcurl^T\matmass_2 \matcurl \matproj_{B^s_b}   \vecv^{n} .
\end{aligned} 
	\end{equation}
	
	\paragraph{Energy conservation}
	

Note that this formulation is also \emph{energy conserving}   at the semi-discrete level with respect to the following discrete energy
$$
E_b = E_b(\rho_h,\m_h,\B_h)= \int \half \frac{\m_h^2}{\rho} + \int \half\B^2.
$$
This can be easily checked, by evaluating the time derivative of the total energy
\begin{equation}\label{eq:Eb}
\frac{d}{dt} E_b  =  - \half \int |\u_h|^2 \frac{\partial  \rho }{\partial t} +  \int \u_h \cdot \frac{\partial \m_h }{\partial t}  + \int \B_h \cdot\frac{\partial \B_h}{\partial t}  = 0
\end{equation}
The equality (\ref{eq:Eb}) is easily proved by considering that:  
using  eq. (\ref{eq:splitalfdisc-2}) with   $\mathbf{C}_h=\mathbf{u}_h\in V_1,$   and eq. (\ref{eq:splitalfdisc-3})  after multiplication by  $\B_h$ and integration over the full domain. A Crank-Nicolson time-discretization $\theta_b=\frac{1}{2}$ would ensure energy conservation also at the fully discrete level.
\paragraph{Magnetic helicity conservation}
 
Thanks to FEEC and the orthogonality preserving cross product (\ref{eq:crossp0}), we can show that magnetic helicity  is conserved, i.e.
$$
\frac{d H}{d t} :=\frac{d  }{d t}  \int \mathbf{A} \cdot \mathbf{B} = 0
$$
 Indeed, after defining the vector potential to be $\mathbf{A}_h\in V_1$ such that $\nabla \times  \mathbf{A}_h$ is defined \emph{strongly} in $V_2$, and we may define
\begin{equation*}
\begin{aligned}
&\nabla \times  \mathbf{A}_h = \mathbf{B}_h, \quad \mathbf{A}_h \in V_1,\mathbf{B}_h \in V_2 \\
&\int  q_h  \nabla \cdot \mathbf{A}_h = 0, \quad \forall q_h \in V_0
\end{aligned}
\end{equation*}  
 one can first show that
\begin{equation}
\begin{aligned}
\frac{d H}{dt}  = 2 \int \mathbf{A}_h \cdot \frac{\partial \mathbf{B}_h}{\partial t}  .
\end{aligned}
\end{equation} 
Then, by multiplying (\ref{eq:splitalfdisc-3}) by $\mathbf{A}_h$ and integrating over the full domain one obtains
\begin{equation}
\begin{aligned}
\frac{d H}{dt}  & = 2 \int \mathbf{A}_h \cdot \nabla\times P_1 \left( \v_h \times P_1(\mathbf{B}_h) \right) \\
&   =  \int \nabla\times  \mathbf{A}_h \cdot P_1 \left( \v_h  \times P_1(\mathbf{B}_h) \right) \\
	& =  \int P_1(\mathbf{B}_h)\cdot  \left( \v_h \times P_1(\mathbf{B}_h) \right) \equiv 0   
\end{aligned}
\end{equation}
One has to notice that, in order to conserve magnetic energy or helicity at the fully discrete level, one has to choose a Crank-Nicolson time discretization ($\theta_b=\frac{1}{2}$), and to solve the non-linear system (\ref{eq:semidisc-weakM}) 
until convergence.
 
	The Alfv\'enic step is completed once both the (dual) barycentric values of the total energy density $\rho E$ and linear momentum $\m$ have been also updated according to the conservative tensor fluxes (\ref{eq:MHD-Fb}), with the aforementioned \emph{conservative} FV scheme (\ref{eq:expfv3d})  applied to the Alfv\'enic sub-system
	\begin{equation}
\partial_t \Q + \nabla \cdot  \Fb  = 0.
\end{equation} 
More specifically, the update of barycentric 	values follows the conservative update formula  for $\Q \in \tilde V_3$
\begin{equation*} 
   {\Q}_{\id,\jd,\kd}^{n+1} = \Q_{\id,\jd,\kd}^n   -  (\matdiv\vecF_b)_{\id,\jd,\kd},  
\end{equation*}
where the numerical fluxes are evaluated with the new FE variables  $\u_e^{n+\theta_v}\in V_1$ and  $\B_f^{n+\theta_b}\in V_2$, i.e.  $\vecF_b=\vecF_b(\rho_h, \rho \u_e^{n+\theta_v}, p_h,  \B_f^{n+\theta_b})$,
 using the analytical fluxes evaluated directly via interpolation of the FE variables  $\u_e^{n+\theta_v}\in V_1$ and  $\B_f^{n+\theta_b}\in V_2$. 
	Also in this case, this strategy guarantees exact conservation and accuracy also for the total-energy density.
	 Moreover, also the linear momentum is exactly conserved because the discrete formula (\ref{eq:FEm}) is itself conservative.

\subsubsection{The resistive part}\label{sec:res}
$$
 \partial_t Q =  \mathcal{L}^\eta ( Q)
$$
As it was mentioned above, a natural compatible discretization of resistive terms refer to a weak formulation of the Faraday equation. Indeed, we look for $\tilde \B_h 
\in \tilde V_1$, so that $\nabla \times\tilde \B_h\in \tilde V_2$ is defined strongly. To do that, the so-called \emph{Hodge operator}, mapping the primal solution space $V_2$ into the dual $\tilde V_1$, is defined, i.e. 
$$
\tilde{\mathbb{H}}_1: V_2 \longrightarrow \tilde V_1,
$$
with $\tilde \B  = \tilde{\mathbb{H}}_1 \B$.
The major requirement for the definition of the Hodge operator  is the preservation of the order of accuracy, but many different choices can be made, e.g. between interpolation, Galerkin projection or others. Note that, in the low order case, a Hodge transformation between $V_2$ and $V_1$ may introduce a very inconvenient Lax-Friedrich-type dissipation, e.g. after cell-averaging node-to-cell points or vice-versa. For this reason, in this work we defined the \emph{dual} $\tilde \B$ as a  vector variable in $\tilde V_2$, which is a discrete $H$-curl space built on the \emph{dual grid}, see  Fig. \ref{fig:dualgrid},  and it is \emph{dual-to} $V_1$ through the Hodge $\mathbb{H}_2$.
 The convenience of this strategy becomes even more evident after observing that, in the low-order case, the spatial representation of the primal $V_2$ d.o.f. coincide with the dual $\tilde V_1$. Indeed, the face-barycenters of the main grid coincide with the edge-barycenters of the dual, see Fig. \ref{fig:dualgrid}. This means that the Hodge is an \emph{identity} operator in the low order case, while it can be regarded as an approximation of the identity in the general high-order case.

 The Galerkin approximation of \eqref{eq:splitres-3} then reads
	   \begin{align}
		   \fract{}{t} \int \tilde \B_h\cdot \mathbf{C}_h  +    \int \eta \nabla\times\tilde \B_h\cdot \nabla \times\mathbf{C}_h =0 ~\forall \mathbf{C}_h \in \tilde V_1\label{eq:weakB0res}
	   \end{align}
 The semi-discretization in time yields  
	   \begin{align}\label{eq:semidisc-weakB0res}
		 \int  \tilde\B _h^{n+1}\cdot \tilde\C_h  + \Delta t\int\eta\nabla \times  \tilde\B_h^{n+\theta_{\tilde{b}}}\cdot \nabla \times \tilde\C_h   = \int  \tilde\B _h^{n}\cdot \tilde\C _h ~\forall \tilde\C_h \in  \tilde V_1
	   \end{align}  
 We notice that the left-hand-side is a symmetric positive definite bilinear form, and it can be solved very efficiently for $\tilde\B_h^{n+1}$.

Finally,  the resistive step may be written in matrix-tensor form as an independent linear and positive definite system for only the degrees of freedom of the dual $ \vecbd \in \tilde{V}_1$:
  \begin{equation}
	(\tilde \matmass_1  + \theta_{\tilde{b}} \Delta t \matcurld^T \tilde \matmass^\eta_2 \matcurld  )\vecbd^{n+1}    =     \vecheb^{n} 
  \end{equation}
	where the known terms are
		$$
		   \vecheb^{n} = \tilde \matmass_1 \vecbd^{n}  - (1-\theta_{\tilde{b}} )   \matcurld^T \tilde \matmass^\eta_2 \matcurld \vecbd^{n},
		$$
	$\tilde\matmass_1$ and $\tilde \matmass^\eta_2$  are the  mass matrix  and the positive weighted mass-matrix of $\tilde{V}_2$, i.e.
\begin{align}
	 (\tilde \matmass_1)_{mn}  &= \int \tilde{\C}_m \cdot\tilde{\C}_n \dd\mathbf{x},      & \tilde{\C}_m, \tilde{\C}_n \in \tilde{V}_1, \\
	 (\tilde \matmass^\eta_2)_{mn}  &= \int \eta \tilde{\D}_m \cdot\tilde{\D}_n \dd\mathbf{x}, &\tilde{\D}_m, \tilde{\D}_n \in \tilde{V}_2,
	\end{align}
	correspondingly.
	  
Then, even if the computed \emph{dual} $\tilde\B_h^{n+1}$ is only \emph{weakly} divergence free, we may use it to update the \emph{primal} $\B\in V_1$ via a \emph{strong} Galerkin discretization that reads
\begin{equation}\label{eq:semidisc-primalBupdate}
	 \B _h^{n+1} =    \B _h^{n} -  \Delta t  \nabla \times  \left( P_1( \eta\nabla \times \tilde\B_h^{n+1})\right) .
\end{equation} 
Here $P_1$ is a Galerkin orthogonal projection to $V_1$.
   Then, by construction  $\nabla \cdot \B_h$ stays zero if it is zero at the initial time (\emph{strongly}). 

	\paragraph{Energy stability}
	

On the other hand, for the resistivity step we have \emph{energy stability}, i.e. at the semi-discrete level the inequality
$$
\frac{d}{dt} \tilde E_\eta  \leq 0  
$$ 
 holds,   with respect to the magnetic energy 
$
\tilde  E_\eta = \tilde E_\eta(\tilde \B_h)= \int \half |\tilde \B|^2_h.
$
Indeed, by evaluating the time derivative of the total energy
$$
\frac{d}{dt} \tilde E_\eta  = \int \tilde \B_h \cdot\frac{\partial \tilde \B_h}{\partial t}  = -  \int \eta \nabla\times\tilde \B_h\cdot \nabla \times \tilde \B_h \leq 0
$$ 
 where we substituted  eq. (\ref{eq:weakB0res}) after setting $ \mathbf{C}_h= \tilde \B_h\in \tilde{V}_1$. 

Obviously, for ideal flows ($\eta=0$) the equality holds and the semi-discrete scheme is said to be energy conserving (or preserving) with respect to $\tilde E_\eta$.  Again, a Crank-Nicolson time-discretization $\theta_b=\frac{1}{2}$ would ensure energy stability (or conservation) also at the fully discrete level.
The resistivity step is completed once the total energy is also updated by applying the aforementioned conservative FV scheme with respect to the conservative fluxes of eq. (\ref{eq:MHD-Fb}). The FV scheme (\ref{eq:expfv3d}) is used here only for the update of the total energy density $\rho E$, while the (dual-) barycentric values of the magnetic field are evaluated directly as cell-average of the FE $\tilde{V}_1$-polynomials.
\paragraph{Adding numerical stabilization}
In this work, we add to the resistivity equation the following artificial dissipative term
\begin{equation}
\nabla \times \nabla_{\tilde{\eta}} \times \tilde \B_h,  
\end{equation}
that we define b as
\begin{align}
( \nabla_{\tilde{\eta}} \times \tilde \B_h)_i &:= \sum_{  jk} \levi_{ijk}   \tilde{\eta}_{j}^i \left( \partial_j \tilde B_k\right):= \sum_{  jk} \tilde\levi_{ijk}     \left( \partial_j \tilde B_k\right),  
\end{align} 
where  $\levi_{ijk}$ is the Levi-Civita symbol.
Inspired by the theory of multi-dimensional-Riemann solvers (see e.g. \cite{balsarahlle2d,balsarahllc2d,balsarahlle3d,BalsaraMultiDRS,MUSIC1,MUSIC2}) and following \cite{SIMHD,Fambri20}, we defined the artificial resistivity to be proportional to the Rusanov viscosity and, similarly to the pressure system, to an error estimator $\epsilon=\epsilon(\Delta x)$, i.e.
$$
\tilde{\eta}_{j}  :=   s_\alpha  \epsilon(\Delta x_j).
$$
Now we want to prove the conditions for this artificial stabilization term to guarantee \emph{energy stability}.

Indeed, after multiplication by a test function $\tilde \C_h = (\tilde C_x, \tilde C_y, \tilde C_z)$ and integration on the spatial domain, one derives the following bilinear form
\begin{align}   \label{eq:resterm} 
 &\int  \tilde \C_h \cdot \nabla \times \nabla_{\tilde{\eta}} \times \tilde \B_h  = 
\int \nabla \times  \tilde  \C_h\cdot \nabla_{\tilde{\eta}} \times \tilde \B_h.  
\end{align} 
By re-writing the bilinear form as acting directly on the components of $\nabla \B$, one may notice that the corresponding coefficient matrix
is \emph{block-diagonal}. Indeed, since the diagonal terms $(\nabla \B)_{ii}$ and $(\nabla \C)_{ii}$  for $i=1,2,3$ do not enter into (\ref{eq:resterm}), we may rewrite integrand of   (\ref{eq:resterm})  as 
 \begin{align}  \label{eq:block}
 \nabla \times  \tilde  \C \cdot \nabla_{\tilde{\eta}} \times \tilde \B_h  \equiv \sum_{\substack{ i=1,2,3\\ (ijk)\in P_{(1,2,3)} }} \begin{pmatrix} \partial_j C_k & \partial_k C_j  \end{pmatrix} \cdot \begin{pmatrix}
\tilde\eta_j  & -\tilde\eta_k \\ -\tilde\eta_j & \tilde\eta_k
\end{pmatrix}  \cdot \begin{pmatrix}
\partial_j B_k \\ \partial_k B_j
\end{pmatrix},\\ \nonumber \qquad \forall \tilde\eta_j, \tilde\eta_k \geq 0.   
\end{align}
where  $P_{(1,2,3)}$ is the collection of all  cyclic permutations of $(1,2,3)$. Observe that for every $i$-th diagonal block one may derive the \emph{sufficient and necessary condition for the symmetric part being non-indefinite}, and in particular for being non-negative definite since $\tilde\eta_l,\tilde\eta_m\geq 0$, which is imposing 
\begin{align}\label{eq:Bstable}
\tilde\eta_j \equiv \tilde\eta_k =:\tilde\eta^i, \quad \text{for $(i,j,k) \in P_{(1,2,3)}$}.
\end{align} 
This can be quickly checked by evaluating the analytical eigenvalues of the symmetric part of the diagonal block (\ref{eq:block}). 
This result is very important because condition  (\ref{eq:Bstable}) is also the \emph{sufficient and necessary condition  to guarantee (magnetic) energy stability} of the discrete equation with physical and artificial resistivity terms, i.e.
\begin{align}
\frac{d}{dt} E_\eta  = \int \tilde \B_h \cdot\frac{\partial \tilde \B_h}{\partial t}  = -  \int \left(  \eta \nabla\times\tilde \B_h +   \nabla_{\tilde{\eta}} \times\tilde \B_h\right) \cdot \nabla \times \tilde \B_h \leq 0,
\end{align}
which holds true also at the fully-discrete level by choosing a Crank-Nicolson time-discretization ($\theta_b=1/2$).
One should notice that condition (\ref{eq:Bstable}) was only conjectured in the finite-volume/finite-difference scheme of \cite{Fambri20} to preserve the symmetry of the resulting algebraic system for $\tilde\B$. 

In the numerical experiments, similarly to the numerical viscosity introduced in Sec. \ref{sec:acoustic} we set 
\begin{equation}
\begin{aligned}
&  \tilde\eta^i:= c_{\tilde\eta} \frac{\max\left( \lambda^{\alpha,j} \Delta x_j, \lambda^{\alpha,k} \Delta x_k\right)}{2}, \quad \text{for $(i,j,k) \in P_{(1,2,3)}$}. 
\end{aligned} \label{eq:sxyz}
\end{equation}
for a specified modulating parameter $c_{\tilde\eta}\geq0$.
\subsection{A semi-implicit time discretization} \label{sec:time}
 It will be convenient to write the right-hand-side of (\ref{eq:PDE1}) in an equivalent quasi-linear form, i.e.
$$
 \mathcal{L}^\alpha ( Q) \equiv  \hat{\mathcal{L}}^\alpha ( Q)\cdot Q.
$$
 and the corresponding discrete update-operator $ \mathcal{H}^{\alpha}$ as
\begin{align}
Q^{n+1}		&= \mathcal{H}^{\alpha}_{\Delta t } \cdot   Q^n.
\end{align}
As it was mentioned above,  see Section \ref{sec:res}, the resistive part is first split and decoupled from the ideal-MHD part. 
This is done by solving  a second-order Strang splitting scheme 
\begin{align}
Q^{n+1}		&= ( \mathcal{H}^{\eta}_{\Delta t/2 } \cdot \mathcal{H}^{v,p,b}_{\Delta t} \cdot \mathcal{H}^{\eta}_{\Delta t/2 }	)    Q^n.
\end{align}
The update-operator $\mathcal{H}^\eta$ follows directly from the definition of the implicit ($\theta$-method) time-in\-te\-gra\-tion described in Section \ref{sec:res}.
The viscous MHD part $\mathcal{H}^{v,p,b}_{\Delta t} $  is a semi-implicit time integration of the type
\begin{equation} \label{eq:semi-implicit}
\frac{Q^{n+1}-Q^n}{\Delta t} =  \mathcal{L}^v ( Q^{n}) +  \mathcal{L}^b ( Q^{n+\theta_b}) +  \mathcal{L}^p ( Q^{n+\theta_p})  
\end{equation}
To solve (\ref{eq:semi-implicit}), following \cite{Fambri20}, we first apply a \emph{recursive} and \emph{operator-splitting algorithm} as  
\begin{align}  
\frac{\tilde{Q}^{n+1,r+1}-Q^n}{\Delta t} &=  \mathcal{L}^v ( Q^{n}) + {\mathcal{L}}^b ( \tilde{Q}^{n+\theta_b,r+1})+ {\mathcal{L}}^p ( Q^{n+\theta_p,r}),\label{eq:SIb}\\
\frac{Q^{n+1,r+1}-Q^n}{\Delta t} &=  \mathcal{L}^v ( Q^{n}) + {\mathcal{L}}^b ( \tilde{Q}^{n+\theta_b,r+1}) + {\mathcal{L}}^p ( Q^{n+\theta_p,r+1}), \label{eq:SIp}
\\
\nonumber  &r=0, \ldots,  R;
\end{align} 
and then linearize in time in the sense of Picard obtaining 
\begin{align}
& \frac{\tilde{Q}^{n+1,r+1,s_b+1}-Q^n}{\Delta t} =   \mathcal{L}^v ( Q^{n}) +  \hat{\mathcal{L}}^b ( Q^{n+\theta_b,r,s_b})\cdot  \tilde{Q}^{n+\theta_b,r+1,s_b+1} + \hat{\mathcal{L}}^p ( Q^{n+\theta_p,r})\cdot  Q^{n+\theta_p,r} \label{eq:SIbPicard}  \\
& \frac{Q^{n+1,r+1,s_p+1}-Q^n}{\Delta t} =  \mathcal{L}^v ( Q^{n}) + \hat{\mathcal{L}}^b ( Q^{n+\theta_b,r,S_b})\cdot  \tilde{Q}^{n+\theta_b,r+1,S_b+1} + \hat{\mathcal{L}}^p ( Q^{n+\theta_p,r,s_p})\cdot  Q^{n+\theta_p,r+1,s_p+1}  \label{eq:SIpPicard} \\
 & \nonumber \qquad s_b=0, \ldots , S_p; \; s_p=0,\ldots ,S_b; \; r=0, \ldots ,  R.
\end{align} 
Here, three different iteration Picard indexes $(r,s_b,s_p)$ have been introduced: index $r$ is for the convergence of the splitting scheme (\ref{eq:SIb}-\ref{eq:SIp}) to the original equation (\ref{eq:semi-implicit}), while indexes $(s_b,s_p)$ are for the convergence of the solution of the linearized equations (\ref{eq:SIb}-\ref{eq:SIp}) to the original nonlinear system (\ref{eq:SIb}-\ref{eq:SIp}).  The final algorithm consists in two nested recursions: one \emph{external} recursion of $r=0$, $\ldots$, $R$ and two \emph{internal} recursion of $s_b=0$, $\ldots$, $S_p$ and $s_p=0$, $\ldots$, $S_b$.
The tilde in $\tilde{Q}$ is used to indicate that it is a temporary (or \emph{predictor})  solution.  This scheme may recall  an \emph{alternating direction implicit} (ADI) method adapted to a three-operator splitting of the type \emph{explicit-implicit-implicit}. In particular, by computing only a single Picard iteration the resulting implicit algorithm can be related to a Douglas-Rachford  ($\theta_b=\theta_p=1$)  or a  Douglas-Kim ($\theta_b=\theta_p=\frac{1}{2}$) scheme, see \cite{DouglasRachford,DouglasKim,Douglas62a,Douglas62b} and \cite{SplittingBook} and references therein for a review on splitting schemes. 

Now two alternative options appears for every iteration index: (i) whether one just performs a fixed number of Picard recursions aiming at preserving the time-accuracy of the scheme, or (ii) one iterates sequentially until convergence of a selected \emph{error}- or \emph{indicator}-norm within a prescribed tolerance $\epsilon$.
In the first case, depending on the non-linearity and stiffness of the algebraic systems, the computational costs may be heavily reduced while maintaining the desired (2nd) order of accuracy in time. 
In the second case, all the conservation properties proved at the semi-discrete level may be also enforced at the fully discrete level and up to a prescribed tolerance. 
As an example, in section \ref{sec:split-alf} we show that stopping the iterations according to a prescribed tolerance, i.e. $S_b=S_b(\epsilon)$, we are able to preserve some physical quantities, e.g. magnetic energy or magnetic helicity.
 We anticipate that, in the numerical tests where it's not specified, we set $R=S_p=1$ and $S_b=0$. 
A rigorous analysis of the proposed time-discretization is out of the scope of this paper.

According to this, a three operator-splitting is proposed as an efficient solver for the decoupled convection-diffusion, acoustic and Alfv\'enic terms. In particular, the second order MUSCL-FV scheme is used to compute the convective part $\mathcal{L}^v$, while the implicit theta-method is used to solve the Alfv\'enic-acoustic terms, $\mathcal{L}^b$ and $\mathcal{L}^p$.
The operator splitting has multiple advantages when dealing with implicit integrators. First, it allows to split the original fully-coupled problem into a sequence of smaller and decoupled resolvents. Indeed, instead of inverting the fully coupled system  (\ref{eq:SIb}), i.e. solving for $Q^{n+1,r+1}$ 
\begin{equation}
\begin{aligned} 
 \left(1 - \hat{\mathcal{L}}_{\theta_b \Delta t }^b ( Q^{n+\theta_b,r}) - \hat{\mathcal{L}}^p_{\theta_p \Delta t} ( Q^{n+\theta_p,r})\right) \cdot {Q}^{n+1,r+1} & = \text{r.h.s.}, 
\end{aligned}
\end{equation}
one has to solve the two Alfv\'enic (\ref{eq:SIb}) and acoustic (\ref{eq:SIp})  systems separately, i.e. solving first for $\tilde{Q}^{n+1,r+1}$ 
\begin{equation} \label{eq:implicit1}
\begin{aligned} 
 \left(1 - \hat{\mathcal{L}}_{\theta_b \Delta t }^b ( Q^{n+\theta_b,r})\right) \cdot \tilde{Q}^{n+1,r+1} & = \tilde{R}^r,  
\end{aligned}
\end{equation}
and then for  ${Q}^{n+1,r+1}$ 
\begin{equation} \label{eq:implicit2}
\begin{aligned}  
 \left(1 - \hat{\mathcal{L}}^p_{\theta_p \Delta t} ( Q^{n+\theta_p,r})\right) \cdot  Q^{n+1,r+1}  & = R^r,  \\
\end{aligned}
\end{equation}
respectively, where 
$$
\tilde{R}^n = Q^n +  \Delta t\mathcal{L}^v ( Q^{n}) + (1-\theta_b) \Delta t\hat{\mathcal{L}}^b ( Q^{n+\theta_b,r})\cdot  \tilde{Q}^{n} +      \Delta t\hat{\mathcal{L}}^p ( Q^{n+\theta_p,r})\cdot  Q^{n+\theta_p,r},
$$
$$
R^r =   Q^n +  \Delta t \mathcal{L}^v ( Q^{n}) + \Delta t \hat{\mathcal{L}}^b ( Q^{n+\theta_b,r})\cdot  \tilde{Q}^{n+\theta_b,r+1} + \Delta t (1-\theta_p) \hat{\mathcal{L}}^p ( Q^{n+\theta_p,r})\cdot  Q^{n}. $$

Ultimately, the non-linear viscous MHD operator  $\mathcal{H}^{v,p,b}_{\Delta t}$  is built as a nested recursive procedure solving the convection-diffusion, sonic and Alfv\'enic steps described in Sec. \ref{sec:split}. 
 The double and nested recursion is due to the non-linearity of the implicit terms in the two consecutive steps in (\ref{eq:implicit1}-\ref{eq:implicit2}). 
In this process, the explicit terms of the convection-diffusion step are computed only once per time-step. 

\smallskip
The standard \emph{hyperbolic-parabolic} $\CFL$ restriction on the computational time-step for Cartesian meshes reads as
\begin{equation} 
 \Delta t = \Delta t(\lambda^h)  := \CFL \frac{1}{ \frac{\max |\lambda^h_x|}{\Delta x} + \frac{ \max |\lambda^h_y|}{\Delta y}+ \frac{ \max |\lambda^h_z|}{\Delta z} + 2 \lambda^p \left( \frac{1}{\Delta x^2} + \frac{1}{\Delta y^2} + \frac{1}{\Delta z^2} \right) }. \label{eq:CFL}
\end{equation} 
Here, $\lambda^h$ is  a properly chosen set of   \emph{hyperbolic}  eigenvalues, $\lambda^p>0$ is the set of  \emph{parabolic}  eigenvalues of the \emph{physical} viscous terms in the PDE.
In general, numerical stability assumes $\lambda^h$ and $\lambda^p$
to be the set of eigenvalues of the Jacobian of the \emph{explicit} fluxes and  the  penalty term that refers to the \emph{explicit} parabolic terms, respectively. In our formulation,  the physical viscosity and heat-conduction are the only parabolic terms solved in the explicit step, see Sec. \ref{sec:convdiff}. Then we chose $\lambda^p=(4/3)\mu /\rho  +\kappa/(c_v \rho)$.  
Regarding $\lambda^h$ there is an important note. 
Indeed, for stability it is important that the spectral radius  $s_h= \max(\lambda^h)$ used to evaluate the time-step size in (\ref{eq:CFL}) must be larger than the spectral radius of the explicit terms in the discrete equations.
 In the presented work, the explicit terms in the discrete equations amount to the central convective and viscous fluxes, 
plus the Rusanov stabilization term in eq.s (\ref{eq:Rusanovx}-\ref{eq:Rusanovz}) which takes the form of a parabolic diffusive term. Then, the corresponding numerical viscosity is proportional to $\frac{1}{2} s_\alpha \Delta x^{N}$,
 by assuming the jumps $\w^+ - \w^-$ converging with the order of accuracy $N$. 
In principle, artificial diffusion may be added to the scheme by selecting different eigenspectra for the Rusanov viscosity $s_\alpha=s_\alpha(\lambda^\alpha)$. 
This may eventually become useful, e.g., to dissipate fast modes in convection dominated flows. Then, a general constraint on the choice of $\lambda^h$ is
$$
\max(\lambda^h) =  s_h \geq s_\alpha = \max(\lambda^\alpha).
$$
Finally, $\CFL\in(0,1]$ is the \lq\lq\CFL\rq\rq\, security parameter, sometime called \emph{Courant number}. In this paper, we may   refer to a $\lambda$-based Courant number, defined with respect to the inviscid equations, i.e. $$
 \CFL_{\lambda} = \frac{\Delta t_\lambda}{s(\lambda) \Delta x}.$$
 For semi-implicit schemes, the $\CFL$ parameter is used as a security factor, while the time-step scale will be mainly identified by the reference characteristics spectra, i.e. $\lambda^h$.
In the line-legends of the plots, first order FV scheme is simply labeled with \lq\lq FV\rq\rq\,, second-order MUSCL TVD FV-scheme with \lq\lq FV2\rq\rq\,.

\section{Numerical validation}
\label{sec:results} 

\subsection{1D Riemann problems}

To test the robustness and accuracy of the method,  in this section we first solve a series of Riemann or shock-tube problems for the MHD equations, i.e. the initial condition is 
$$
(\rho,   \u, p, \B ) = \left\{ \begin{array}{cr} (\rho_L,   \u_L, p_L, \B_L ) & \text{for $x \leq 0$} \\  (\rho_R,   \u_R, p_R, \B_R )   &  \text{for $x > 0$} \end{array} \right.
$$
 
In particular, we want to check that, even selecting Courant number based on the pure convection instead of the full magneto-acoustic waves, the waves can still be well resolved. 
The correct wave-speeds are resolved even by using large time-steps, with the faster modes being solved with  more numerical diffusion with respect to the low-modes. 
Note that the scheme is agnostic of the full fan of the MHD modes and, consequently, it can be considered as a \emph{1-wave solver}.   
In this sense, our novel semi-implicit hybrid FV-FEEC scheme is shown to be very robust and low dissipative if compared to other low order 1-wave explicit solvers. Thanks to spectral analysis provided by \cite{roebalsara}, full knowledge of the eigenspectrum of ideal MHD is given and exact or approximate Riemann solvers may be built on that.

In the next, different numerical configurations will be used to verify the capabilities of our novel semi-implicit hybrid FV-FEEC method. 
Consider that, only different couples of variables are plotted. These may differ between the different Riemann problems, in the aim of reporting the most significant results. Many of the following test configurations comes from the following reference papers \cite{BrioWu,RyuJones,DaiWoodward}, where the spectral properties are deeply discussed. 
The different initial conditions are reported in table \ref{tab:rp}. The specific heat capacity ratio is assumed to be $\gamma=5/3$ (monoatomic ideal gases) for all the tests.
 The $\CFL$ security parameter is set to $0.9$. The spatial resolution for the $V_0$ space is  set to $\Delta x= 1/1000$ for all the tests RP1-5, only for the isolated stationary contact wave we set $\Delta x= 1/100$.
In the plots, the implicit weight of the $\theta$-method are sometimes (only when considered interesting) written as $\theta = ( \theta_b, \theta_p)$, and their default value is assumed to be equal to one. Similarly, the default Rusanov viscosity is assumed to be evaluated with respect to the convective eigenvalues, i.e. $s_\alpha \equiv s_v =s(\lambda^{v})$, as well as the time-step size $\Delta t = \Delta t(\lambda^{v})$.
Again, these are just default values, and they may vary only if they are explicitly specified in the text.
A theoretical investigation about the stability or robustness of the presented 3-splitting FV-FEEC scheme depending on the implicit weights $\theta_p$ and $\theta_b$ is out of the scopes of this manuscript.

Consider that all the following shock-dominated problems have been solved using time-steps of the convective scale $\Delta t(\lambda^v)$ with $\CFL$ security parameter equal to 0.9. In literature, stable results are usually obtained using time-steps $\Delta t(\lambda^{\MHD})$, and meanwhile lower $\CFL$ parameters.

In the plots, the \lq\lq exact\rq\rq\, solutions come from S.A.E.G. Falle \cite{fallemhd,falle2}, while the so-called \lq\lq reference\rq\rq\, solutions are obtained by our semi-implicit hybrid FV-FEEC on a very refined grid.

\begin{table}[!t]
 \caption{Initial states left (L) and right (R) for the primitive variables for the selected  Riemann problems of the ideal classical MHD equations. 
 In all cases $\gamma=5/3$. Final time of the simulations $t_f$ and initial position of the discontinuity $x_d$ are also reported. } 
\begin{center} 
 \begin{tabular}{!{\extracolsep{-5pt}}r c ccc c ccc cc!{}}
 \hline
 Case & $\rho$ & $u_x$ & $u_y$ & $u_z$ & $p$ & $B_x/\sqrt{4 \pi}$ & $B_y/\sqrt{4 \pi}$ & $B_z/\sqrt{4 \pi}$  & $t_f$ &  $x_d$       \\ 
 \hline    
	RP0 L: &  1   		&  0				& 0   		& 0     		&  $10^3$		&  $100$		& 0 												& $100$			&  $10^3 $  &  0	\\  
			R: &  0.125		&  0				& 0   		& 0     		&  $10^3$		& $100$ 			& 0 												& $100$  		&  		 &  		\\ 			\hline
	RP1 L: &  1				&  0    		& 0   		& 0    			&  1				& $3/4\sqrt{4 \pi}$ 									&  $\phantom{+}\sqrt{4 \pi}$						& 0      									&  0.1   &  0		\\
			R: &  0.125		&  0				& 0   		& 0    			&  0.1     	& $3/4\sqrt{4 \pi}$ 									& $-\sqrt{4 \pi}$  										& 0      									&  		 &  		\\ 			\hline
	RP2 L: &  1.08		&  1.2			& 0.01   	& 0.5  			&  0.95    	& $2$				&  $3.6    $	& $2    $		&  0.2   &  -0.1	\\
			R: &  0.9891	&  -0.0131	& 0.0269 	& 0.010037	&  0.97159 	& $2$				&  $4.0244 $	& $2.0026  $&  		 &  		\\ 			\hline
	RP3 L: &  1.7			&  0    		& 0   		& 0     		&  1.7			&$ 3.899398$ &  $3.544908$ 	& $0.0$      							&  0.04   &  0		\\
			R: &  0.2			&  0    		& 0   		& -1.496891	&  0.2   		& $3.899398$ &  $2.785898$ 	& $2.192064$	&  		 &  		\\ 			\hline
	RP4 L: &  1   		&  0    		& 0   		& 0  				&  1    		& $1.3\sqrt{4 \pi}$ 									&  $\phantom{+}\sqrt{4 \pi}$  	 				& 0           						&  0.16   &  0		\\
			R: &  0.4			&  0    		& 0   		& 0  				&  0.4     	& $1.3\sqrt{4 \pi}$ 									&  $-\sqrt{4 \pi}$  										& 0   										&  		 &  		\\ 			\hline
	RP5 L: &  1   		&  0    		& 0   		& 0  				&  $\sqrt{4 \pi}$   		& $\sqrt{4 \pi}$ 										&  $ \phantom{+}\sqrt{4 \pi}$					& 0           						&  0.03   &  0	\\ 		 
			R: &  0.2			&  0    		& 0   		& 0  				&  0.2     	& $\sqrt{4 \pi}$ 										&  $\sqrt{4 \pi}\cos(3)$								& $\sqrt{4 \pi}\sin(3)$              &  		 &  		\\ 		 
 \hline
 \end{tabular}
\end{center} 
 \label{tab:rp}
\end{table}

\begin{figure}[!t]
\begin{center}
\includegraphics[width=0.56\textwidth]{\MyFigFolder/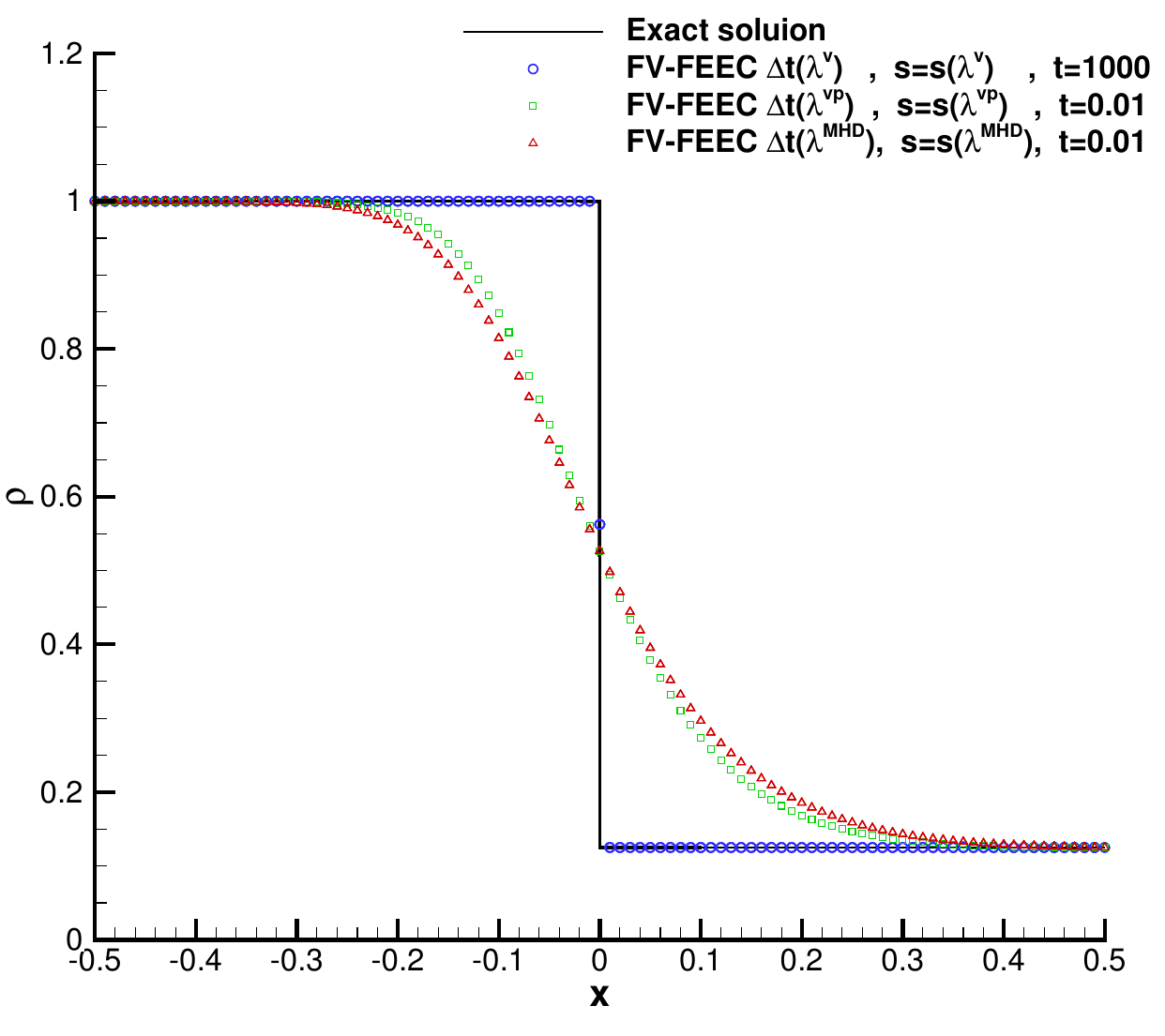}
\caption{Exact and numerical solutions 
 for the Riemann problem RP0  ($\Delta x= 1/100$). The low Mach isolated, stationary contact wave is exactly preserved even at long time $t_f=1000$ for our hybrid FV-FEEC method ($s(\lambda^v)$). The FV-FEEC solution obtained by setting a higher numerical diffusion, i.e. $s=s(\lambda^b)$ or $s=s(\lambda^{\MHD})$, is also plotted at time $t=0.01$. 
(see colored version online) } 
\label{fig:rp0}
\end{center}
\end{figure}

\smallskip 
\paragraph{RP0} The first  test RP0 represents an isolated, stationary contact  discontinuity in a \emph{low} Mach regime. 
Thanks to the chosen implicit discretization, the stabilization terms of the explicit fluxes depends only on the fluid velocity.
 For this reason, the contact discontinuity is exactly preserved independently of the sonic, and the Alfv\'enic Mach numbers. In Fig. \ref{fig:rp0}, we compare three different numerical results: the original FV-FEEC solution at long time $t_f=1000$, where the time-step and the Rusanov viscosity are chosen to depend only on the convective eigenvalues, i.e. $\Delta t(\lambda^{v})$, $s(\lambda^{v})$; 
two FV-FEEC solutions   where the time-step \emph{and} Rusanov viscosity are evaluated with respect to the Alfv\'enic and the full MHD eigenvalues, $\lambda^{b}$ and $\lambda^{\MHD}$, respectively, see (\ref{eq:lambdab}) and (\ref{eq:mhd.wavespeeds}). Consider that the original scheme with  $s=s(\lambda^{v})$ naturally preserve the contact discontinuity independently of the time-scales and the sonic or Alfv\'enic Mach numbers. Indeed, this test can be solved even in one single time-step. Note also that the last \emph{viscous} results do not differ much, because for this magnetohydrodynamic system the spectral radius of the Alv\'enic and acoustic systems are comparable.

\begin{figure}[!t]
\begin{center}
\includegraphics[width=0.48\textwidth]{\MyFigFolder/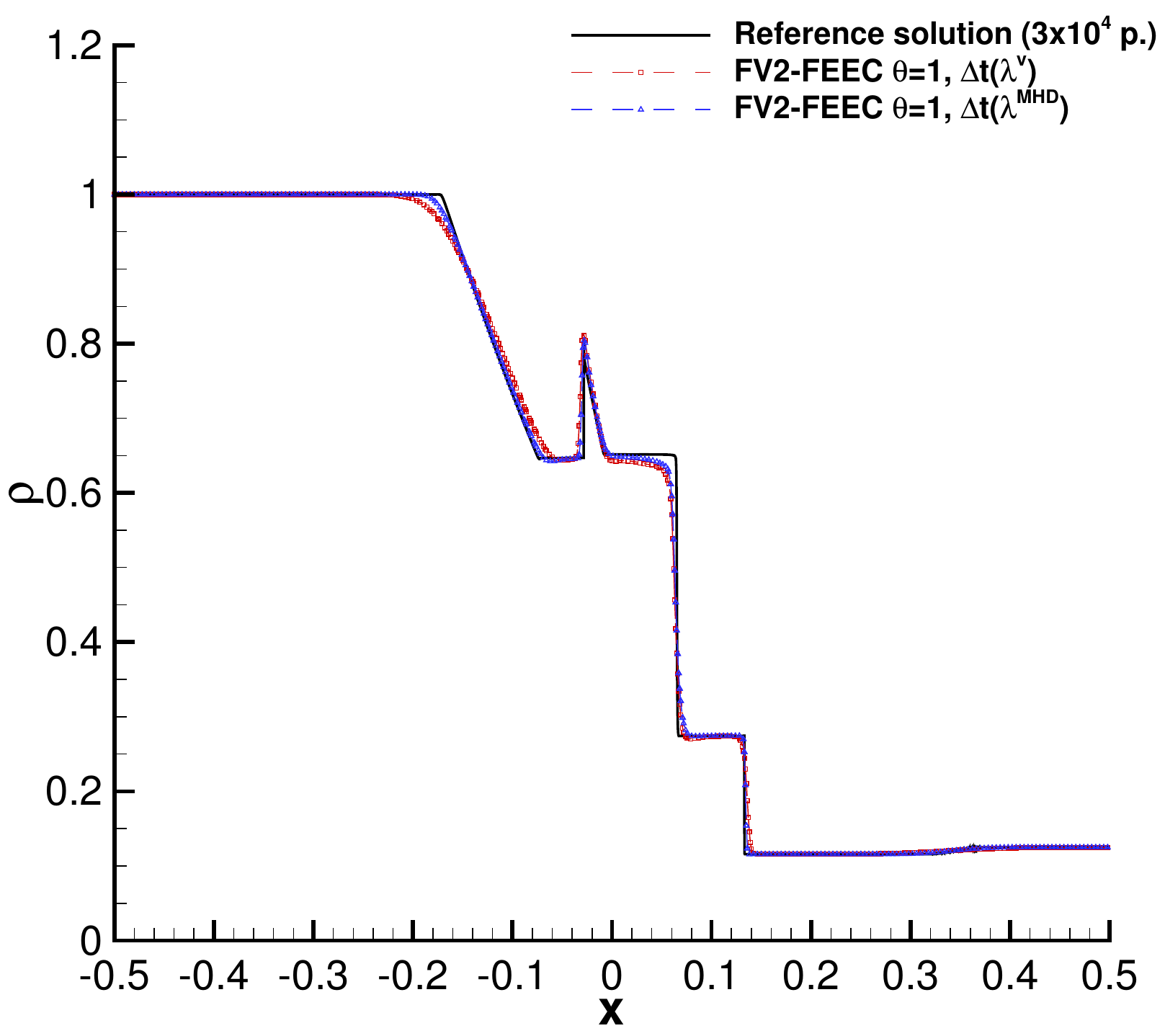}  
\includegraphics[width=0.48\textwidth]{\MyFigFolder/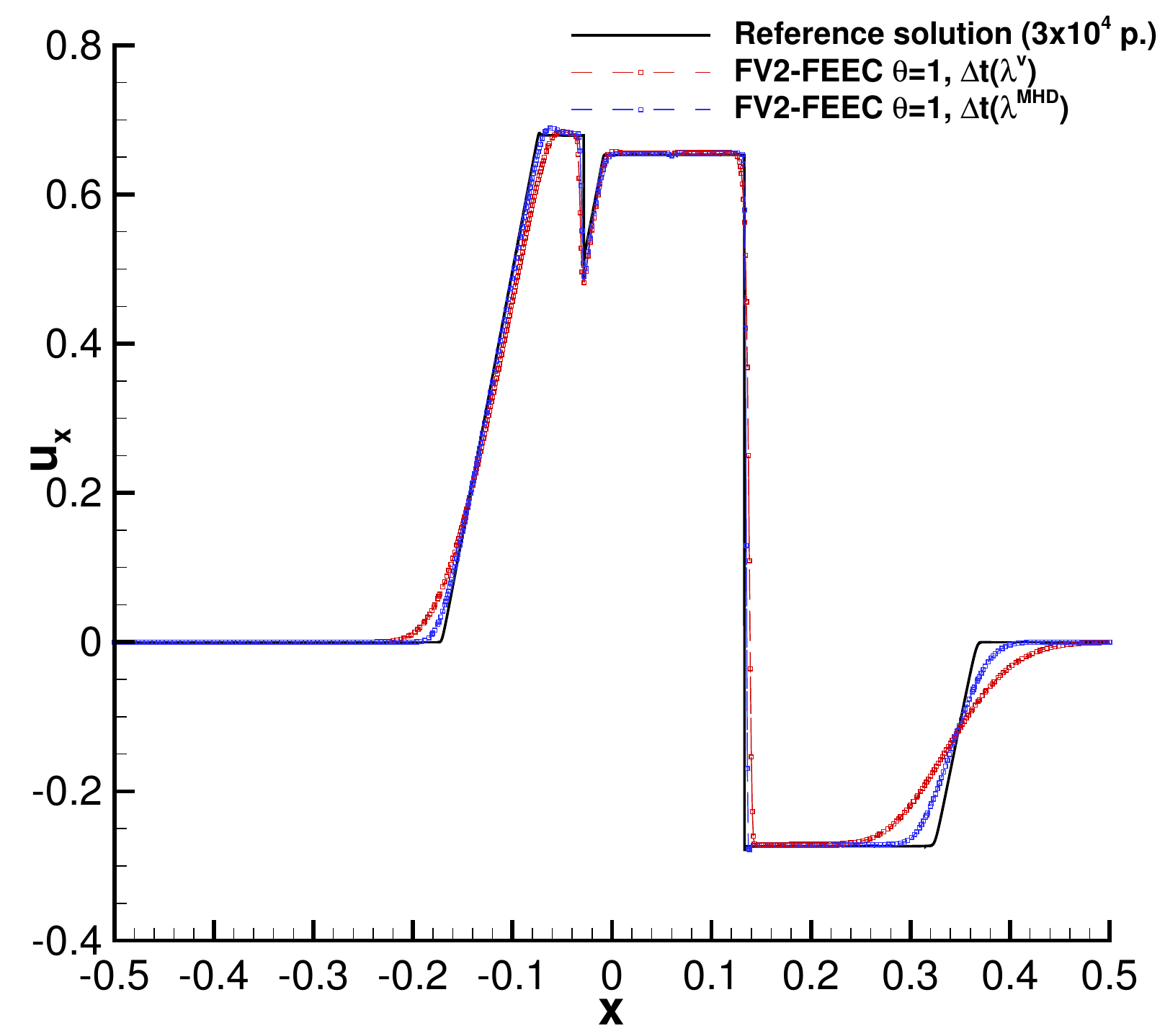}  
\caption{Exact and numerical solution 
 for the MHD Riemann problem RP1 at $t=0.1$. Density (left) and horizontal velocity $u_x$ (right) are plotted, comparing the numerical solutions obtained after choosing different time-step scales, i.e. the convection scale $\Delta t(\lambda^{v})$ (red squares), and the MHD scale $\Delta t(\lambda^{\MHD})$ (blue deltas). (see colored version online) 
} 
\label{fig:rp1}
\end{center}
\end{figure}

\smallskip 
\paragraph{RP1} The second test is the well known Brio-Wu Riemann problem \cite{BrioWu}, which can be seen as an MHD counterpart of the Sod shock-tube problem in compressible gas-dynamics. In this test we have a fast-rarefaction, a slow compound wave (left-traveling), a contact discontinuity, a slow shock and a fast rarefaction (right-traveling). The compound wave is a non-classical wave and its existence is due to the non-convexity of the MHD equations, but we refer to  \cite{torrilhon,Torrilhon2004} for a more detailed discussion on the topic. 
All the waves are correctly approximated. Note that the faster waves are more affected by numerical diffusion mainly due to the time-error when the time-step is larger than the time-scale of the faster modes $\Delta t(\lambda^v) >  \Delta t(\lambda^{\MHD})$, see e.g. the plot of the horizontal velocity in Fig. \ref{fig:rp1}.
This is a good news, meaning that we have time-convergence and stability even when solving time-steps larger than the physical time-scales of the problem.

\begin{figure}[!t]
\begin{center}
\includegraphics[width=0.48\textwidth]{\MyFigFolder/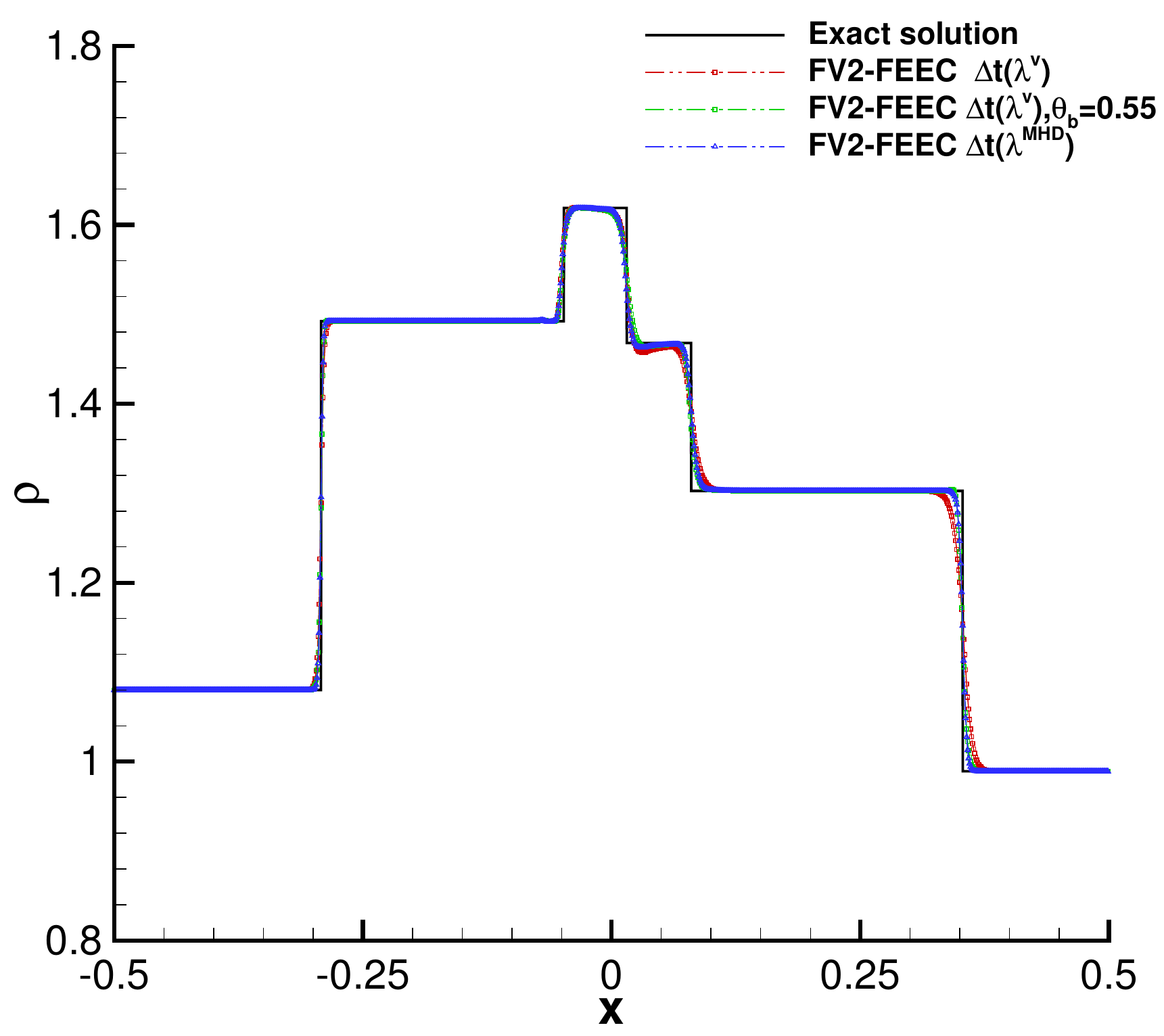}    
\includegraphics[width=0.48\textwidth]{\MyFigFolder/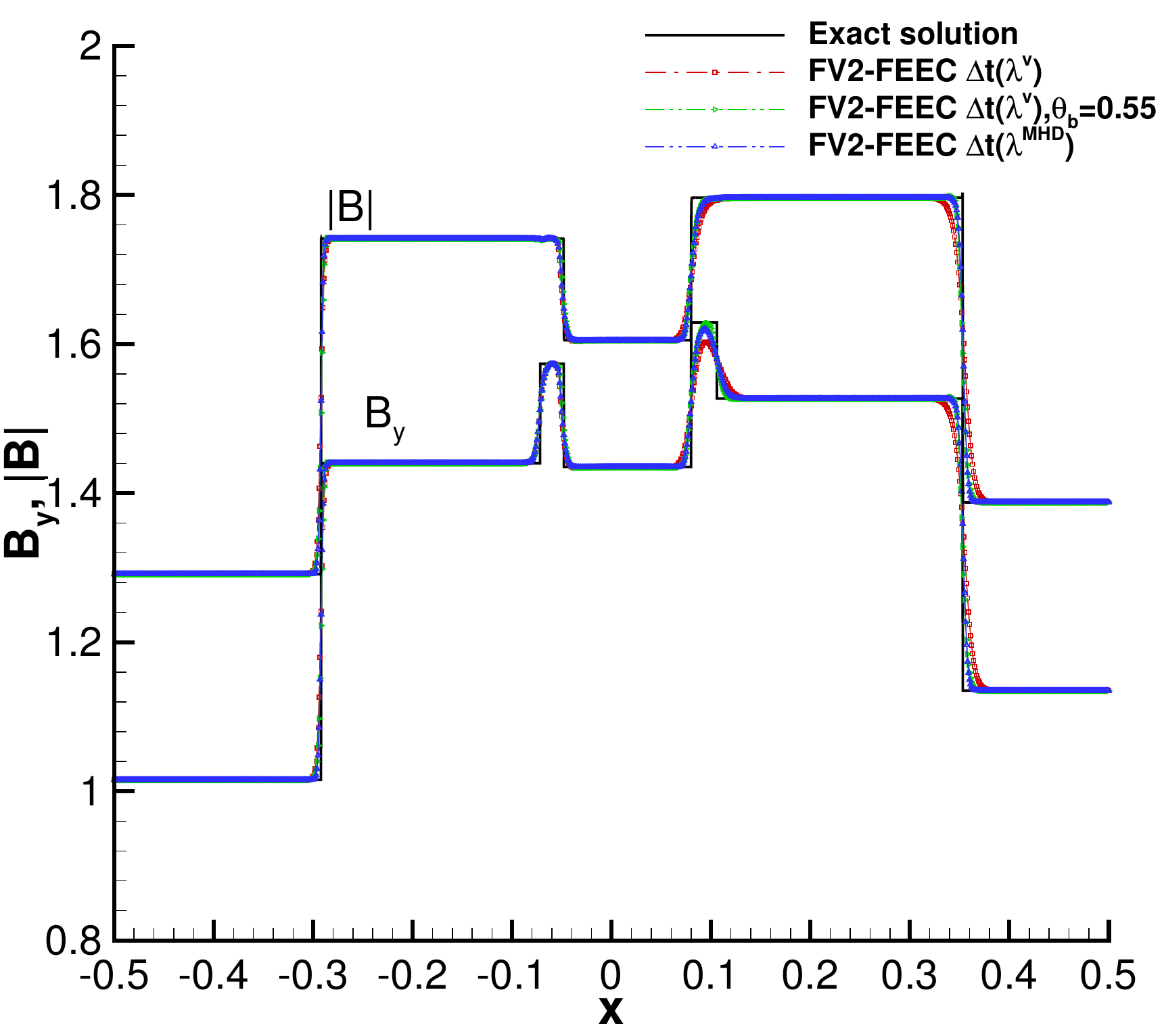}   
\caption{Exact and numerical solution 
 for the Riemann problem  RP2 at $t=0.2$. Density (left column) and magnetic field $B_y$ and $|B|$ (right column) are plotted, comparing the numerical solutions obtained after choosing different time-step scales, i.e. the convection scale $\Delta t(\lambda^{v})$ with $\theta_b=1$ (red squares),  $\theta_b=0.55$ (green squares), and the MHD scale $\Delta t(\lambda^{\MHD})$ (blue deltas). The vertical component $B_y$ is plotted next to the magnitude $|B|$ to show the correct approximation of the two rotational waves. (see colored version online) 
} 
\label{fig:rp2}
\end{center}
\end{figure}

\smallskip
\paragraph{RP2} 
Test RP2 is instead the special case where only discontinuous waves are produced, see \cite{DaiWoodward}: two (left- and right- traveling) fast shocks, followed by two rotational discontinuities, two slow shocks, and the central contact wave. 
All the seven discontinuities have been correctly approximated by our hybrid FV-FEEC solver. 
Numerical results for the density and magnetic field are plotted in Fig. \ref{fig:rp2}, where vertical component $B_y$ is plotted next to the magnitude $|B|$ to check the correct approximation of the two rotational waves. In this case, three different numerical set up are compared, which differ in the time-discretization: the coarser  ($\Delta t(\lambda^v)$), and the finer ($\Delta t(\lambda^{\MHD})$) time discretizations and the results obtained exploiting the higher-order of accuracy of the $\theta$-method still using the coarser time-discretization with $\theta_b=0.55$. The results of the $\theta$-method are very promising. 
Indeed, the $\theta$ method allows for faster convergence of the implicit solver \cite{GreenspanCasulli}, since the off-diagonals are proportional to $\theta^2$, see eq. (\ref{eq:disc-weakm0}), meanwhile being higher-order accurate. 
For this test, the $\theta$-method gave comparable and sometimes better approximation of the seven discontinuities, even compared with the highly resolved case $\Delta t(\lambda^{\MHD})$.

\begin{figure}[!t]
\begin{center}
\includegraphics[width=0.48\textwidth]{\MyFigFolder/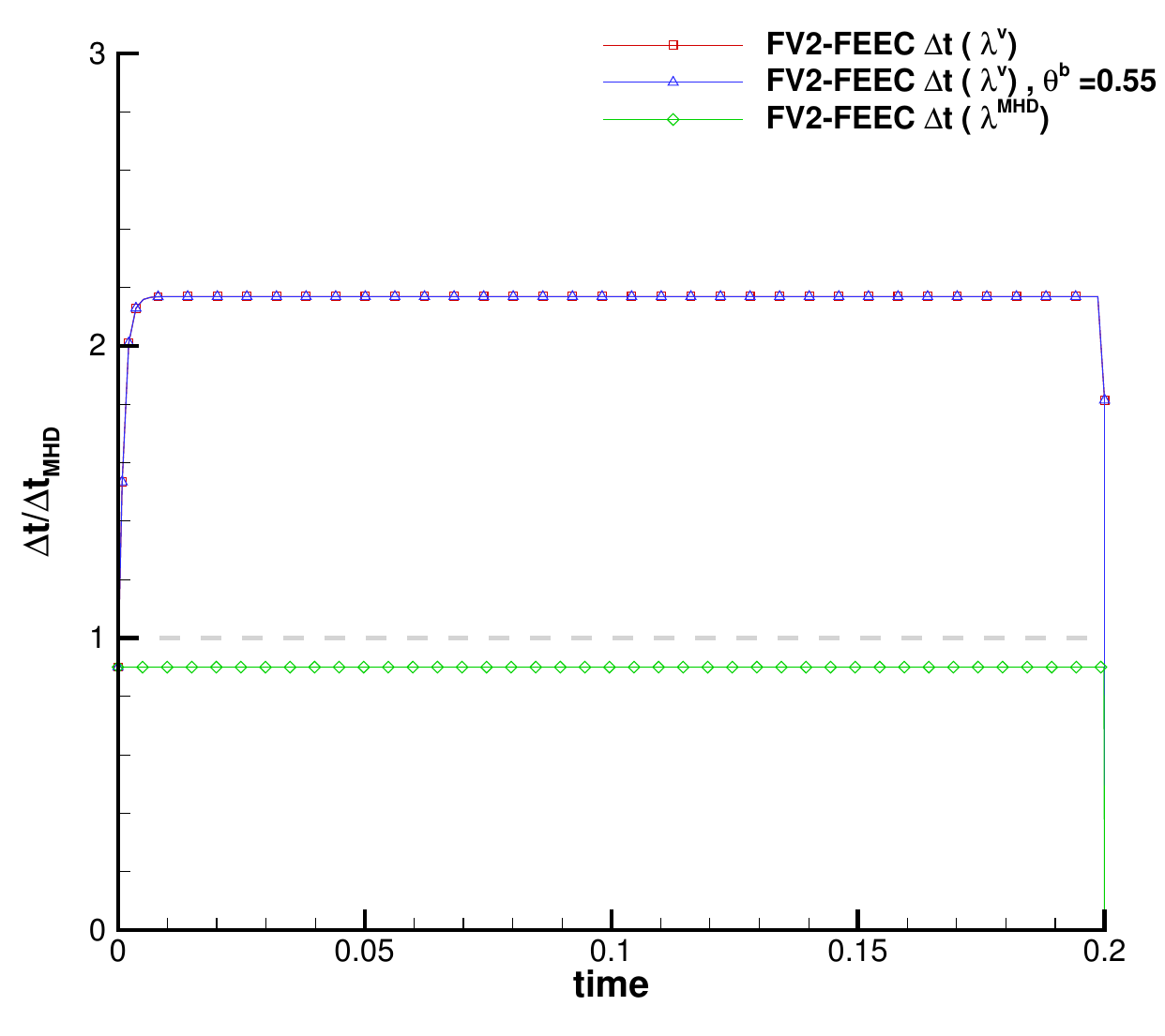}    
\includegraphics[width=0.48\textwidth]{\MyFigFolder/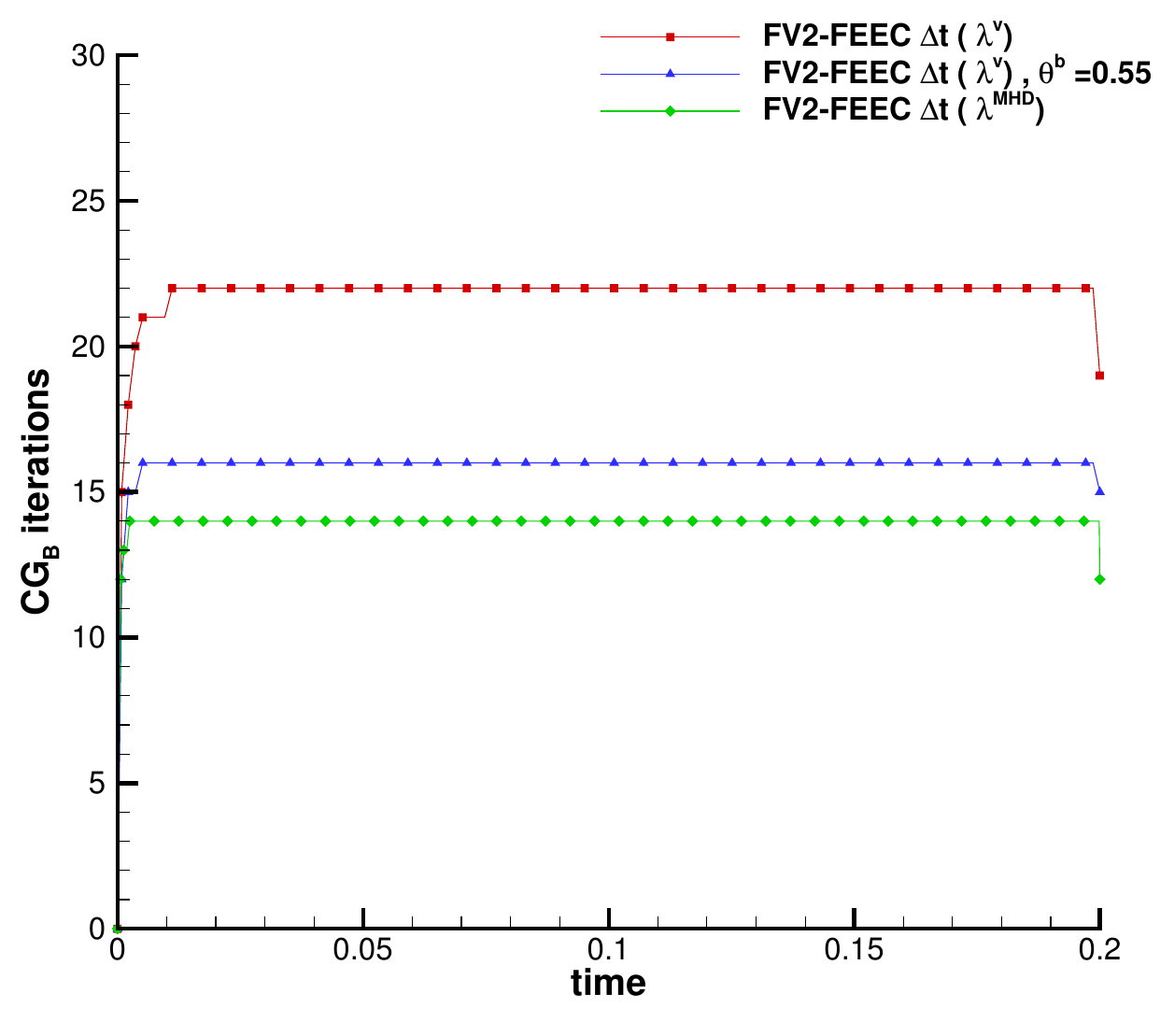}   \\
\includegraphics[width=0.48\textwidth]{\MyFigFolder/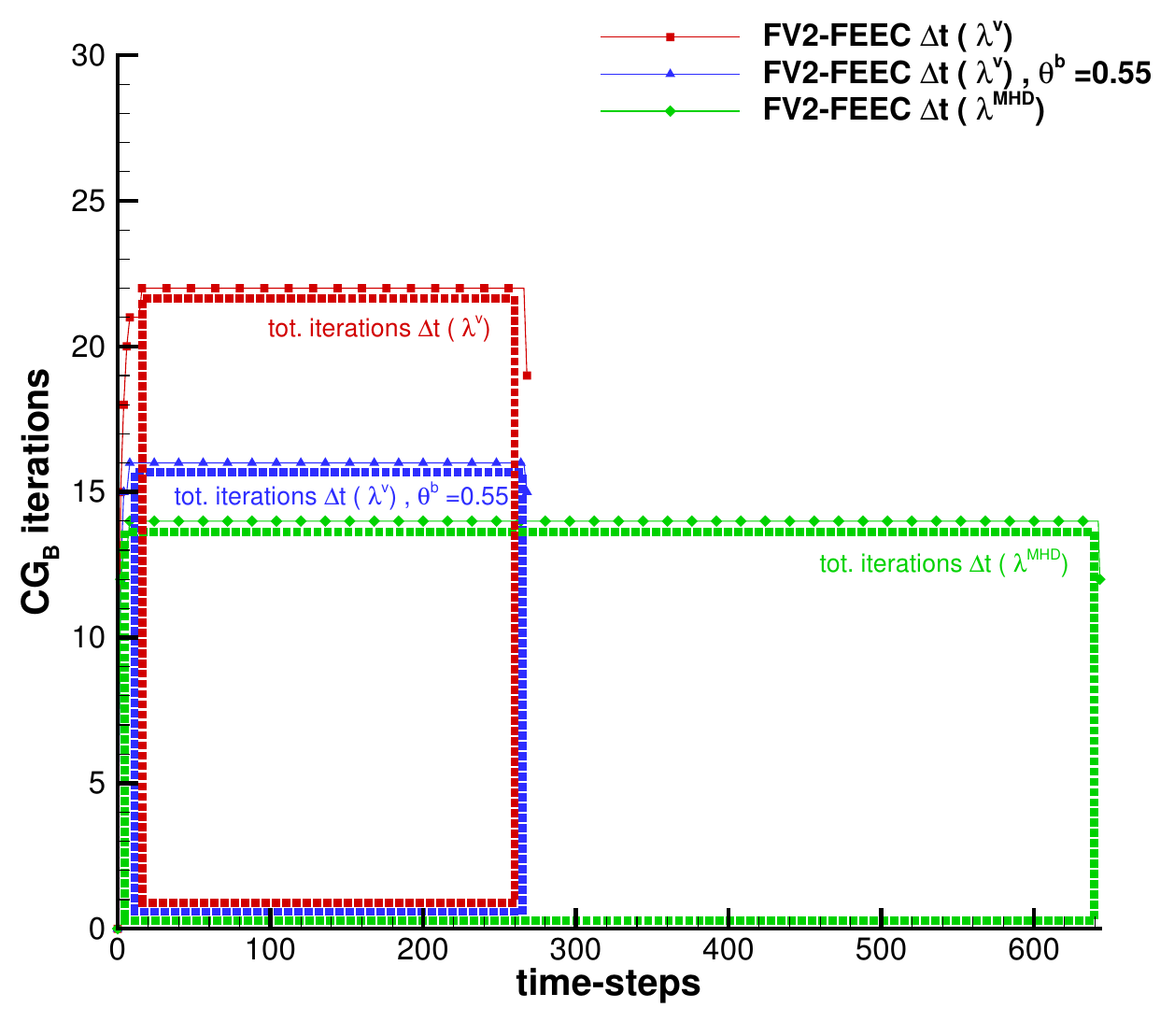}    
\includegraphics[width=0.48\textwidth]{\MyFigFolder/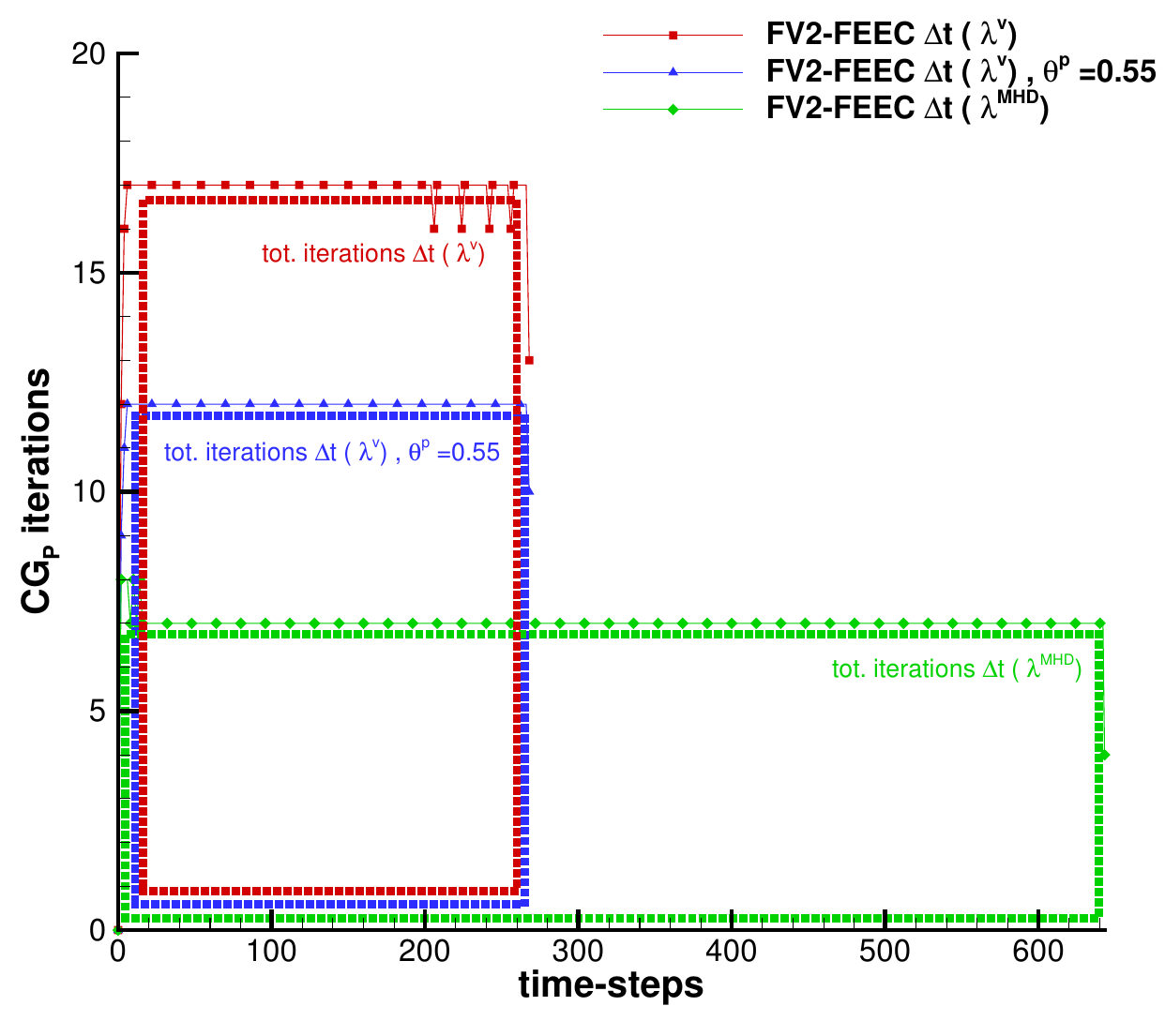}   \\
\caption{Comparison of the numerical performances obtained by solving the  RP2 test and tuning the time-discretization parameters, i.e. the time-step size $\Delta t$ and tuning the $\theta$-method parameter. First row: the time evolution of so-called \emph{effective Courant number} $\Delta t/\Delta t(\lambda^{\MHD})$ \cite{Fambri20} is plotted at the left; at the right, the time evolution of the number of iterations of the conjugate gradient, needed to solve the Alfv\'enic algebraic system. 
Second row: the \lq\lq time-step\rq\rq\, evolution of the number of iterations needed to solve the Alfv\'enic algebraic system (left) and the acoustic algebraic system (right). The subtended (colored) area gives an estimate of the total number of CG iterations required to complete the simulation.   (see colored version online) 
} 
\label{fig:rp2b}
\end{center}
\end{figure}

\smallskip
We drew this simple MHD shock-tube test to compare the computational advantages of our semi-implicit hybrid FV-FEEC method obtained by tuning these parameters. Note that the advantages of using implicit methods are evident in the low Mach regime, where the $\MHD$ time-scale $\Delta t(\lambda^{\MHD})$ is prohibitive, see e.g. RP0 where one single large time-step was sufficient to complete the simulation without artificial dissipation.
In Fig. \ref{fig:rp2} the time evolution of the \emph{effective Courant number} of our semi-implicit simulations $\Delta t/\Delta t(\lambda^{\MHD})$ \cite{Fambri20}, next to three different plots of the number of  iterations needed for the convergence of the conjugate-gradient solver. Remember that these data consider only solvers \emph{without any preconditioning}.
The precise number of iterations is not important since it may change in different implementations, e.g. the choice of the tolerance-to-convergence or the definition of the functional to be minimized. 
The important info we may get from these plots is the following: obviously, implicit methods allows for effective Courant numbers larger than 1, and as a consequence fewer time-steps are required to complete the simulations; less obvious and more interesting are the plots of the CG iterations:
 first, the higher computational efficiency of the $\theta$-method is verified by the numerical results, i.e. less CG iterations are required if $\theta<1$; 
second, choosing the maximum Courant number, i.e. the maximum time-step, does not necessarily mean to get the best computational performance. Indeed, in the second row of Fig. \ref{fig:rp2b}, the number of CG iterations are plotted as a function of the time-step number. 
In particular, the area bounded by the graphs gives an estimate of the total number of iterations required in the full simulation. 
One can see that, for this test RP2, only using the $\theta$-method gives a clear advantage, while only modulating the effective Courant number, i.e. choosing between $\Delta t(\lambda^v)$ or $\Delta t(\lambda^{\MHD})$, do not change the effective performance.
 Obviously, a spectral analysis of the algebraic system is needed to estimate the maximum number of iterations that are necessary to reach convergence, depending on the time-step size. In this way one could give rigorous estimate of the gain/loss of computational performances, eventually designing an optimal preconditioning that may minimize the total number of iterations 
\begin{align*}
 \text{CG}_{\text{tot.}} & = \text{CG}_{\text{iter.}} \times n^ \circ_{\text{time steps}}  \sim \frac{ \text{CG}_{\text{iter.}}(N_p, \theta \Delta t (N_p)) }{\Delta t (N_p)}   
\end{align*}
where $\text{CG}_{\text{iter.}}=\text{CG}_{\text{iter.}}(N_p, \theta \Delta t (N_p))$ is the average number of iterations, and the dependency of the time-step on the total number of discretization point $N_p$ is displayed.
This can surely be an interesting topic of our future research.

\begin{figure}[!t]
\begin{center}
\includegraphics[width=0.48\textwidth]{\MyFigFolder/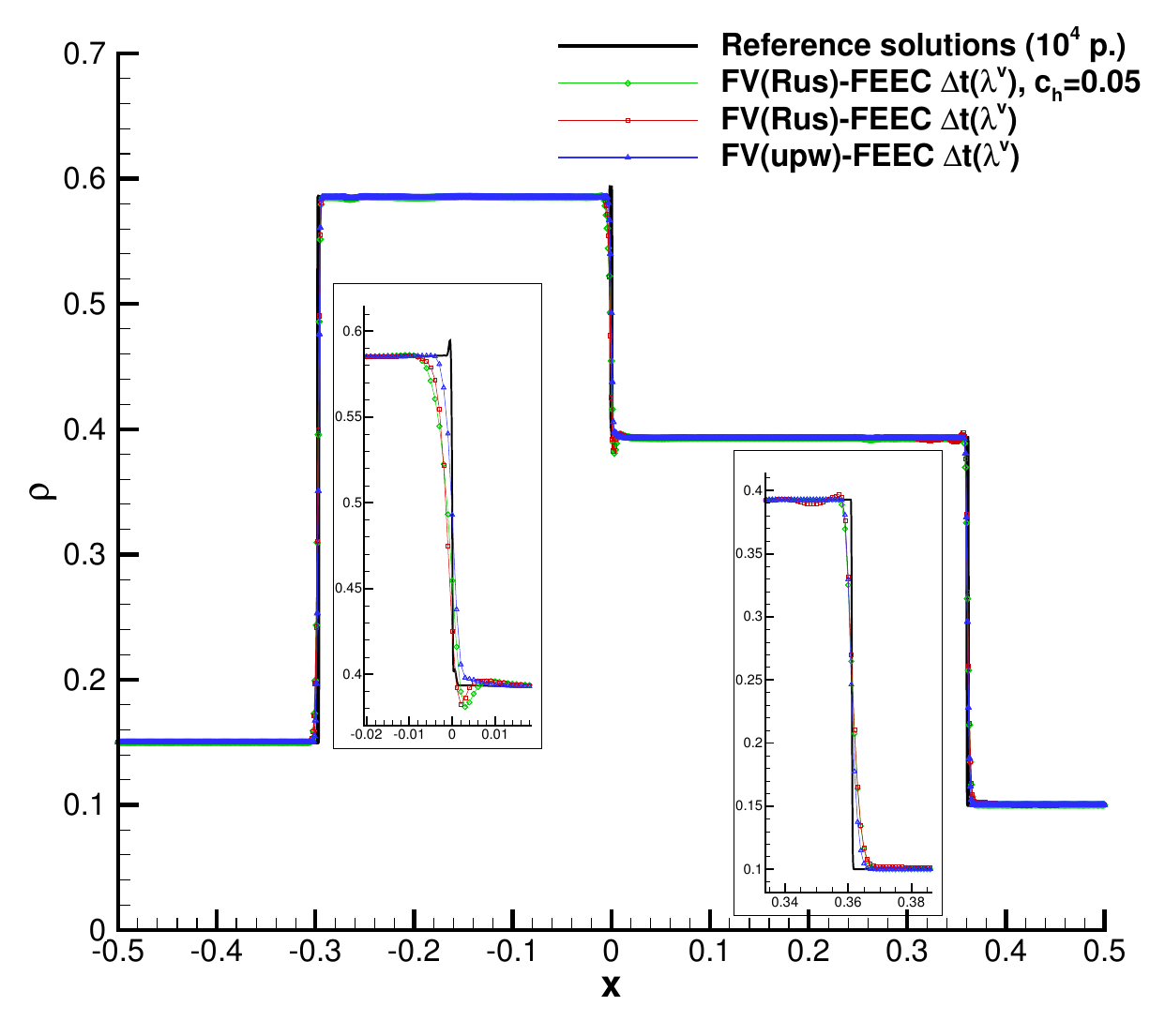}    
\includegraphics[width=0.48\textwidth]{\MyFigFolder/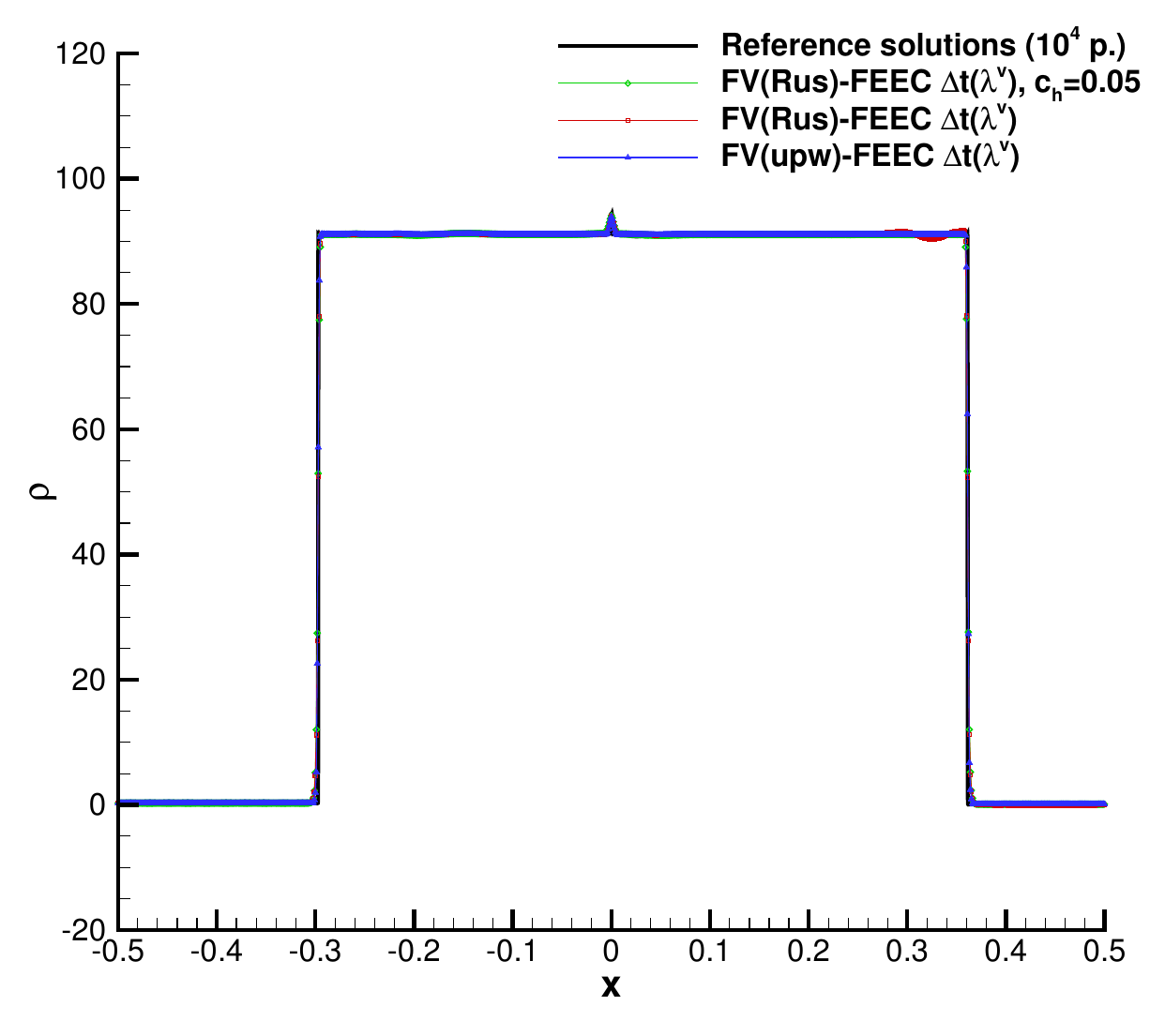}   
\caption{Reference and numerical solution 
 for the Riemann problem  RP3 at $t=0.04$.  Density (left column) and  magnetic field component $B_y$ (right column) are plotted, comparing the numerical solutions obtained after choosing different numerical convective fluxes, i.e. Rusanov flux with (green diamonds) and without (red squares) a very small amount of implicit Rusanov-like stabilization, and the alternative \emph{upwind} flux (blue deltas). The additional stabilization is defined by $c_h=0.05$. All the three line-plots almost overlap: the only visible difference is the oscillatory behavior of the results obtained with Rusanov convective fluxes without implicit stabilization (red squares), see the right-traveling shock. (see colored version online) 
} 
\label{fig:rp3}
\end{center}
\end{figure}

\smallskip
\paragraph{RP3}  
The RP3 test presents the evolution of high Mach shocks: one  left- and one right-traveling magnetosonic (fast shocks) at Mach numbers $12$ and $14.9$, respectively. 
Then, a tangential discontinuity. For this test, the piecewise linear reconstruction was not used to improve the robustness of the simulation. 
 In this regard, for future applications the authors believe that adaptive limiting techniques should be used to stabilize high-Mach shocks without graving on the accuracy of the smooth regions of the flow. 
Instead, we compared the results obtained by choosing the two different alternative numerical fluxes, i.e. the Rusanov  (\ref{eq:Rusanovx})-(\ref{eq:Rusanovz}) with and without a small amount of Rusanov stabilization in the implicit pressure solver ($c_h=0,0.05$, see eq. \ref{eq:hdef}), and the \emph{upwind} (\ref{eq:Upwindx}) fluxes, see Fig. \ref{fig:rp3}. In particular, the artificial Rusanov stabilization in the implicit solver (obtained with $c_h=0.05$ in eq. \ref{eq:hdef}) as well the upwind flux  were two sufficient alternative ways to remove some oscillations. Those oscillations were due to time-accuracy, indeed also reducing the $\CFL$ parameter was also giving almost oscillation-free results.
Note that the \emph{upwind} flux gave a much better approximation of the central contact wave, see the zoom windows of the density plot in Fig. \ref{fig:rp3}. These results suggest that a better choice for the numerical flux could  significantly improve the quality of the numerical results.
Regarding the pressure plot, all the three line-plots almost overlap: the only visible difference is the oscillatory behavior of the results obtained with Rusanov convective fluxes without implicit stabilization (red squares), see the right-traveling shock in Fig. \ref{fig:rp3}.

\begin{figure}[!t]
\begin{center}
\includegraphics[width=0.48\textwidth]{\MyFigFolder/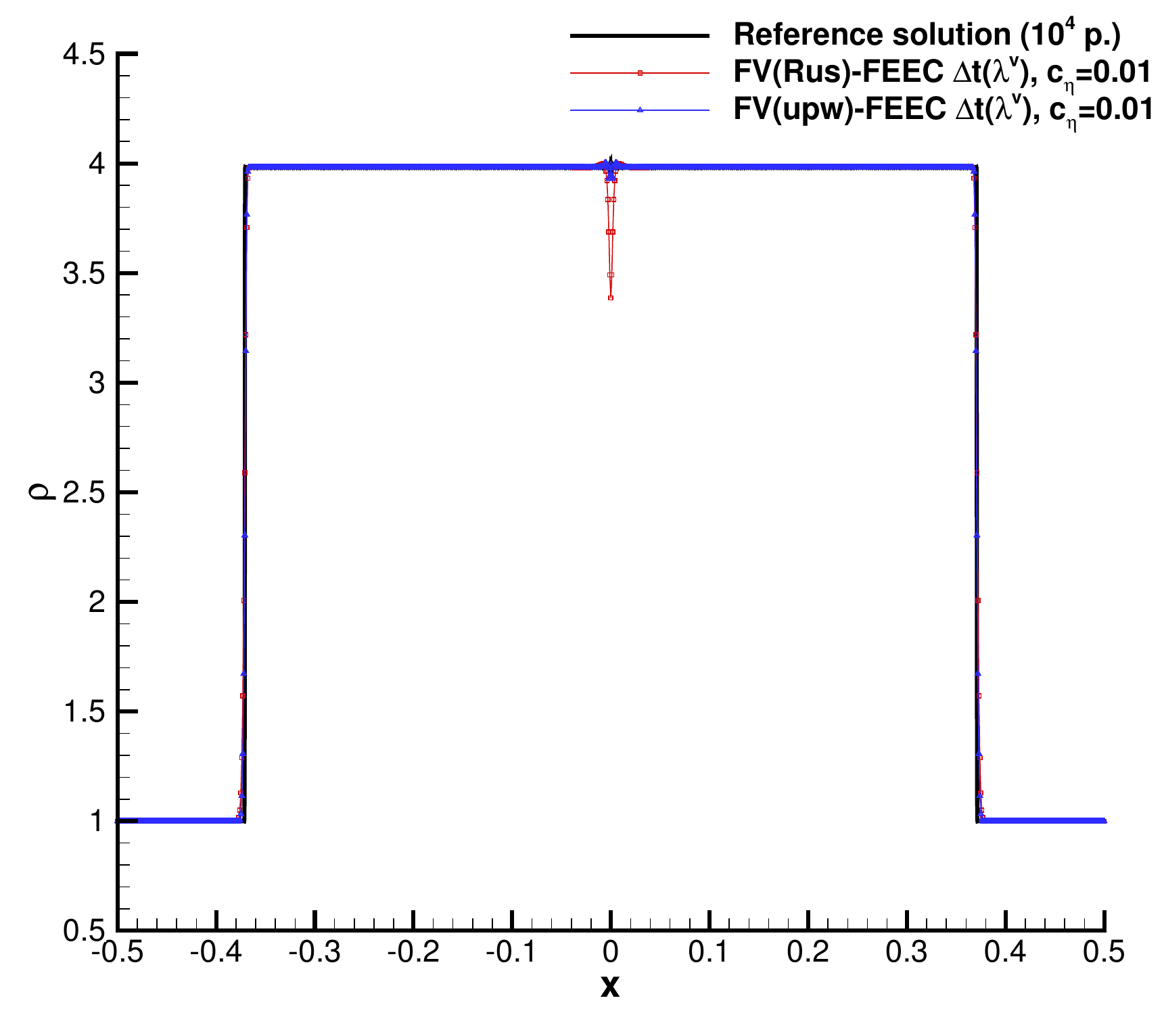}    
\includegraphics[width=0.48\textwidth]{\MyFigFolder/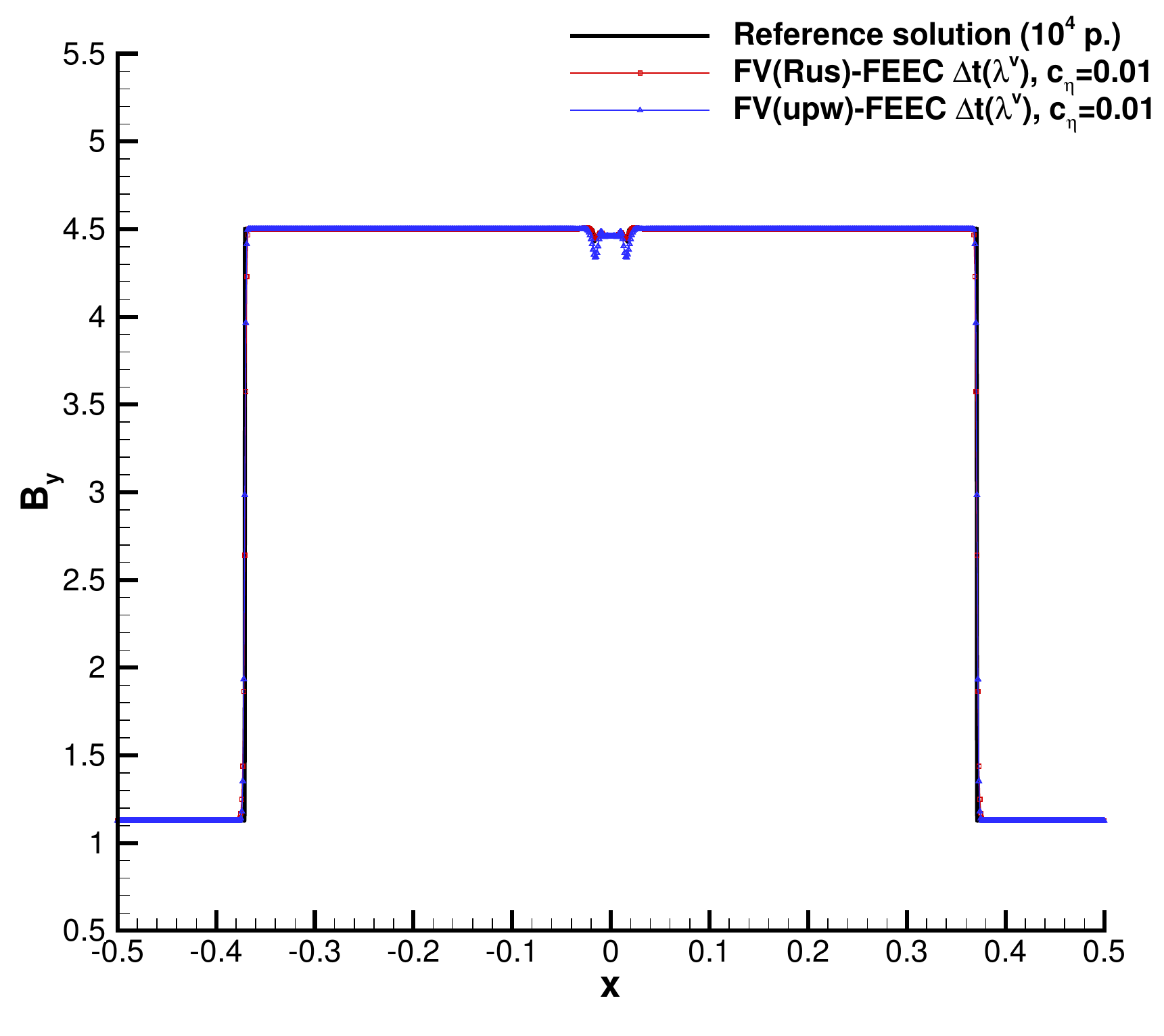}   
\caption{Reference and numerical solution 
 for the Riemann problem  RP4 at $t=0.03$. Density (left column) and  magnetic field component $B_y$ (right column) are plotted, comparing the numerical solutions obtained after choosing different numerical convective fluxes, i.e. Rusanov flux (red squares),  $\theta_b=0.55$ (green squares), and the alternative \emph{upwind} flux (blue deltas). A very small amount of artificial Rusanov-resistivity has been added by setting $c_\eta=0.01$. (see colored version online) 
} 
\label{fig:rp5}
\end{center}
\end{figure}

\smallskip
\paragraph{RP4 (Mach 25.5)}  
The RP4 is a test configuration that generates \textbf{two Mach $\mathbf{25.5}$ fast shocks}. Also in this extreme situation, the waves a correctly resolved by our hybrid FV-FEEC solver. Numerical results for the density and the vertical component of the magnetic field $B_y$ are plotted in Fig. \ref{fig:rp5} next to the reference solution.
A very small amount of numerical stabilization has been added by setting $c_\eta=0.01$, see eq. (\ref{eq:sxyz}). 
Also in this case, the alternative \emph{upwind} numerical flux seems to give higher-accurate results.
 

\subsection{2D problems}
\subsubsection{Convergence study: Alfv\'en wave test} 

To verify the order of accuracy of our semi-implicit hybrid FV-FEEC method, the low Mach Alfv\'en wave test is here solved by sequentially refining the grid. Indeed, the Alfv\'en wave test is a sinusoidal time-dependent solution of the MHD equations, and it is often used in the literature to measure the experimental order of accuracy of new numerical schemes.
The solution is three-dimensional but, for simplicity, it is reduced to 2d by setting 
the direction of propagation  to be on the $xy$-plane, but still non-aligned with the grid, i.e.
$$  \mathbf{n} = ( n_x, n_y, n_z ) = \left( 1, 2, 0 \right)/\sqrt{5}.$$
The solution is given by  
\begin{equation}
\left\{
  \begin{array}{rl}
\rho &= 1, \\
\v &= \alpha \left(  - n_y \cos \varphi,   n_x \cos \varphi, \sin \varphi  \right),  \\
p &= 10^2\\
\B &=   \left( n_x + n_y \alpha  \cos \varphi, n_y - n_x \alpha  \cos \varphi, -\alpha \sin \varphi  \right),
\end{array} \right. \label{eq:Alfven}
\end{equation} 
with
$$  \varphi = 2 \pi \frac{1}{n_y} \left[ n_x\left( x - n_x t\right) +  n_y \left( y - n_y t\right)\right].  $$
Starting from the initial time $t=0$, the simulation is run until the final time $t_f= \sqrt{5}/2$. Tables \ref{tab:Alfven}-\ref{tab:AlfvenHeli} collect all the convergence data obtained by our semi-implicit hybrid FV-FEEC scheme with and without using the orthogonality-preserving cross product (\ref{eq:crossp0}). The second order of accuracy is confirmed by the experimental data. For this test the security parameter is set to $\CFL=1/2$,  while the effective Courant number depends on the Alfv\'enic time-scales, i.e. $\Delta t = \Delta t( \lambda^{b})$. 
Then, a second-order \emph{non-limited} (non-TVD) MUSCL-Hancock reconstruction is used to evaluate the numerical fluxes of the finite-volume (convective) step, and a  Crank-Nicolson time discretization $\theta_p=\theta_b=1/2$ is solved for the FEEC (Alfv\'enic and acoustic) steps.  

For the purpose of validating the conservation of magnetic helicity, we ran again  the same test without the constant background magnetic field. In Figure \ref{fig:helicity}, the time evolution of numerical magnetic helicity  (see Appendix \ref{app:A}) obtained with and without the helicity preserving property is plotted: (i) using the orthogonality-preserving cross product (\ref{eq:crossp0}) and setting the number of iterations of the Alfv\'enic step $S_b=0$; (ii) using the orthogonality-preserving cross product (\ref{eq:crossp0}) and setting the number of iterations of the Alfv\'enic step $S_b=S_b(\epsilon)$; (iii) or using the cross product (\ref{eq:crossp0}) and setting $S_b=S_b(\epsilon)$. For this test, we used  
\begin{equation}
\label{eq:Sb}
S_b  = S_b(\epsilon) \quad \text{ such that} \quad \left\|\mathbf{B}^{S_b(\epsilon)+1}-\mathbf{B}^{S_b(\epsilon)} \right\|/ \left\|\mathbf{B}^{S_b(\epsilon)} \right\| \leq \epsilon \ll 1,
\end{equation}
 with $\epsilon=10^{-14}$. Again we have Crank-Nicolson in time $\theta_b=\frac{1}{2}$.  Results show that magnetic energy and helicity can be effectively better conserved when using a tolerance-dependent recursion in the Picard iterations, if combined with a proper orthogonality preserving cross product. On the other hand, the use of the orthogonality preserving cross product has no substantial effect in terms of conservation if only one iteration is performed ($S_b=0$), 
 see Figure \ref{fig:helicity}. We also plotted the results obtained by a purely explicit scheme obtained by setting $\theta_p=\theta_b=0$ and adding the corresponding stabilization terms that are typical of explicit schemes: in the Rusanov flux of the FV scheme, and in the resistivity step.

\renewcommand{\arraystretch}{0.8}
 \begin{table}[!t] 
 \centering
 \begin{tabular}{!{\extracolsep{-5pt}}ccccccccc!{}}
   \multicolumn{9}{c}{\textbf{Low Mach Alfv\'en-wave test}} \\
   \hline
& $N_{\text{element}}$ &  $L_1$ error &  $L_2$ error & $L_{\infty}$ error & $L_1$ or. & $L_2$ or. & $L_\infty$ or. &    Th. \\
   \hline
   \hline 
   \cline{2-8}
   \multirow{4}{*}{$v_x$} 
& $ 20^2$	& 1.383E-01	& 7.710E-02	& 6.009E-02	&  ---  &  ---	&  --- &  \multirow{4}{*}{\textbf{2}} \\
& $ 40^2$	& 3.445E-02	& 1.967E-02	& 1.535E-02	& 2.00	& 1.97	& 1.97 & 			\\ 
& $ 80^2$	& 8.447E-03	& 4.790E-03	& 3.602E-03	& 2.03	& 2.04	& 2.09 & 			\\ 
& $160^2$	& 2.146E-03	& 1.216E-03	& 9.144E-04	& 1.98	& 1.98	& 1.98 & 			\\ 
   \cline{2-8}
   \multirow{4}{*}{$v_y$} 
& $ 20^2$	& 7.303E-02	& 4.064E-02	& 3.211E-02	&  ---  &  ---	&  --- &  \multirow{4}{*}{\textbf{2}} \\
& $ 40^2$	& 1.459E-02	& 8.409E-03	& 6.766E-03	& 2.32	& 2.27	& 2.25 &  		\\ 
& $ 80^2$	& 3.949E-03	& 2.247E-03	& 1.717E-03	& 1.89	& 1.90	& 1.98 &  		\\ 
& $160^2$	& 9.059E-04	& 5.184E-04	& 3.962E-04	& 2.12	& 2.12	& 2.12 &  		\\ 
   \cline{2-8}
   \multirow{4}{*}{$B_x$} 
& $ 20^2$	& 5.000E-01	& 2.787E-01	& 2.218E-01	&  ---  &  ---	&  --- &  \multirow{4}{*}{\textbf{2}} \\
& $ 40^2$	& 1.193E-01	& 6.803E-02	& 5.321E-02	& 2.07	& 2.03	& 2.06 &  		\\ 
& $ 80^2$	& 2.969E-02	& 1.682E-02	& 1.258E-02	& 2.01	& 2.02	& 2.08 &  		\\ 
& $160^2$	& 7.388E-03	& 4.192E-03	& 3.181E-03	& 2.01	& 2.00	& 1.98 &  		\\ 
   \cline{2-8}
   \multirow{4}{*}{$B_y$} 
& $ 20^2$	& 2.499E-01	& 1.399E-01	& 1.103E-01	&  ---  &  ---	&  --- &  \multirow{4}{*}{\textbf{2}} \\
& $ 40^2$	& 5.946E-02	& 3.400E-02	& 2.688E-02	& 2.07	& 2.04	& 2.04 & 		\\ 
& $ 80^2$	& 1.479E-02	& 8.397E-03	& 6.375E-03	& 2.01	& 2.02	& 2.08 & 		\\ 
& $160^2$	& 3.683E-03	& 2.092E-03	& 1.565E-03	& 2.01	& 2.01	& 2.03 & 		\\  
   \hline
 \end{tabular}
\caption{$L_1$, $L_2$ and $L_\infty$ errors and convergence rates for the low Mach Alfv\'en wave test obtained with the \lq\lq default\rq\rq\, cross product.} \label{tab:Alfven}
 \end{table}

\renewcommand{\arraystretch}{1.0}

\renewcommand{\arraystretch}{0.8}
 \begin{table}[!t] 
 \centering
 \begin{tabular}{!{\extracolsep{-5pt}}ccccccccc!{}}
   \multicolumn{9}{c}{\textbf{Low Mach Alfv\'en-wave test} (orth.-pres. cross product)} \\
   \hline
& $N_{\text{element}}$ &  $L_1$ error &  $L_2$ error & $L_{\infty}$ error & $L_1$ or. & $L_2$ or. & $L_\infty$ or. &    Th. \\ 
   \cline{2-8}
   \multirow{4}{*}{$v_x$}  
& $ 20^2$	& 1.372E-01	&  7.656E-02	&  5.636E-02	&  ---	 &  ---	  &  ---	& \multirow{4}{*}{\textbf{2}} \\ 
& $ 40^2$	& 3.452E-02	&  1.927E-02	&  1.410E-02	&  1.99	 &  1.99	&  2.00 &					\\ 
& $ 80^2$	& 8.443E-03	&  4.677E-03	&  3.316E-03	&  2.03	 &  2.04	&  2.09 &					\\ 
& $160^2$	& 2.144E-03	&  1.188E-03	&  8.309E-04	&  1.98	 &  1.98	&  2.00 &					\\ 
   \cline{2-8}                                                               
   \multirow{4}{*}{$v_y$}                                                    
& $ 20^2$	& 7.289E-02	&  4.062E-02	&  3.005E-02	&  ---	 &  ---	  &  ---  & \multirow{4}{*}{\textbf{2}} \\
& $ 40^2$	& 1.462E-02	&  8.174E-03	&  6.068E-03	&  2.32	 &  2.31	&  2.31 & 				\\ 
& $ 80^2$	& 3.950E-03	&  2.188E-03	&  1.561E-03	&  1.89	 &  1.90	&  1.96 & 				\\ 
& $160^2$	& 9.063E-04	&  5.022E-04	&  3.557E-04	&  2.12	 &  2.12	&  2.13 & 				\\ 
   \cline{2-8}                                                               
   \multirow{4}{*}{$B_x$}                                                    
& $ 20^2$	& 4.980E-01	&  2.767E-01	&  2.028E-01	&  ---	 &  ---	  &  ---  & \multirow{4}{*}{\textbf{2}} \\
& $ 40^2$	& 1.195E-01	&  6.664E-02	&  4.873E-02	&  2.06	 &  2.05	&  2.06 & 				\\ 
& $ 80^2$	& 2.967E-02	&  1.643E-02	&  1.160E-02	&  2.01	 &  2.02	&  2.07 & 				\\ 
& $160^2$	& 7.384E-03	&  4.090E-03	&  2.865E-03	&  2.01	 &  2.01	&  2.02 & 				\\ 
   \cline{2-8}                                                               
   \multirow{4}{*}{$B_y$}                                                    
& $ 20^2$	& 2.500E-01	&  1.388E-01	&  1.013E-01	&  ---	 &  ---	  &  ---  & \multirow{4}{*}{\textbf{2}} \\
& $ 40^2$	& 5.961E-02	&  3.329E-02	&  2.446E-02	&  2.07	 &  2.06	&  2.05 & 				\\ 
& $ 80^2$	& 1.480E-02	&  8.196E-03	&  5.841E-03	&  2.01	 &  2.02	&  2.07 & 				\\ 
& $160^2$	& 3.683E-03	&  2.040E-03	&  1.436E-03	&  2.01	 &  2.01	&  2.02 & 				\\ 
   \hline
 \end{tabular}
\caption{$L_1$, $L_2$ and $L_\infty$ errors and convergence rates for the low Mach Alfv\'en wave test obtained with the  orthogonality-preserving cross product.} \label{tab:AlfvenHeli}
 \end{table}

\renewcommand{\arraystretch}{1.0}
\begin{figure} 
\centering 
			\includegraphics[width=0.49\textwidth]{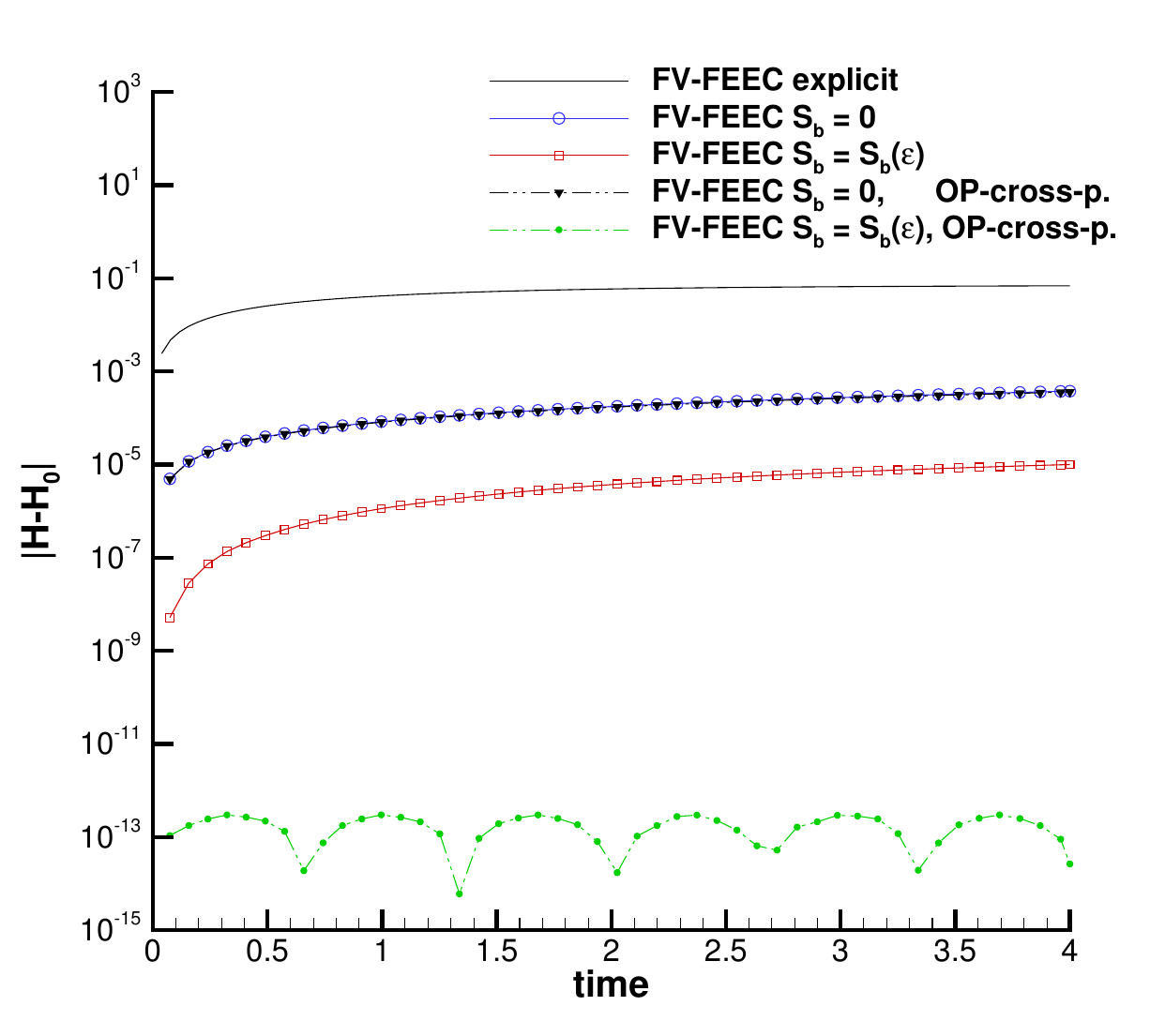} 
			\includegraphics[width=0.49\textwidth]{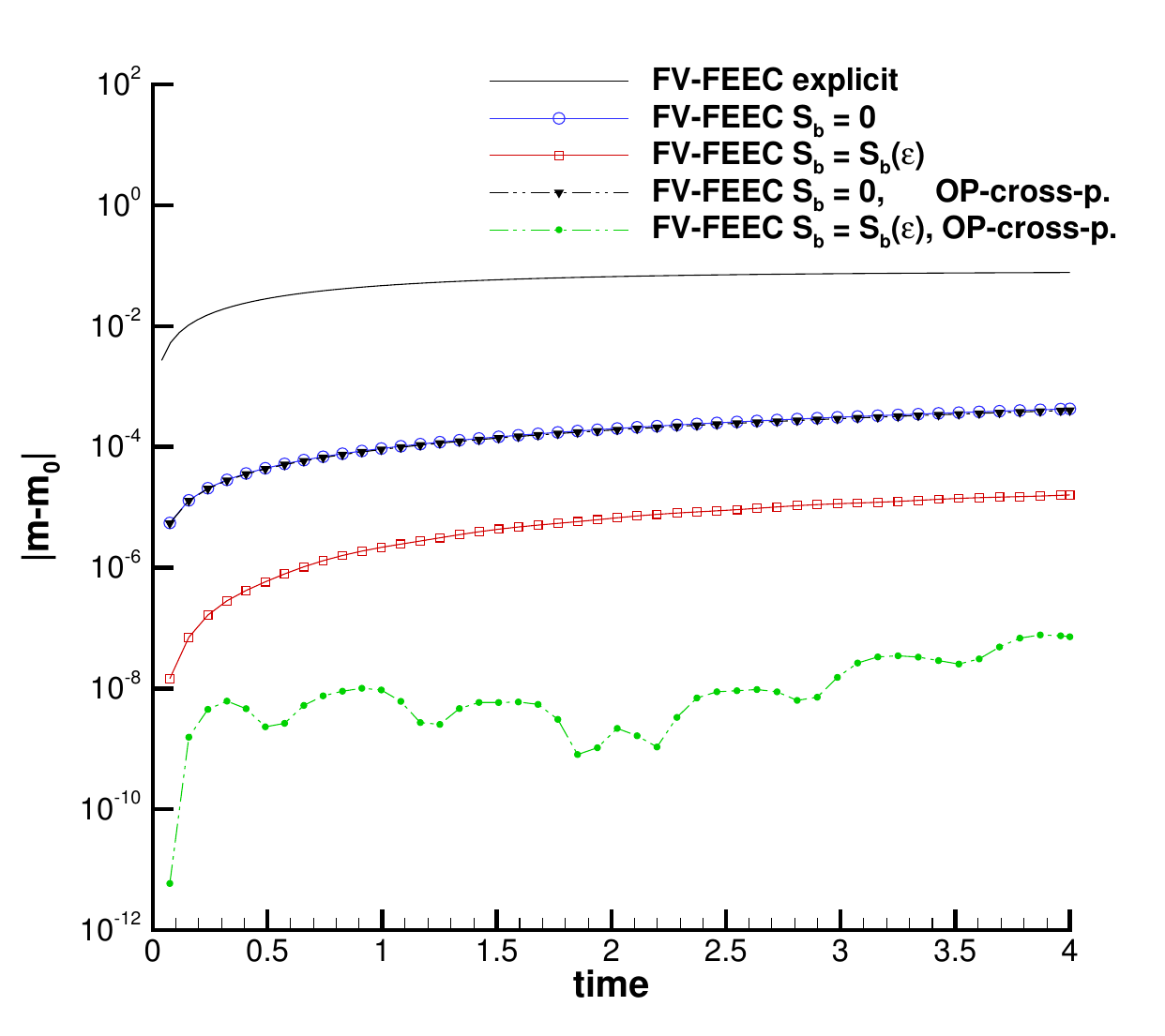}    
\caption{Comparison of the time evolution of the variation of the numerical magnetic helicity $|H-H_0|$ (left) and magnetic energy $|m-m_0|$ (right) with respect to the initial state for the Alfv\'en wave test (with zero back-ground magnetic field) computed on a 2d domain with our hybrid FV-FEEC scheme. Here, we compare the conservation properties of four different numerical setups: by defining $S_b=S_b(\epsilon)$ or $S_b=0$; by using the orthogonality preserving (OP-cross-p.) cross product or the \lq\lq default\rq\rq\, option; or by integrating the PDE explicitly in time} \label{fig:helicity}
\end{figure}

\subsubsection{Stationary MHD vortex} 
\begin{figure} 
\centering 
\begin{tabular}{c}
			\includegraphics[width=0.32\textwidth]{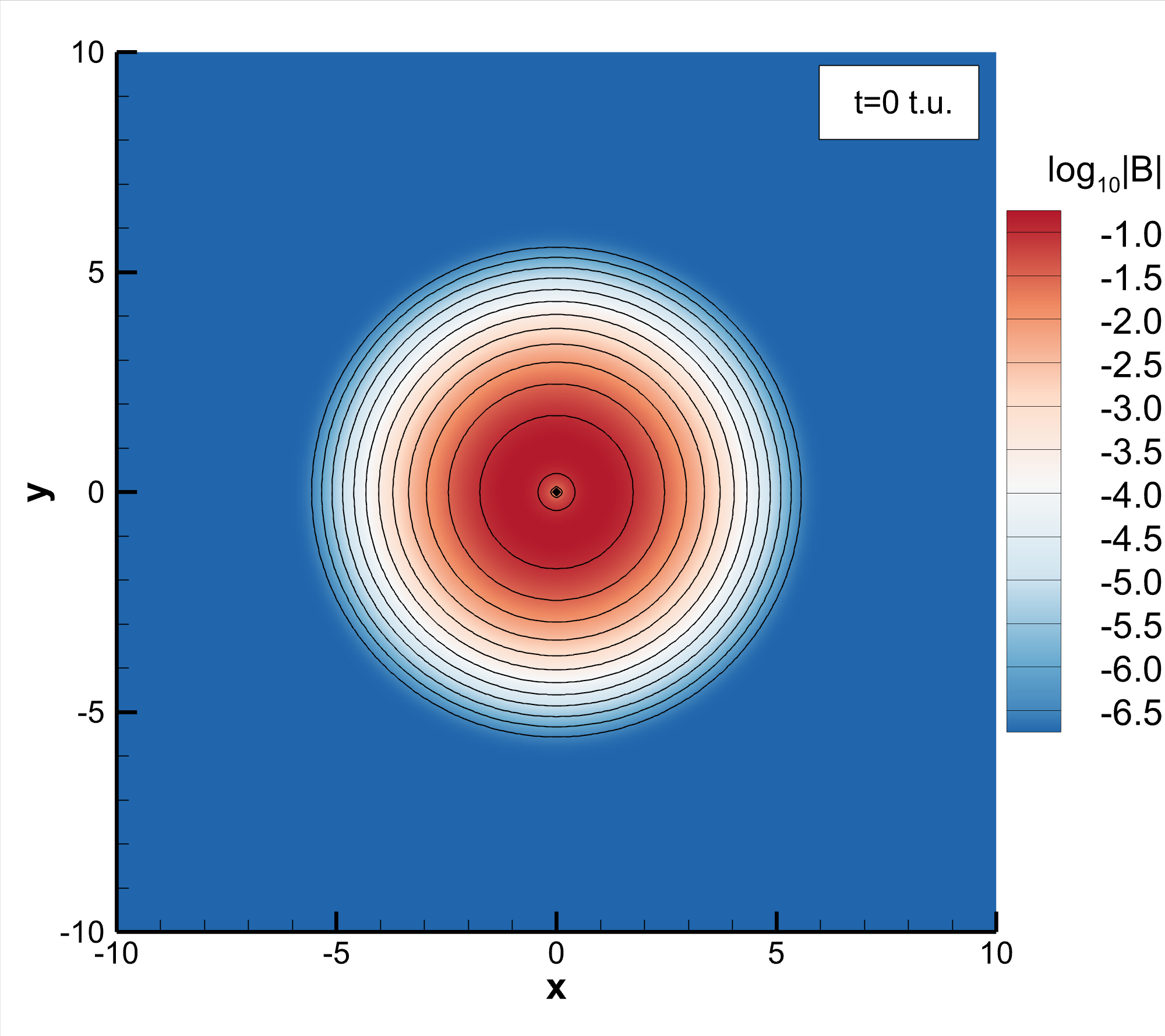} 
			\includegraphics[width=0.32\textwidth]{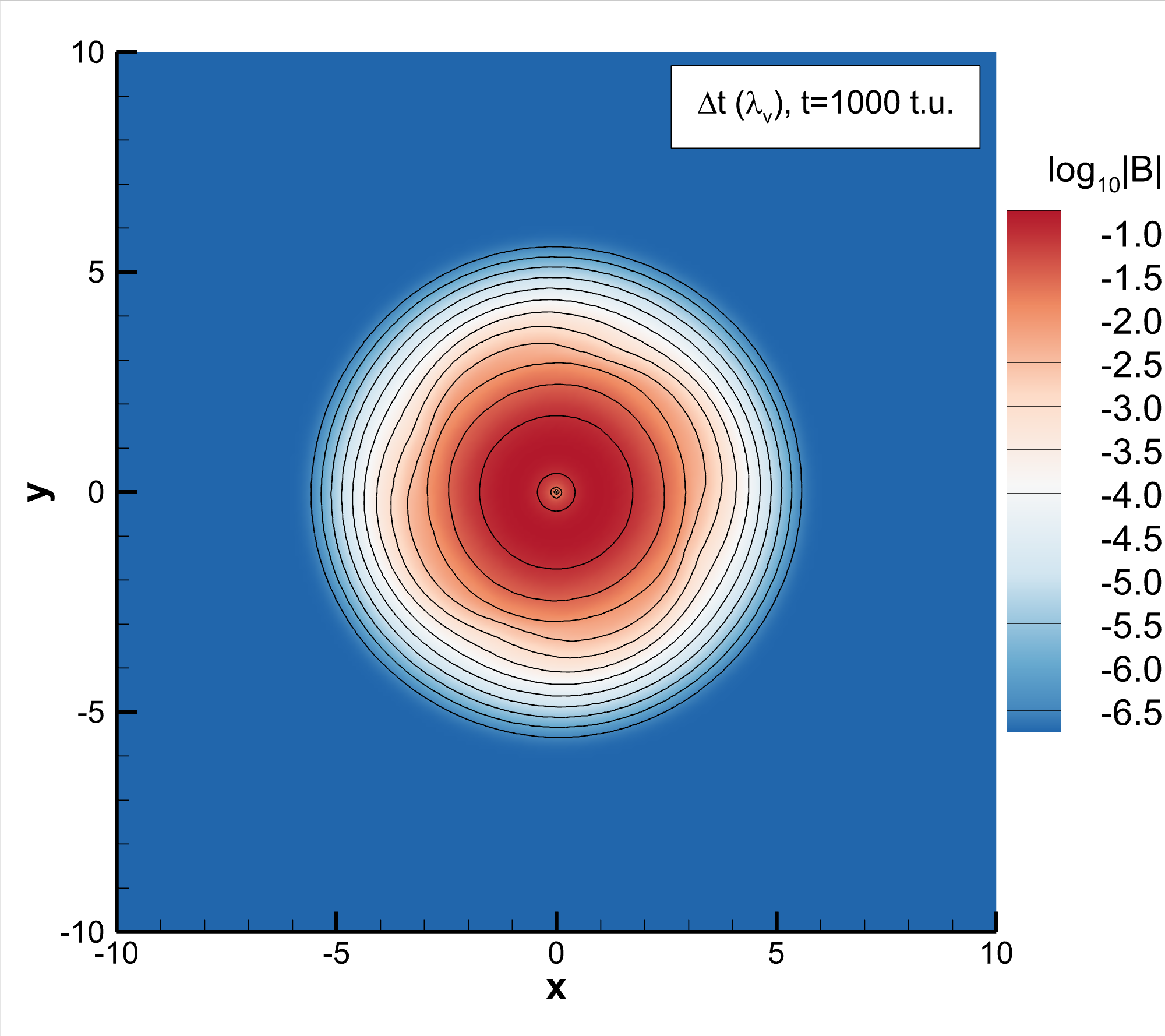} \\ 
			\includegraphics[width=0.32\textwidth]{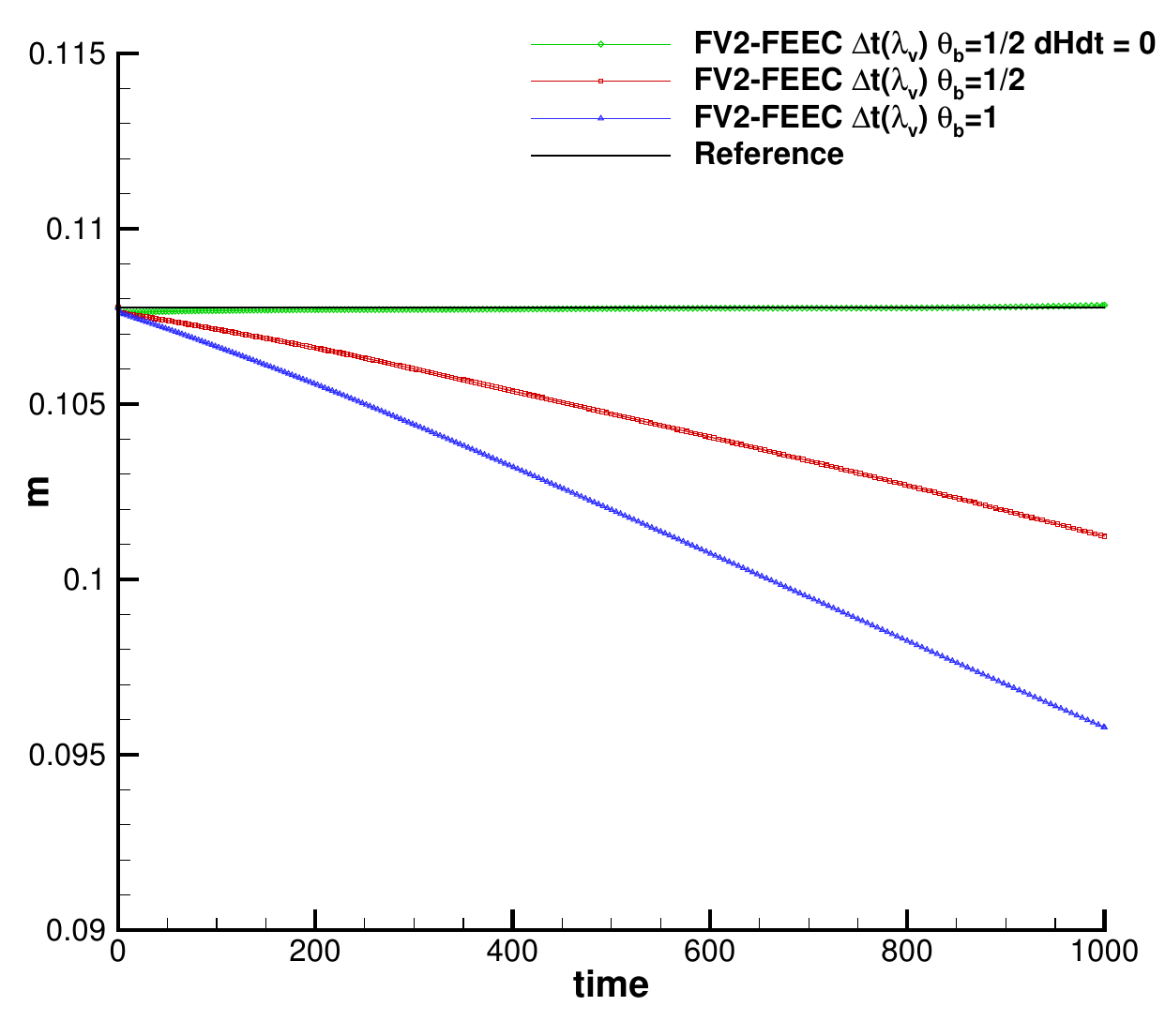} 
			\includegraphics[width=0.32\textwidth]{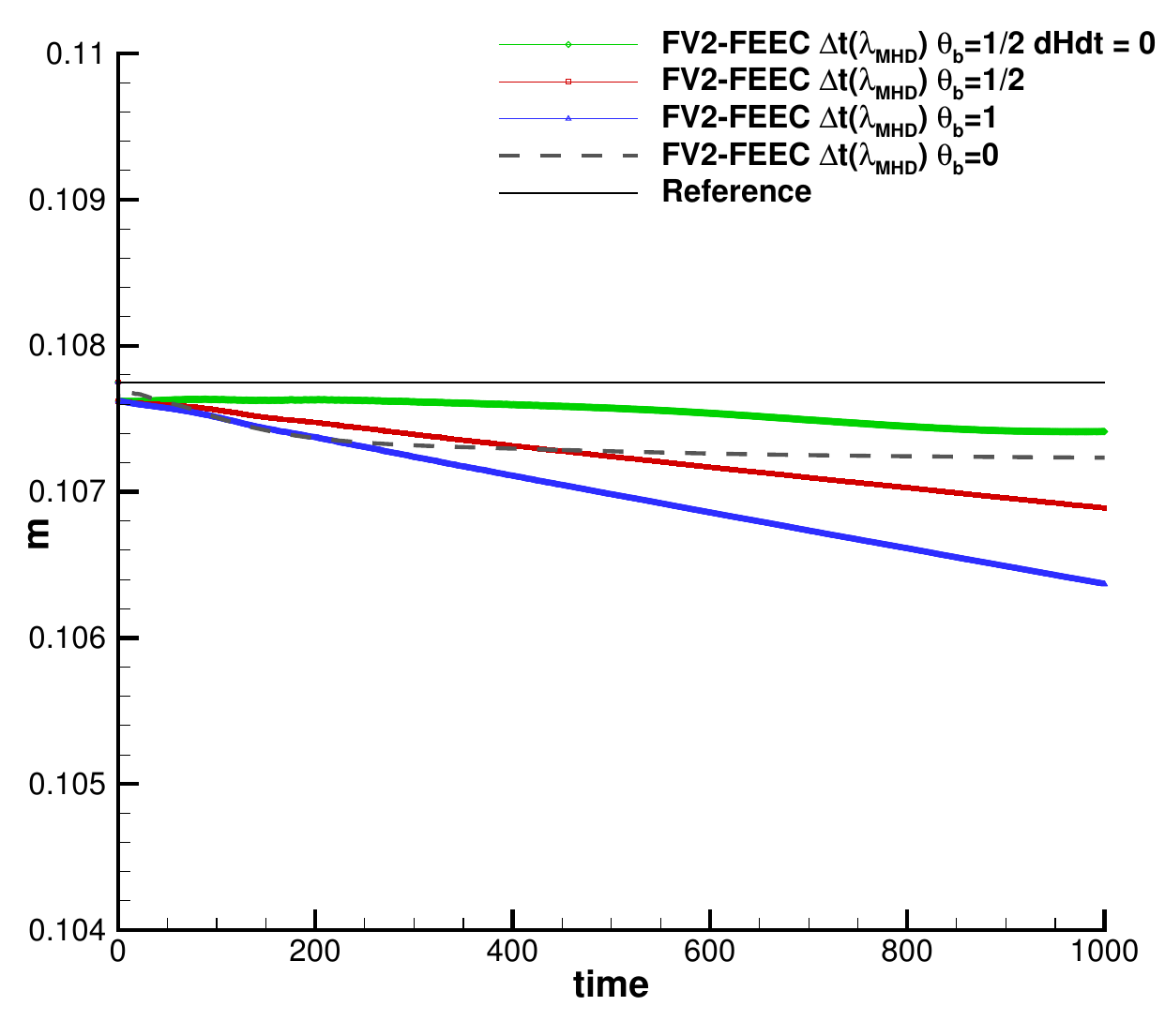} 
			\includegraphics[width=0.32\textwidth]{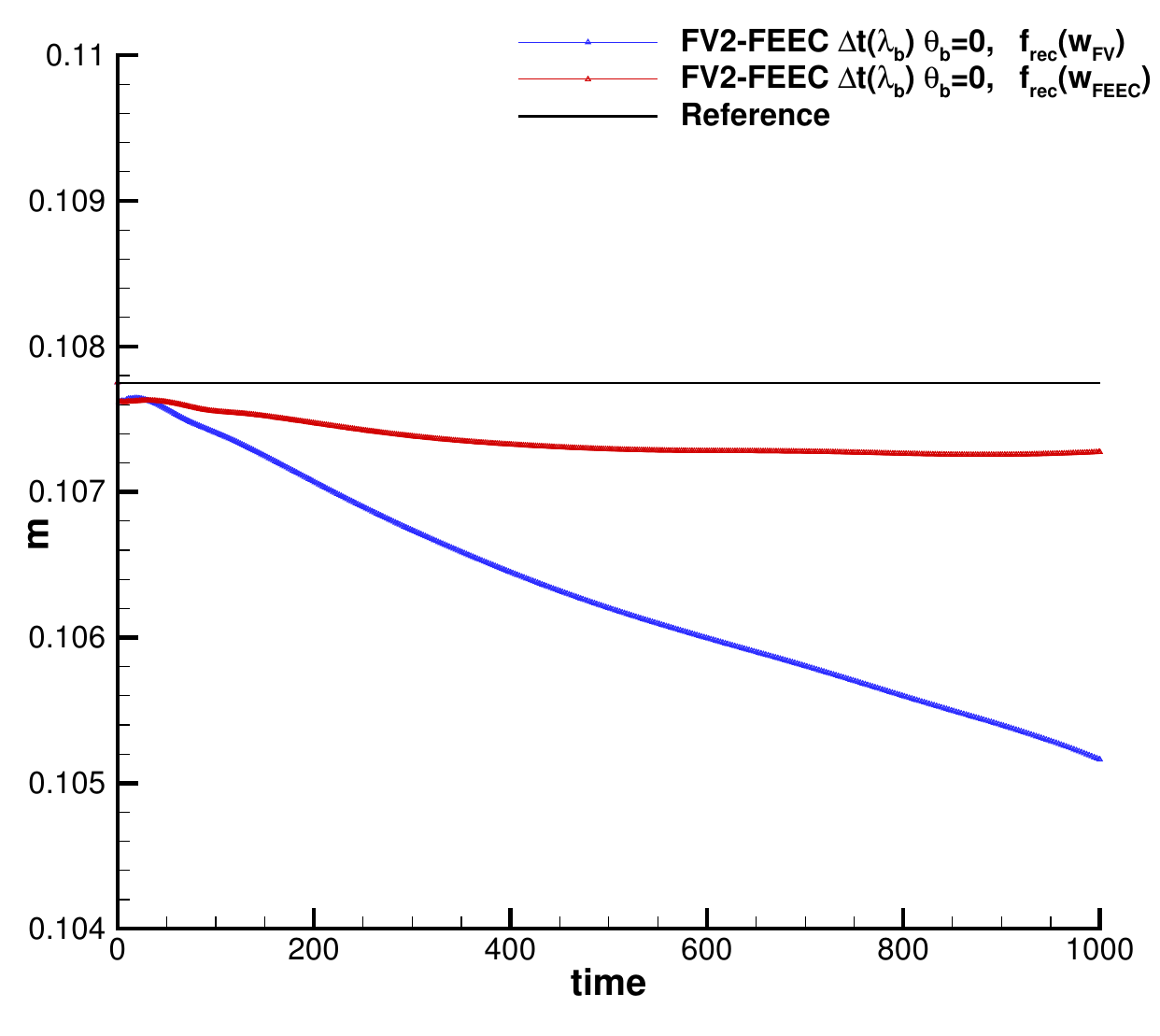}  
			\end{tabular}
\caption{Numerical results for the stationary isodensity MHD vortex problem. In the first row, the contour plot of $\log_{10}(|\B|)$  are plotted for the initial condition at $t=0$ t.u (left), and the final time $t=1000$ t.u. (right). The logarithm is plotted to highlight the main changes. \\ In the second row,  the time evolution of the magnetic energy $m$ computed by using different set-up is plotted: first (left), by choosing discrete time-steps scale of $\Delta t(\lambda_v)$ with different implicit weights $\theta_b=1,1/2$ and the helicity conserving discretization  (labeled in the line-legend with \lq\lq $dHdt=0$\rq\rq\,); second (center), the MHD time scale $\Delta t(\lambda_{\MHD})$ is used, and we also compare the fully explicit case $\theta_b=1,1/2,0$; third (right), the same time discretization with $\Delta t(\lambda_b)$, but using a standard TVD reconstruction, or using directly the FEEC polynomials to evaluate the numerical finite-volume fluxes.    
} \label{fig:IV}
\end{figure}

Here, we present the results of our hybrid FV-FEEC scheme against a long-time simulation of the stationary two-dimensional stationary isodensity vortex. In this test, at the initial time, denoting by $r^2=x^2+y^2$, the primitive variables are given by
\begin{equation} 
\left\{
  \begin{array}{rl}
	\rho & = 1,  												\\ 
	\v   & = \frac{v_0}{2\pi} e^{\frac{1}{2} (1-r^2)} \left( -y,x,0\right),   \\ 
	p & = \frac{1}{8\pi}\left(\frac{A_0}{2\pi}\right)^2 \left(1-r^2\right) e^{\left(1-r^2\right)} - \frac{1}{2} \left(\frac{v_0}{2\pi}\right)^2  e^{\left(1-r^2\right)}, \\ 
	\B & = \frac{A_0}{2\pi} e^{\frac{1}{2} (1-r^2)} \left( -y,x,0\right)
  \end{array}
	\right.
\end{equation}
corresponding to a vector potential 
$$
\mathbf{A}=(0,0, \frac{A_0}{2\pi} e^{\frac{1}{2}(1-r^2)}),
$$
 see \cite{Balsara2004} for other details. The simulation is complete when the final time   $t_f=1000$ t.u. is reached. 
Remember that, in this work, \emph{well-balancing} techniques, see e.g. \cite{hybridhexa}, are not implemented. Then, this kind of equilibrium is \emph{not} expected to be \emph{exactly} preserved and the conservation properties of different numerical set-up can be compared.
In particular, different time-discretization are tested by choosing between time-step sizes of the \emph{convective} time-scale $\Delta t(\lambda_v)\sim 0.158$ t.u., the Alfv\'enic time-scale $\Delta t(\lambda_b)\sim 0.080$ and the full MHD scale  $\Delta t(\lambda_{\MHD})\sim 0.017$ t.u.. Moreover, we also compare different choices of the implicit weight parameter $\theta_b=1,1/2,0$, i.e. we also tested the fully  explicit (Alfv\'enic) case. Some selected interesting data are plotted in Fig. \ref{fig:IV}. It is important to underline the fact that, even for the explicit case ($\theta_b=0$) the simulation was stable even \emph{without additional stabilization terms} in the Alfv\'enic system, i.e. $c_\eta=0$. For the Crank-Nicolson $\theta_b=1/2$, we also compared the results by activating the helicity conserving strategy, i.e. using the orthogonality preserving product and waiting for convergence of the solution of the nonlinear Alfv\'enic system.

First, the higher conservation capabilities of the  magnetic-helicity and -energy conserving options are experimentally confirmed, see the left one-dimensional plot in Fig. \ref{fig:IV}. 
Obviously, after refining the time-step size, the time-resolution is higher and the differences becomes smaller and smaller, see the central one-dimensional plot in Fig. \ref{fig:IV}. For this test we also compared different choices for the polynomial reconstruction, i.e. between the FEEC reconstruction (\ref{eq:FEECrec}) and the finite-volume reconstruction (\ref{eq:MUSCLrec}), see the right one-dimensional plot in Fig. \ref{fig:IV}. 
As expected, the FEEC polynomials are \lq\lq smoother\rq\rq\,, preserve more structures and have better conservation properties, while the chosen finite-volume reconstruction is confirmed to be more diffusive.

\subsection{Ideal Orszag-Tang vortex} \label{sec:OT}

\begin{figure} 
\centering 
 \begin{tabular}{ccc}   	 
 	\begin{adjustbox}{angle=90} \text{\;\;\;\;\;FV-FEEC}\end{adjustbox}&	\includegraphics[width=0.45\textwidth]{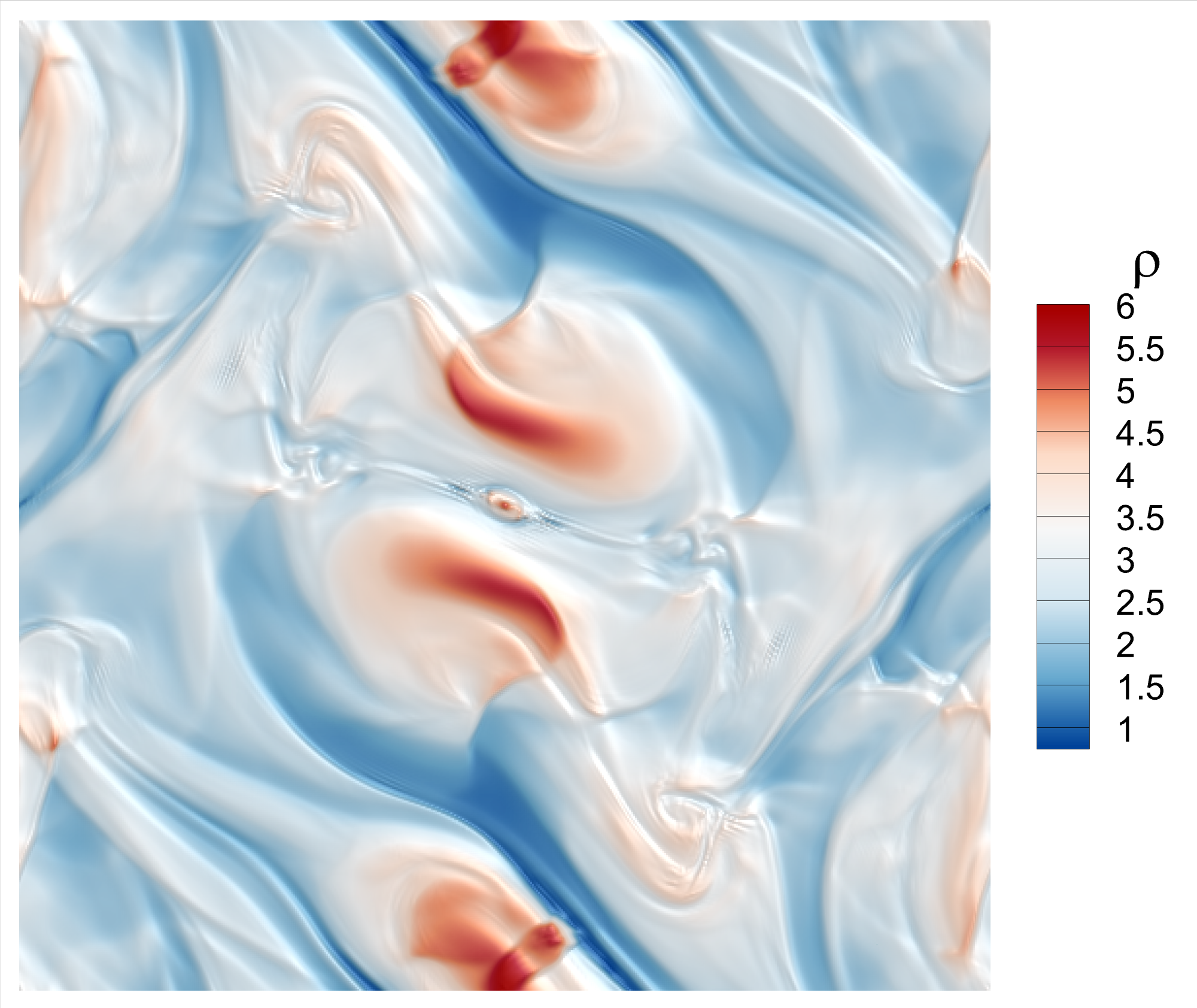} &  
 		\includegraphics[width=0.45\textwidth]{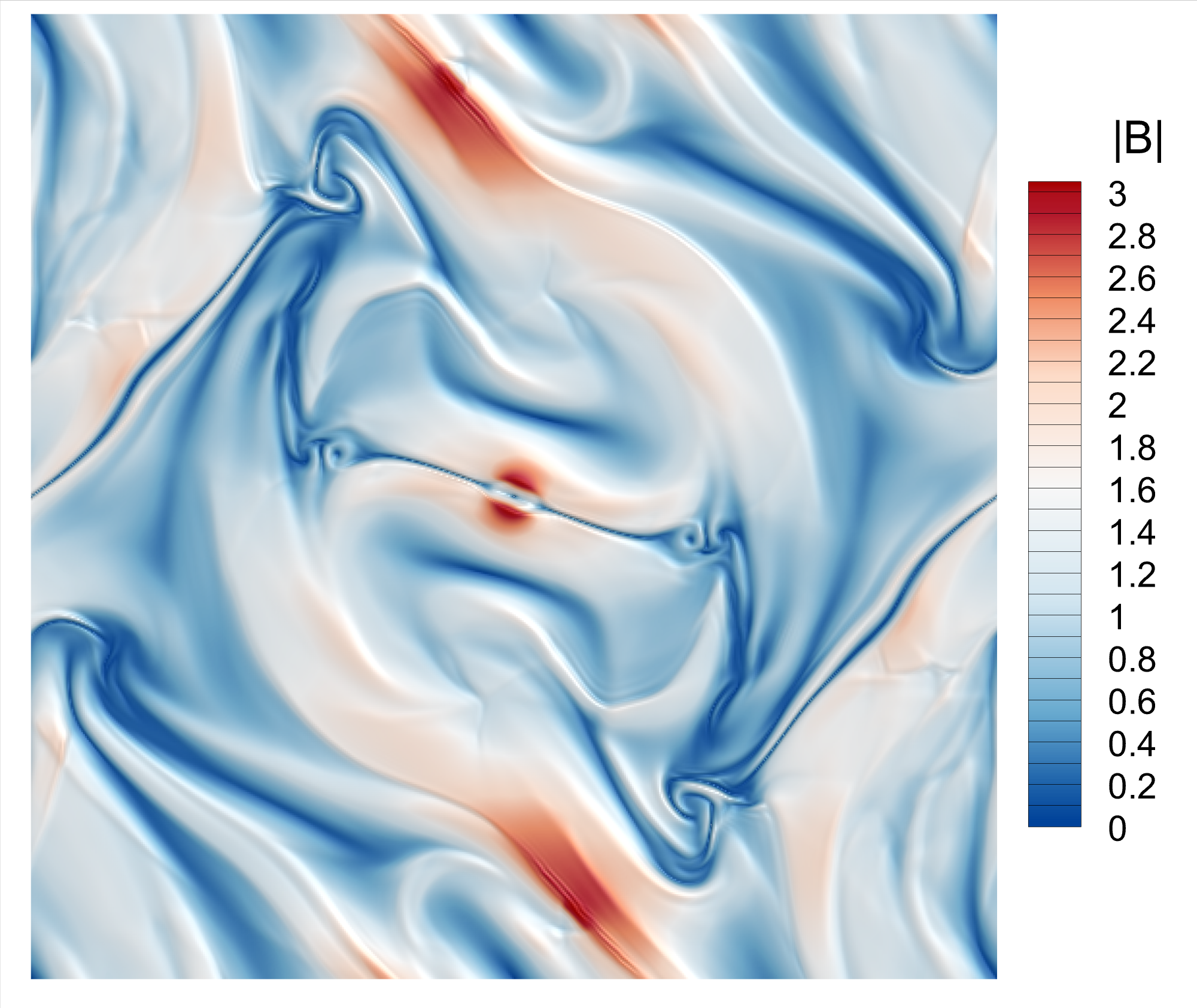}  \\  
 	\begin{adjustbox}{angle=90}\centering \text{reference (ADER-DG)}\end{adjustbox} &
 \includegraphics[width=0.45\textwidth]{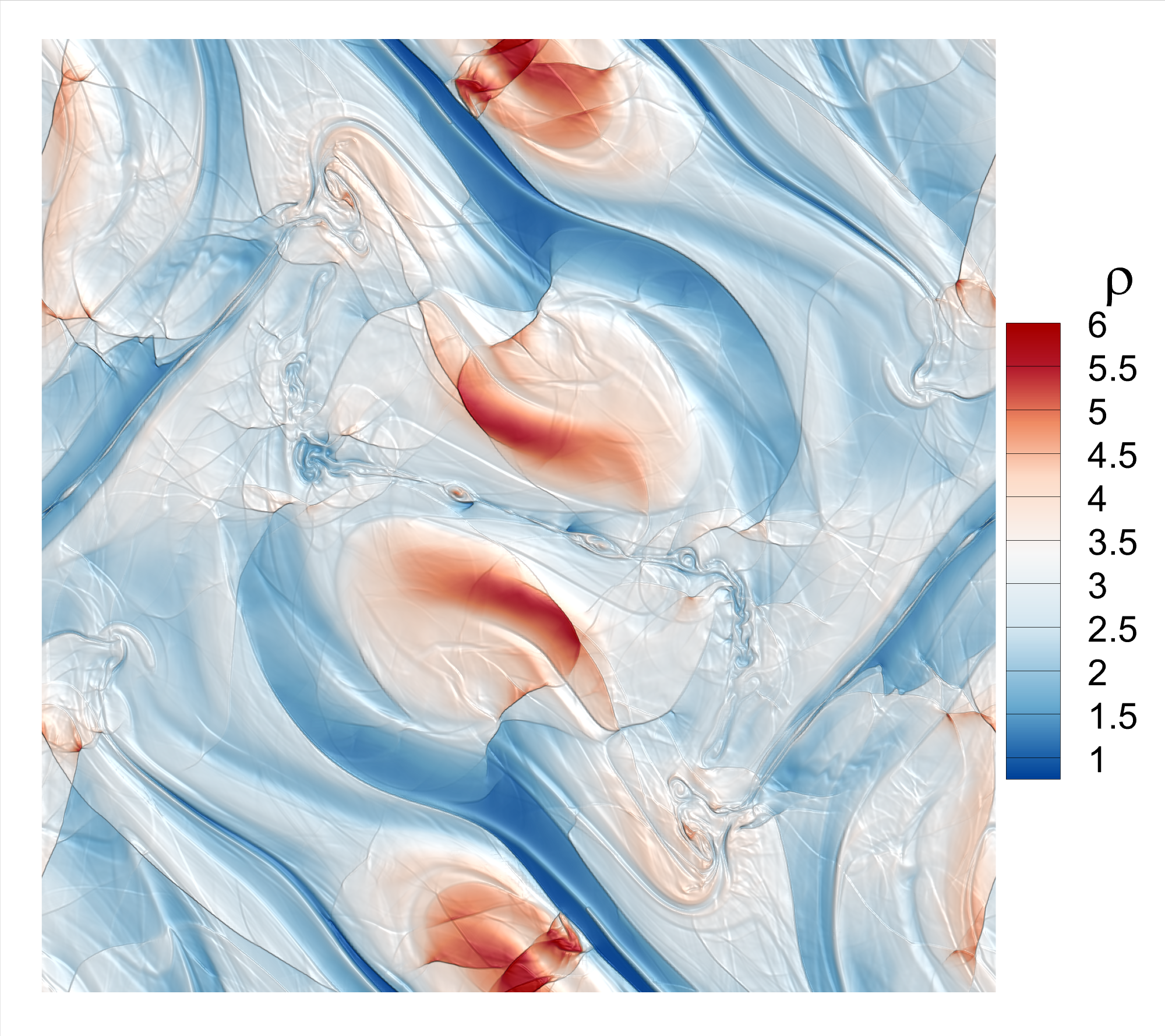} &  
 \includegraphics[width=0.45\textwidth]{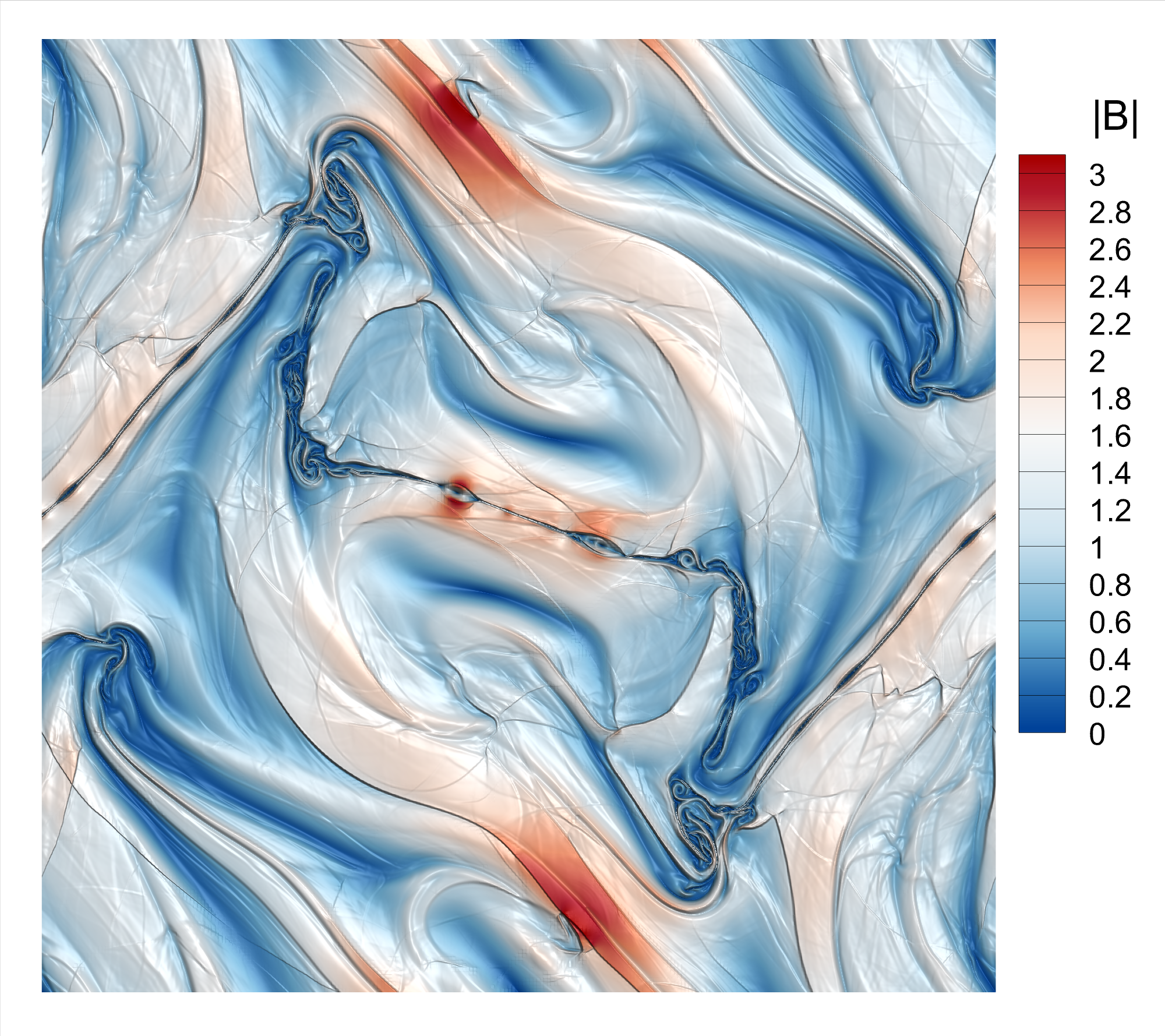} 
 	\end{tabular}
\caption{Numerical results for the ideal Orszag-Tang vortex system at time $t=5$ t.u. obtained with our semi-implicit hybrid FV-FEEC method on a $\Delta x = \Delta y= 1/500$ grid (first row), compared with a reference solution computed at the aid of a sixth order accurate ADER-DG scheme with sub-cell FV-WENO limiting (second row), see \cite{Zanotti2015c}. The discrete solutions for the fluid density $\rho$ (left) and magnetic field $|\B|$ (center) are plotted.} \label{fig:OT}
\end{figure}

\begin{figure} 
\centering    
 		\includegraphics[width=0.32\textwidth]{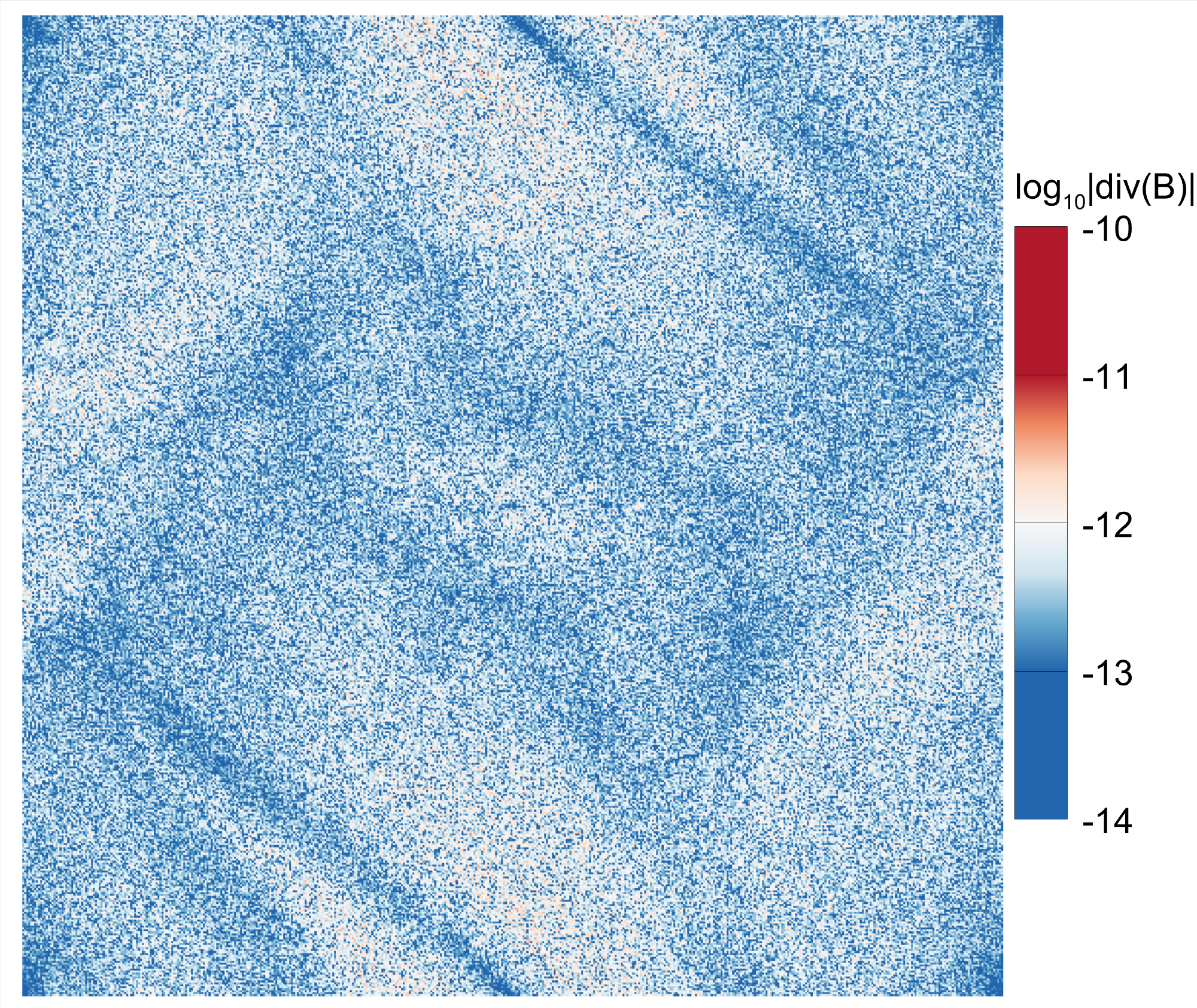}  
 		\includegraphics[width=0.32\textwidth]{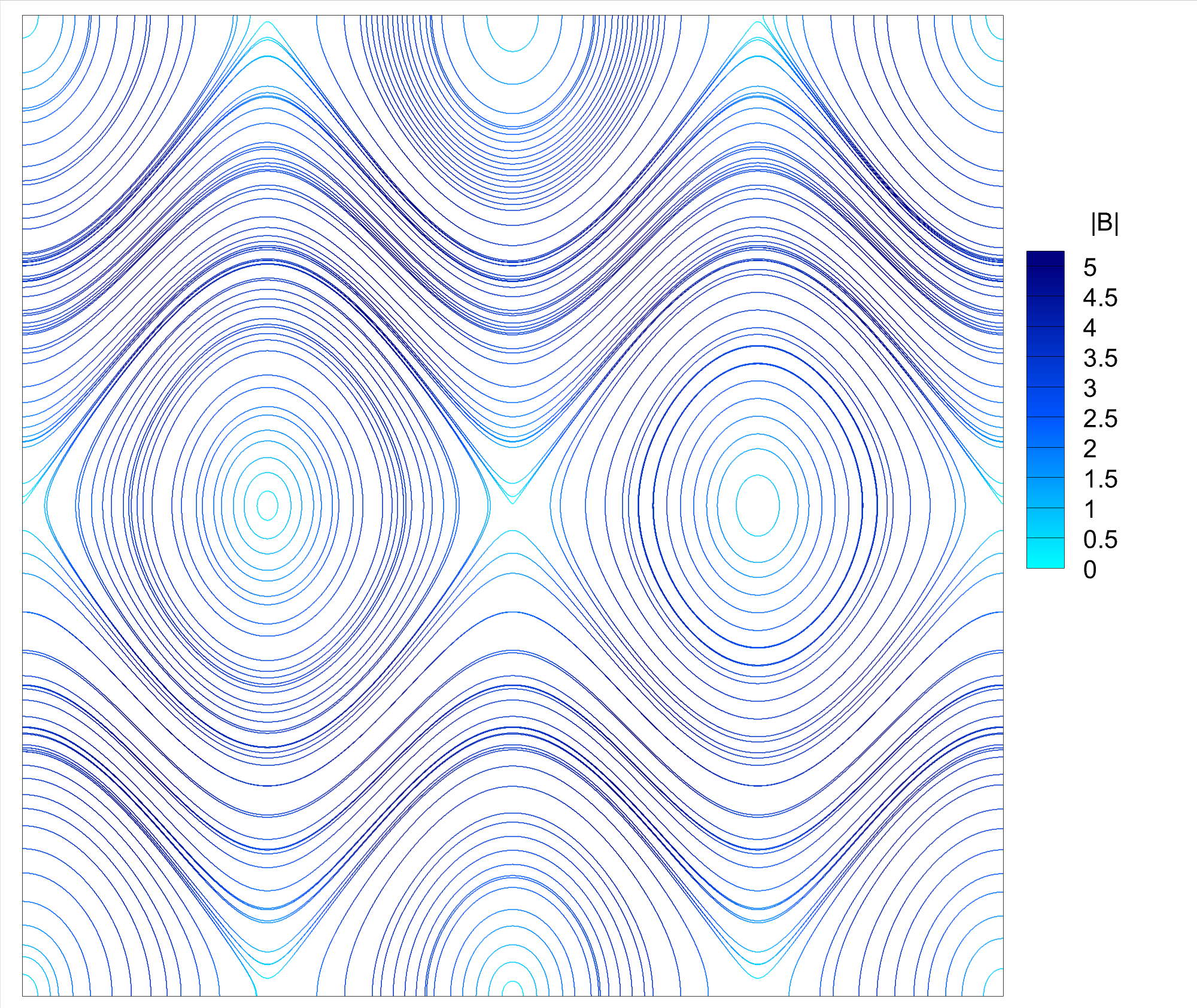} \includegraphics[width=0.32\textwidth]{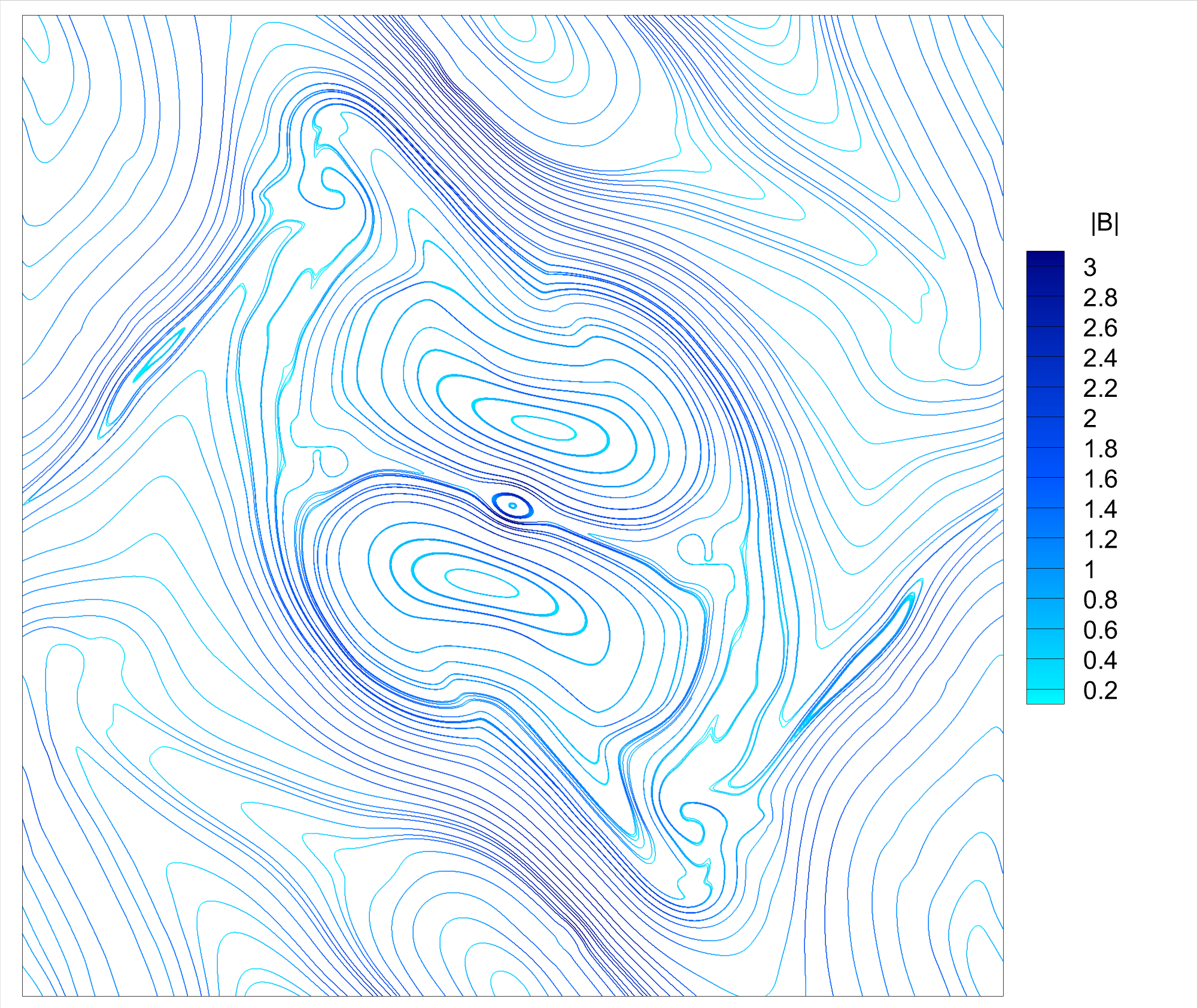} 
\caption{From left to right: numerical divergence error $\log_{10}|div(\B)|$ at time $t=5$ and magnetic field streamlines for the ideal Orszag-Tang vortex system at times $t=0,5$ t.u. obtained with our semi-implicit hybrid FV-FEEC method on a $\Delta x = \Delta y= 1/500$ grid. Streamlines are colored by the magnetic field magnitude.} \label{fig:OTdivB}
\end{figure}

Another test to verify the robustness and accuracy of our novel semi-implicit hybrid FV-FEEC method against shock problems in two-spatial dimension is the nonlinear decay of the ideal Orszag-Tang vortex system, first proposed by \cite{OrszagTang}. For an analysis of the flow, we refer the reader to   \cite{PiconeDahlburg} and \cite{DahlburgPicone}, while the vortex configuration corresponds to the one chosen in 
 \cite{JiangWu}. The initial condition is given by 
\begin{equation} \label{eq:OT}
\left\{
  \begin{array}{rl}
	\rho & =   \gamma^2,  												\\ 
	\v   & = \left(  - \sin\left(y\right), \sin \left(x \right), 0 \right),    \\ 
	p & =   \gamma, \\ 
	\B & = \left(  - \sin\left(y\right), \sin \left(2x \right), 0\right),
  \end{array}
	\right.
\end{equation}
and, with $\gamma = 5/3$ this corresponds to having an MHD-\emph{beta}  $\beta = 2 p_0 /|B_0|^2= 10/3$. Starting from the smooth solution (\ref{eq:OT}), shock waves are generated and only robust shock-capturing methods can reach long-time simulations.
In Fig. \ref{fig:OT}, the numerical results obtained with our semi-implicit hybrid FV-FEEC scheme are compared to a reference solution, which is computed with the aid of a sixth order accurate ADER-DG scheme with sub-cell FV-WENO limiting (second row), see \cite{Zanotti2015c}. Note that the simulation is purely inviscid and   mesh convergence cannot be reached, but all the main patterns are well reproduced and compatible also with other published results in the literature. For this two-dimensional and shock-dominated test, divergence errors are still very low, of the order $10^{-11}/10^{-12}$, see   Fig. \ref{fig:OTdivB}. Finally, in Fig. \ref{fig:OTdivB} also the streamlines of the magnetic field of initial and final times are plotted. In principle, divergence error would cause magnetic streamlines to cross, while numerical viscosity (or better resistivity) would facilitate magnetic reconnection, and for purely inviscid simulations the number of magnetic island should remain constant. 

\subsubsection{MHD Rotor test}\label{sec:MHDRotor}
The so-called MHD rotor problem was proposed by \cite{BalsaraSpicer1999} for testing the numerical approximation of strong torsional Alfv\'en waves in a two-dimensional problem for ideal MHD. The so-called \lq\lq rotor\rq\rq\, consists in a dense cylinder, which is rapidly spinning at an initial time, immersed within a lower density background fluid and initially constant magnetic field. Then, the initial state is 
\begin{equation} \label{eq:MHDRotor}
\left\{
  \begin{array}{rl}
	 (\rho, \v ) &= \left\{ \begin{array}{ll}
						\left(10, \omega \times \mathbf{x} \right) & \text{if} \  |\mathbf{x}|<R \,; \\
						\left(1, (0,0,0) \right)  &  \text{otherwise}\,, \end{array}
\right.\\
	p & =   1, \\ 
	\B & = \left(  2.5,0, 0\right),
  \end{array}
	\right.
\end{equation} 
The rotor spin angular velocity is $\omega = 10 \hat{\mathbf{z}}$.
In this work, a \emph{narrow} linear tapering is used to smooth the radial discontinuity. Following \cite{Zanotti2015c,Fambri20}, we choose a \emph{narrow} tapering in the range $|\mathbf{x}|\in[R,R+5\%R]$ so that we may provide a reference ADER-DG solution.  The numerical solution may slightly differ with some other published results in the literature because many references use instead a \emph{wider} linear tapering in the range $|\mathbf{x}|\in[R,R+15\%R]$. In this sense, our test is more difficult to be solved because the initial discontinuity is only slightly ($5\%$R) smoothed.

\begin{figure}
\centering 
\begin{tabular}{lr}
			\includegraphics[width=0.47\textwidth]{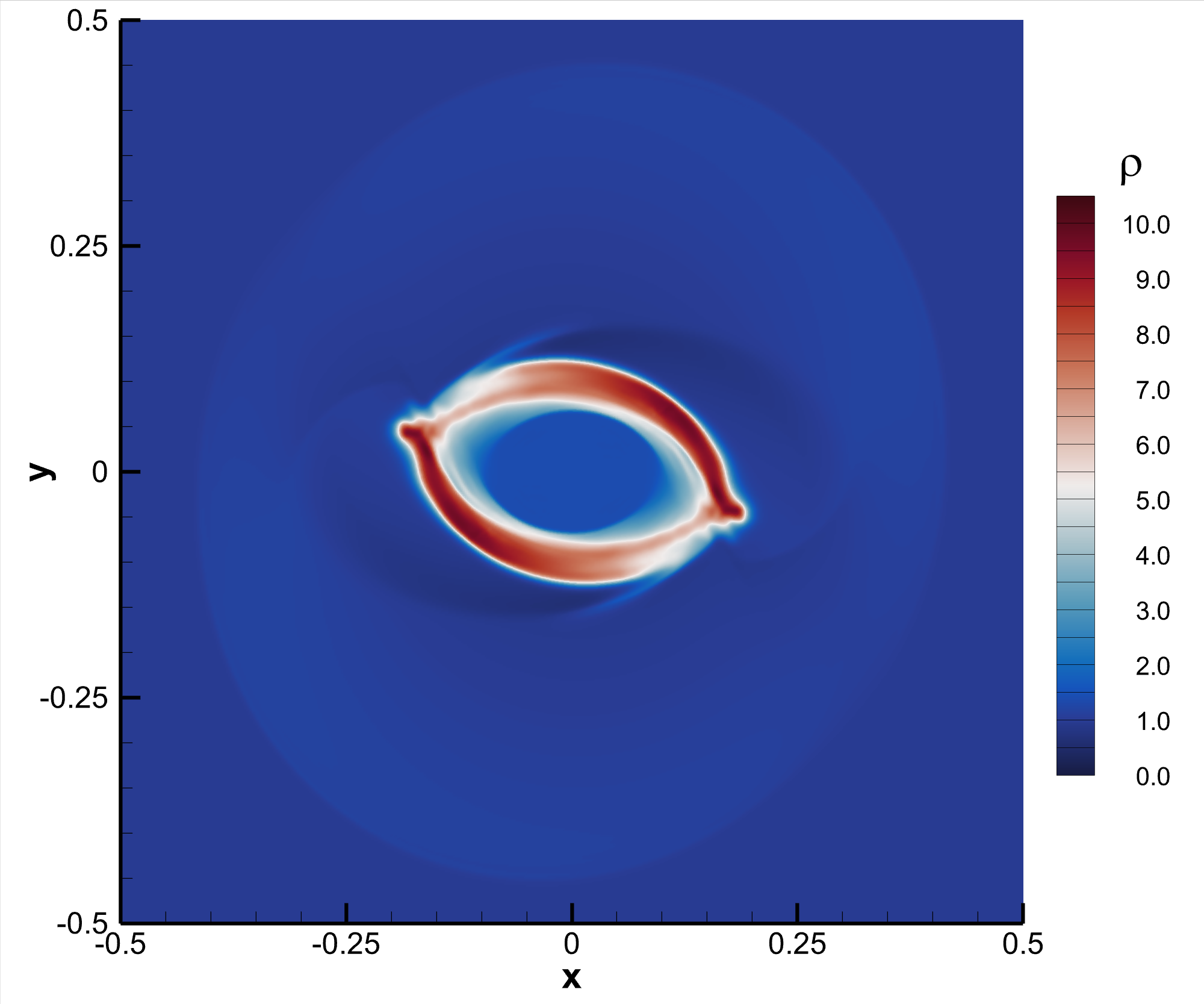} & 
			\includegraphics[width=0.47\textwidth]{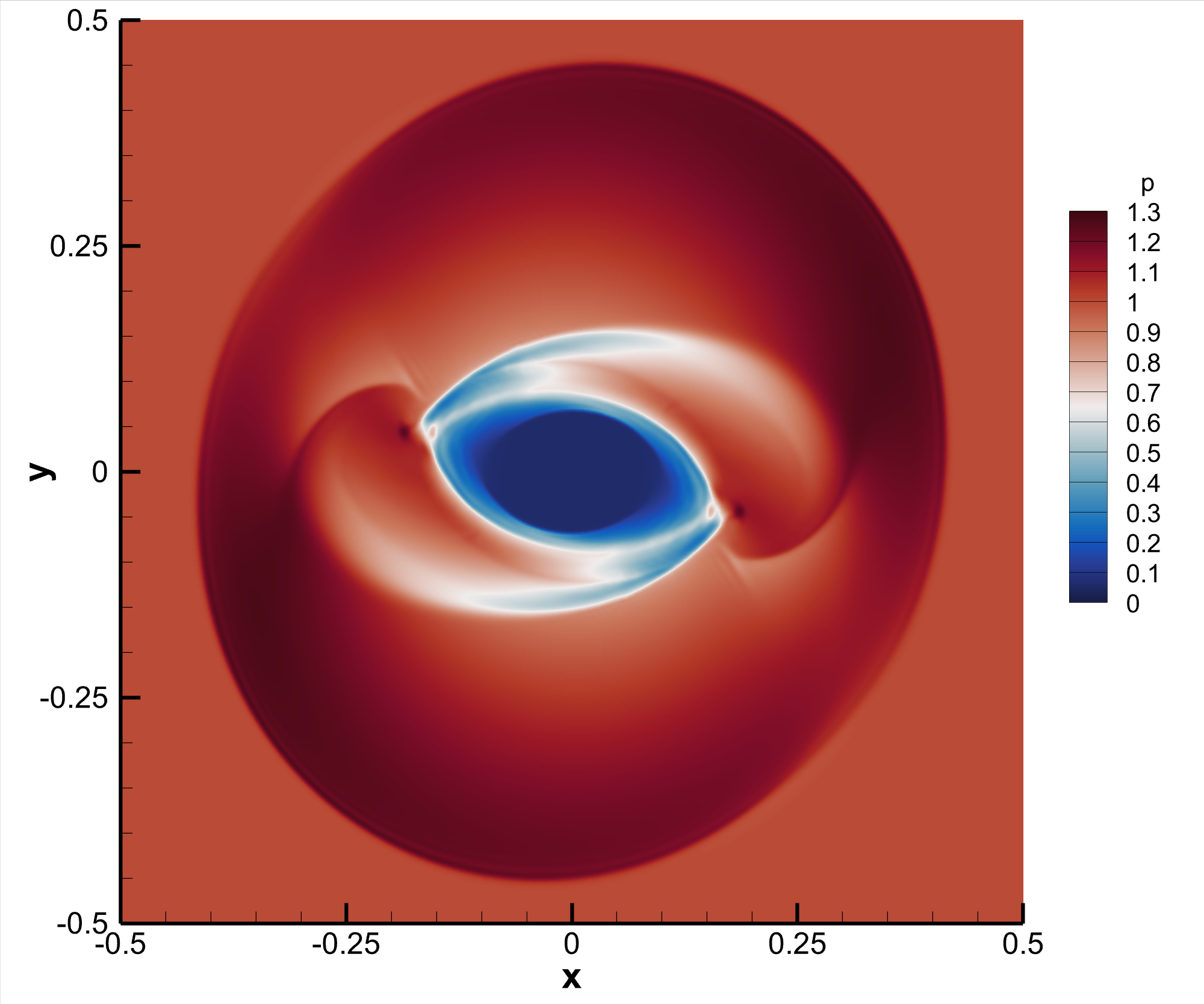}
			\\                                                         
			\includegraphics[width=0.47\textwidth]{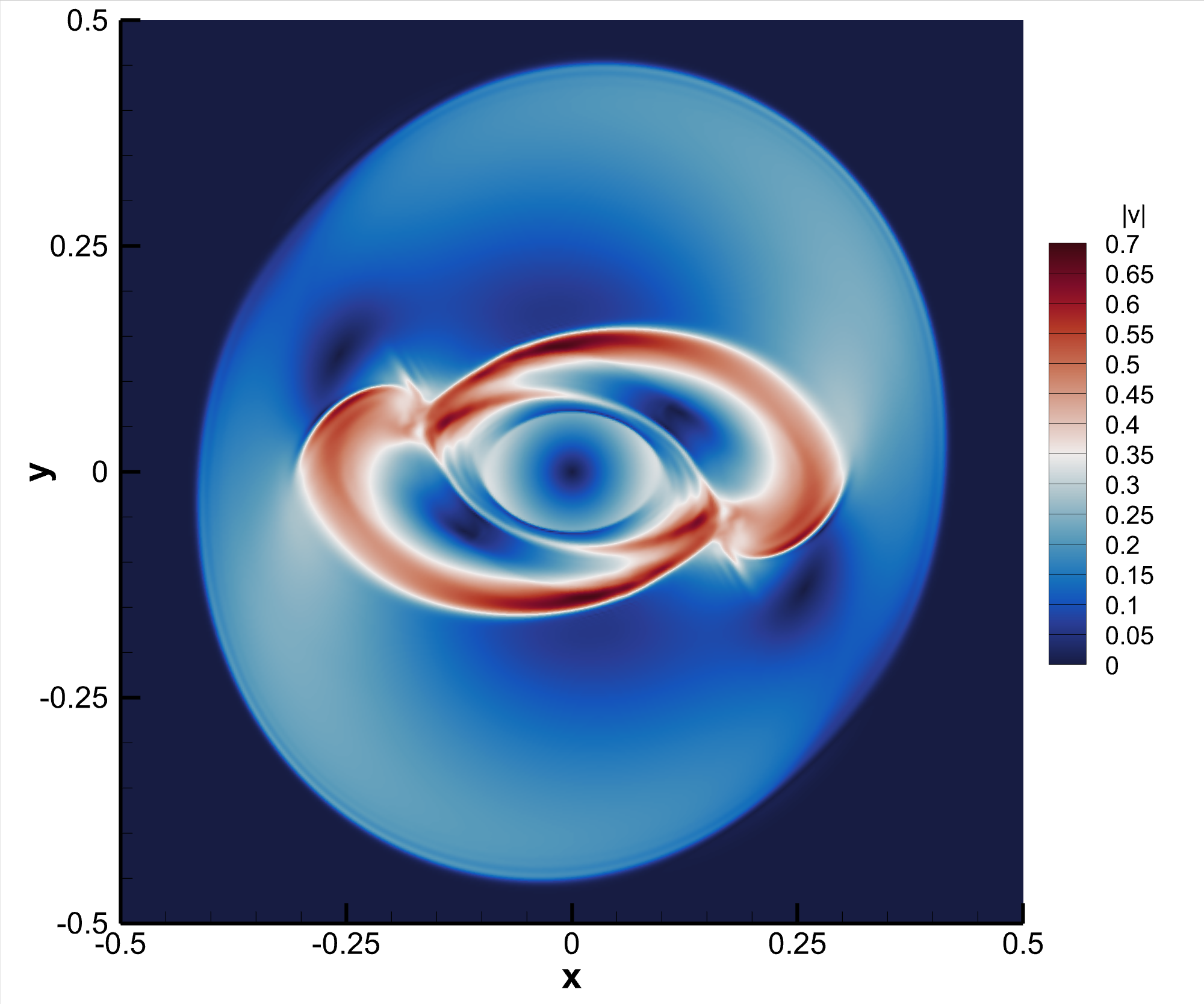} & 
			\includegraphics[width=0.47\textwidth]{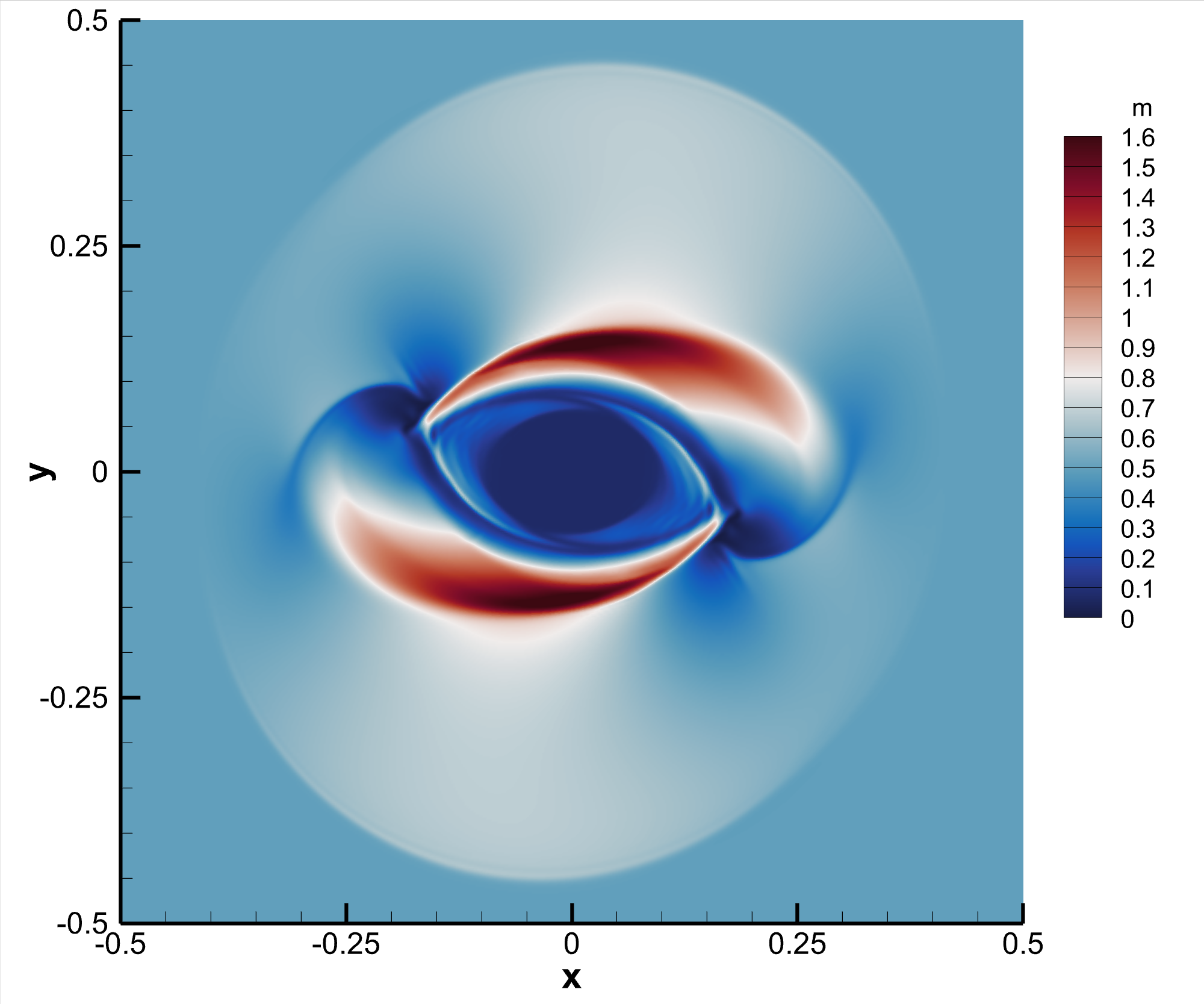} 
	\end{tabular}
\caption{Numerical results for the ideal MHD rotor problem at time $t=0.25$ t.u. obtained with our semi-implicit hybrid FV-FEEC method on a $\Delta x =1/500$ grid, with  security parameter $\CFL=1/4$: two-dimensional contour plots for the  fluid density (top-left), pressure (top-right), absolute value of the velocity  (bottom-left) and magnetic pressure (bottom-right) are shown. In this case, the simulation is ran with the \emph{helicity-preserving} scheme.} \label{fig:Rot}
\end{figure}

Fig. \ref{fig:Rot} shows the contour plots for the numerical solution obtained on a $\Delta x =1/500$ grid, with  security parameter $\CFL=1/4$. In this case, the simulation is ran with the \emph{helicity-preserving} scheme.
Consider that for this test \emph{no-artificial} stabilization term was added to the discrete equations, i.e. $c_h=c_\eta=0$. Moreover, a higher time-resolution is obtained by setting the implicit weights of our FEEC steps to $\theta_p =0.6$ for the acoustic step and Crank-Nicolson $\theta_b=1/2$ for the Alfv\'enic step.

\begin{figure} 
\centering 
\begin{tabular}{cr}
\multicolumn{2}{c}{\includegraphics[width=0.3\textwidth]{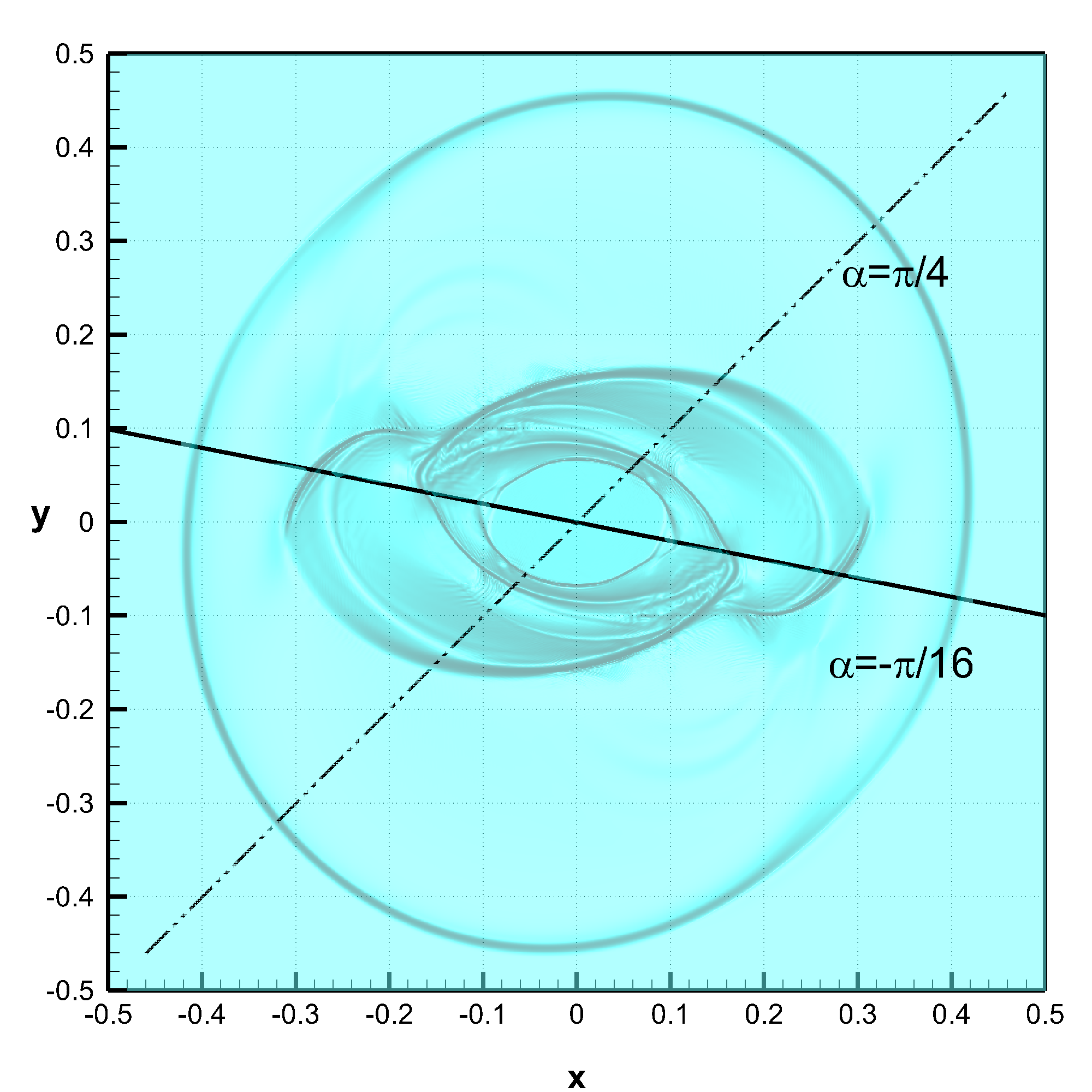} }\\
			\includegraphics[width=0.46\textwidth]{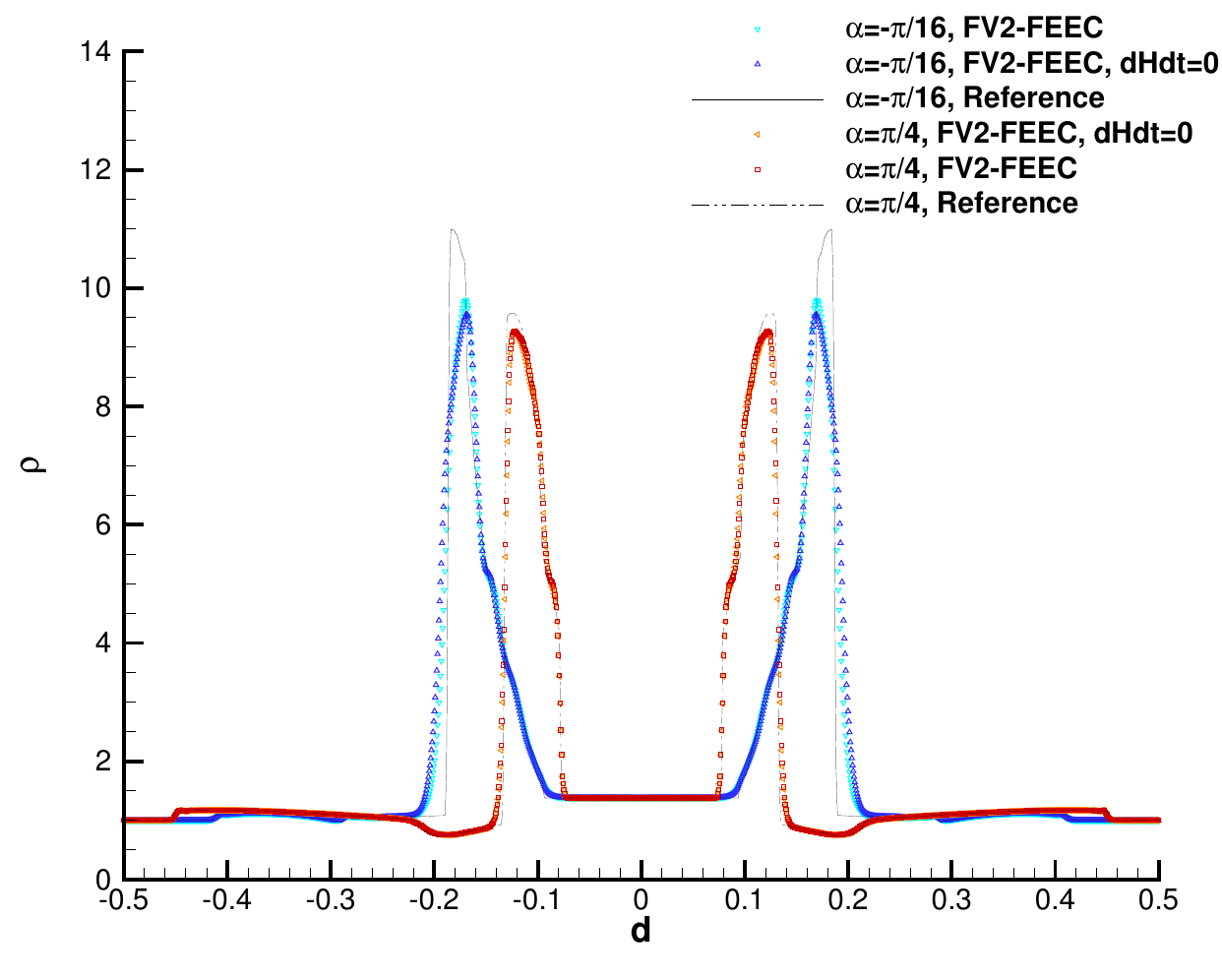} & 
			\includegraphics[width=0.46\textwidth]{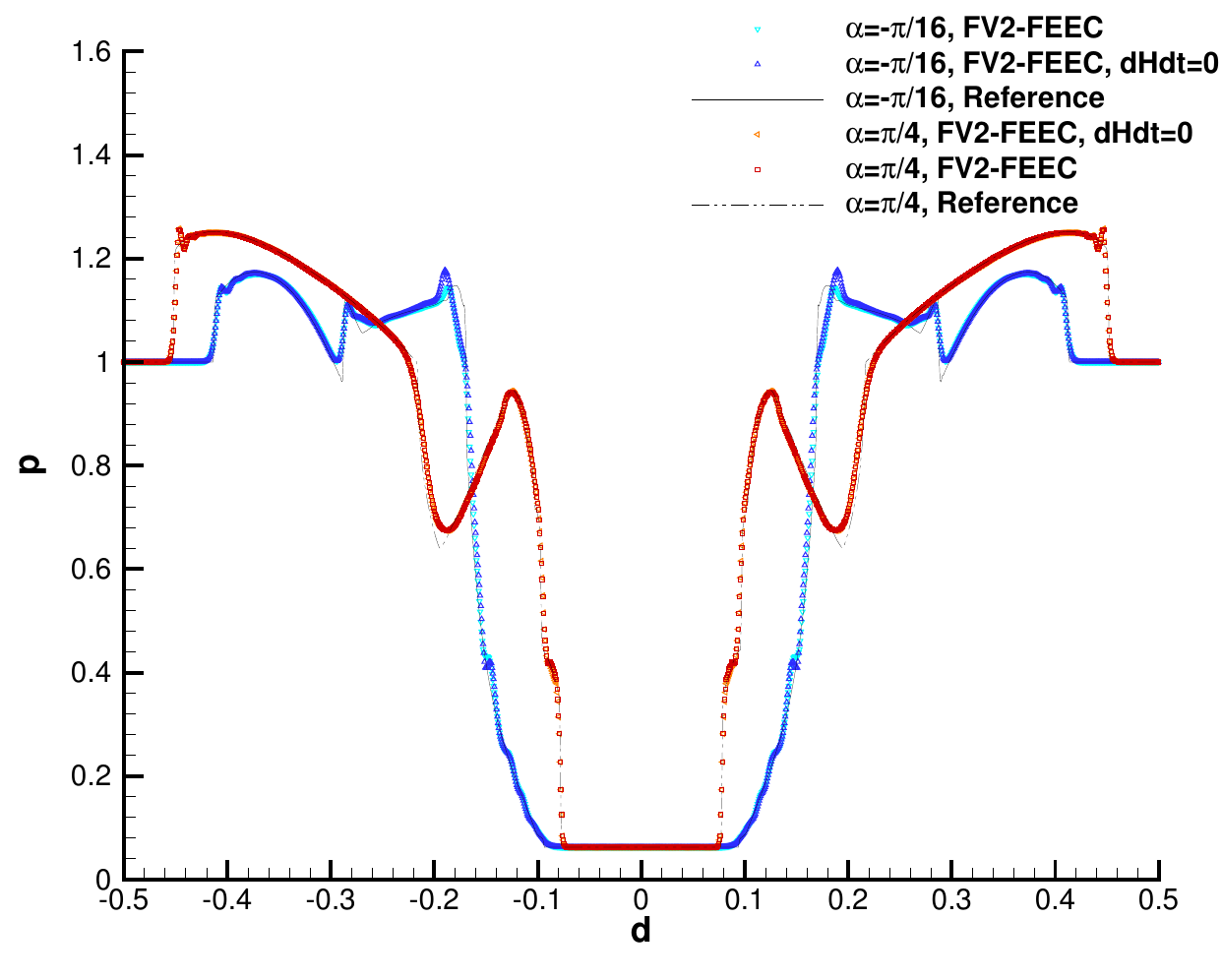}
			\\                         
			\includegraphics[width=0.46\textwidth]{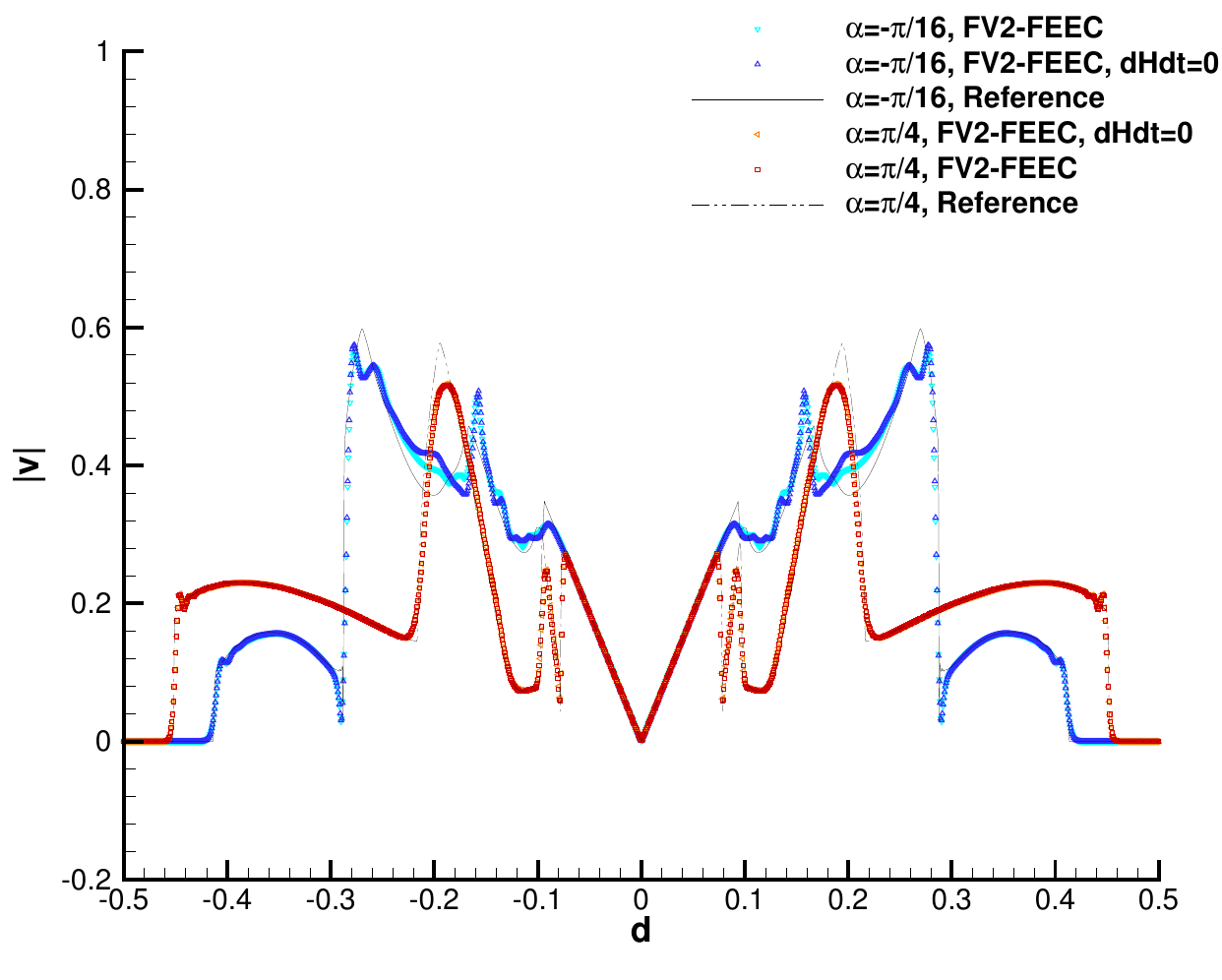} & 
			\includegraphics[width=0.46\textwidth]{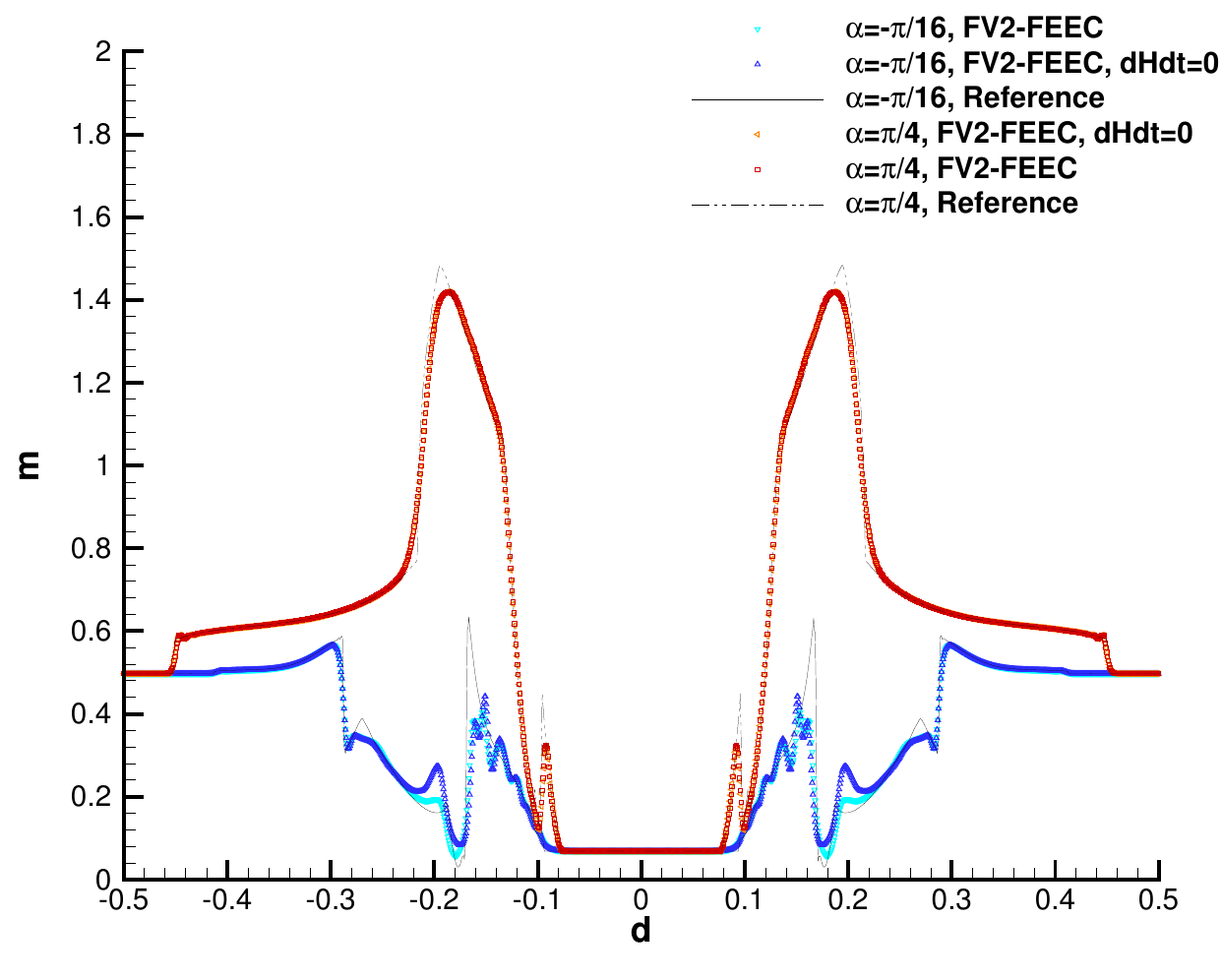} 
	\end{tabular}
\caption{Numerical results for the ideal MHD rotor problem at time $t=0.25$ t.u. obtained with our semi-implicit hybrid FV-FEEC method on a $\Delta x =1/500$ grid, with  security parameter $\CFL=1/4$: the numerical solution is interpolated along the lines $y/x=tg(\alpha)$, $\alpha=\pi/4,-\pi/16$. At the top, the interpolating lines are plotted in the full two-dimensional physical domain. 
Then, the interpolated numerical and reference solution are plotted for the fluid-density (top-left), pressure (top-right), absolute value of the velocity (bottom-left) and magnetic pressure (bottom-right). The reference solution is computed with the aid of a high-order ADER-DG-$\mathbb{P}_5$ scheme with an a-posteriori sub-cell limiter, see \cite{Zanotti2015c}. The numerical solution with and without numerical helicity conservation are compared. The helicity conserving solution is labeled by $dHdt=0$. } \label{fig:Rot_1d}
\end{figure}

\begin{figure} 
\centering 
			\includegraphics[width=0.32\textwidth]{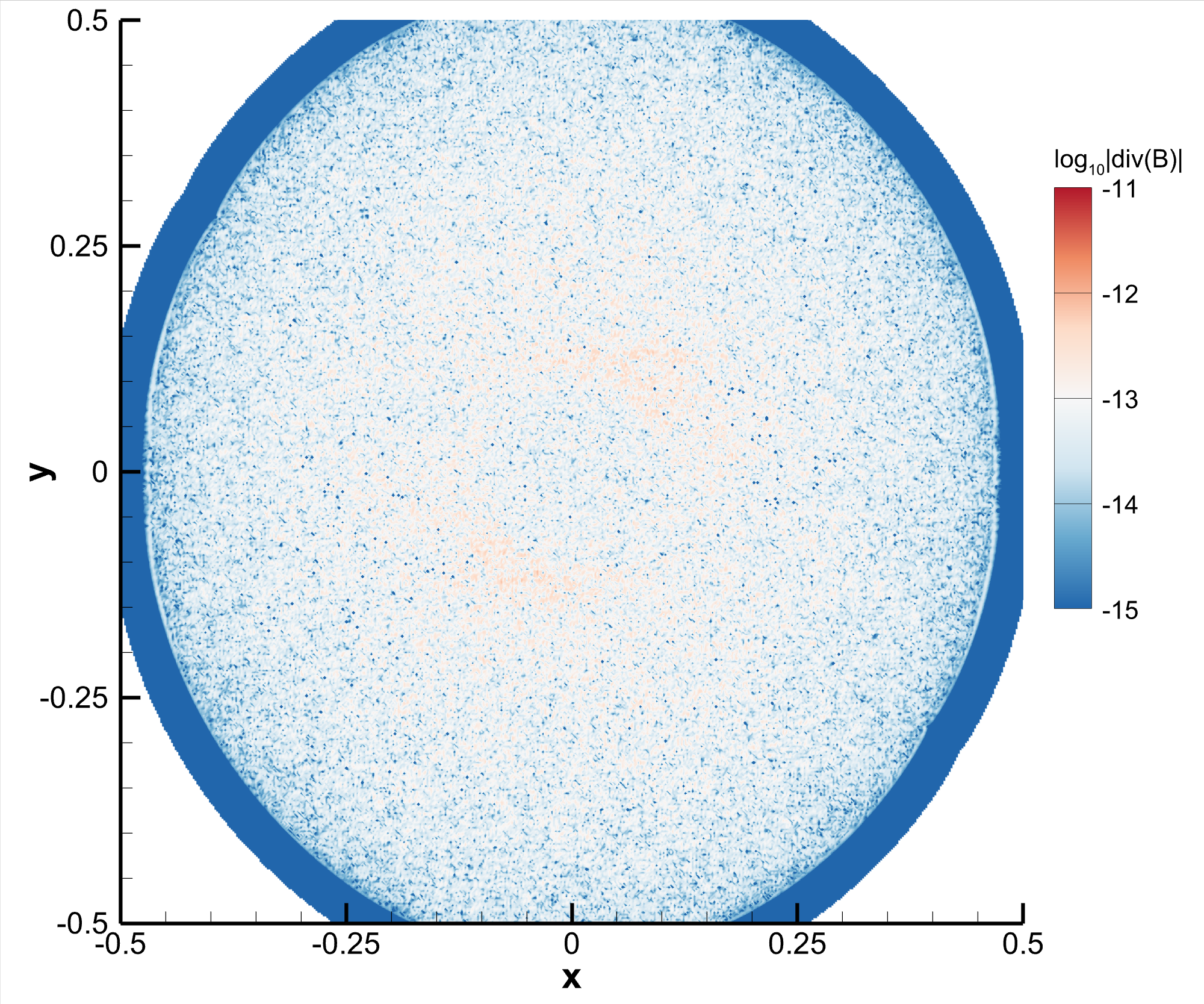}  
			\includegraphics[width=0.32\textwidth]{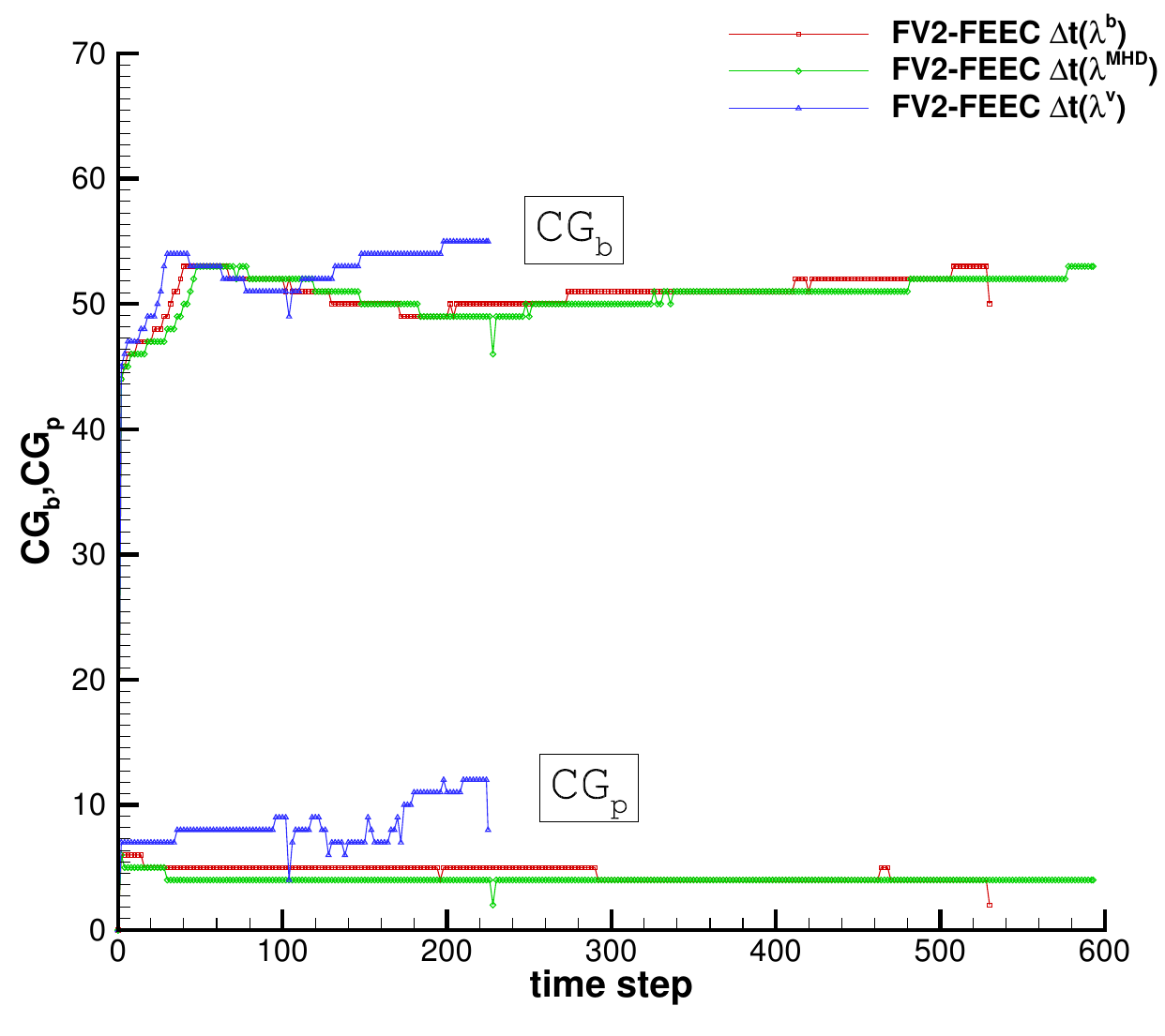}  
			\includegraphics[width=0.32\textwidth]{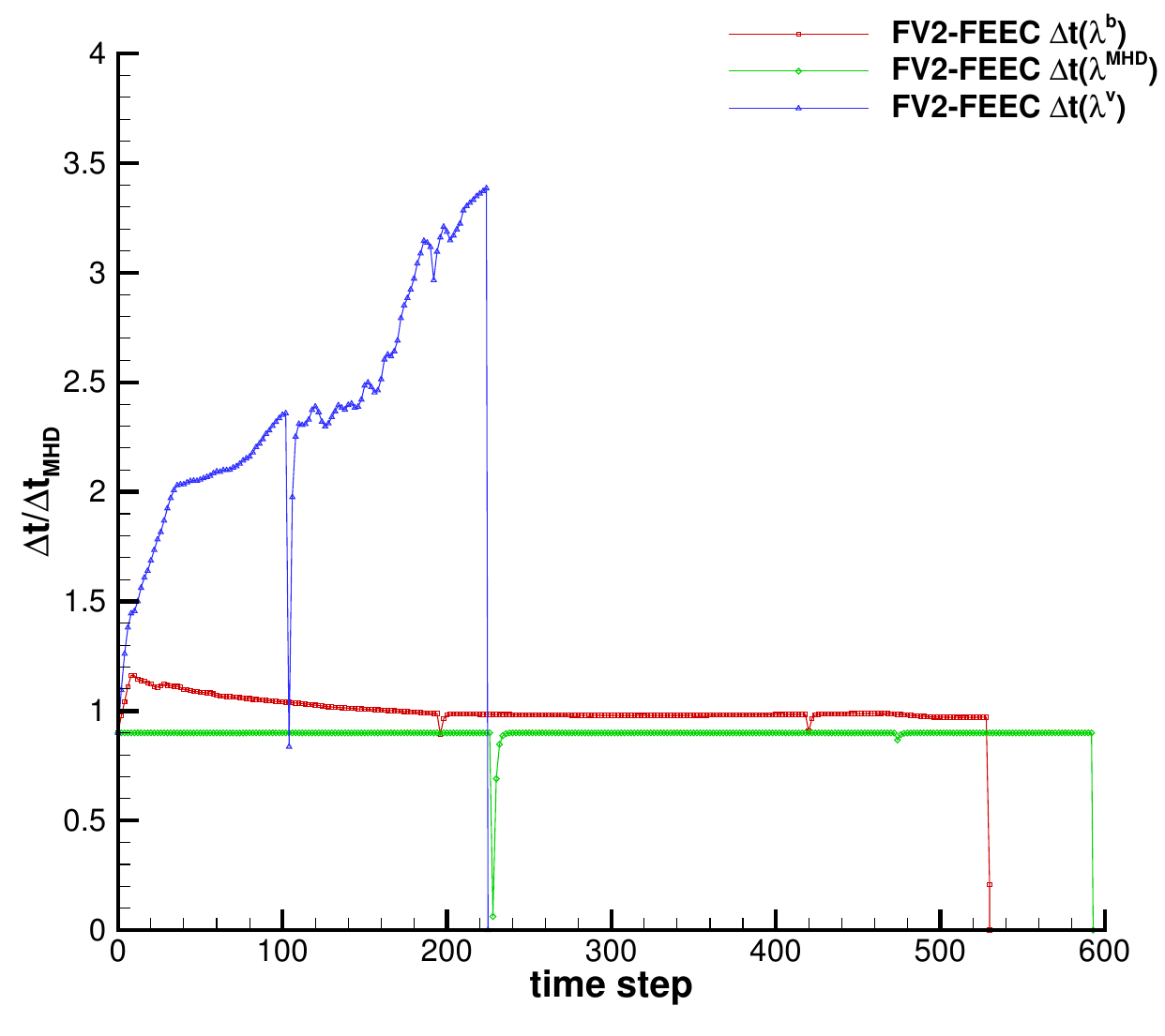}  
\caption{Numerical results for the ideal MHD rotor problem obtained with our semi-implicit hybrid FV-FEEC method on a $\Delta x =\Delta y =1/500$ grid.
At the left,  the divergence error  at time $t=0.25$ t.u. for our helicity conserving method is plotted. For convenience, we set a cut-off for $ div (\B) <10^{-32}$. At the center, the time-step evolution of the number of iterations required by the non-preconditioned and matrix-free conjugate-gradient methods for the Alfv\'enic (CG$_b$) and the acoustic (CG$_p$) system. At the right, the time evolution of the effective Courant number $\Delta t/\Delta t(\lambda^{\MHD})$. }\label{fig:Rot_divB}
\end{figure}

In Fig. \ref{fig:Rot_1d}, the computed FV-FEEC solution, with or without using the magnetic-helicity and -energy preserving strategy, is interpolated along the lines $y/x=tg(\alpha)$, $\alpha=\pi/4,-\pi/16$, and compared with the higher-order accurate ADER-DG (with sub-cell FV limiting) solution of \cite{Zanotti2015c}: all waves are well captured and the numerical results are quite compatible with other published results, despite the current implementation being only second-order accurate. The divergence-error $|div(\B)|$ are of the order $10^{-12}$, see Fig. \ref{fig:Rot_divB}. Fig. \ref{fig:Rot_divB} shows also the time evolution of the number of iterations required by the non-preconditioned matrix-free conjugate gradient  method for the Alfv\'enic (CG$_b$) and acoustic (CG$_p$). One can see that, for this test, the number of CG iterations are almost constant, and using larger time-step definitely helps in computational efficiency.

\subsection{Viscous and resistive Orszag-Tang vortex}\label{sec:VROT}

\begin{figure} 
\centering 
 \begin{tabular}{cc}   	 
 		\includegraphics[width=0.45\textwidth]{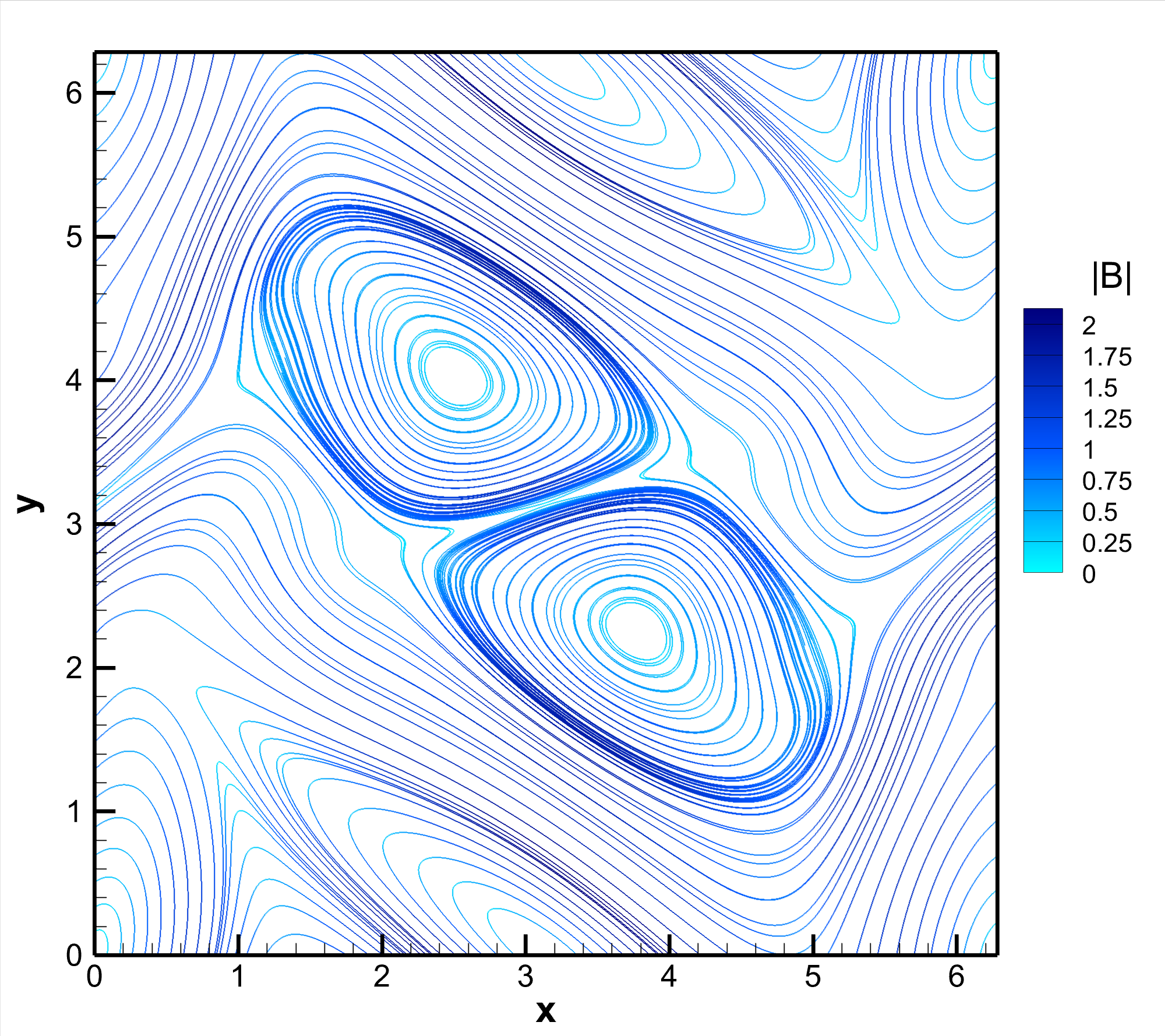} &  
 		\includegraphics[width=0.45\textwidth]{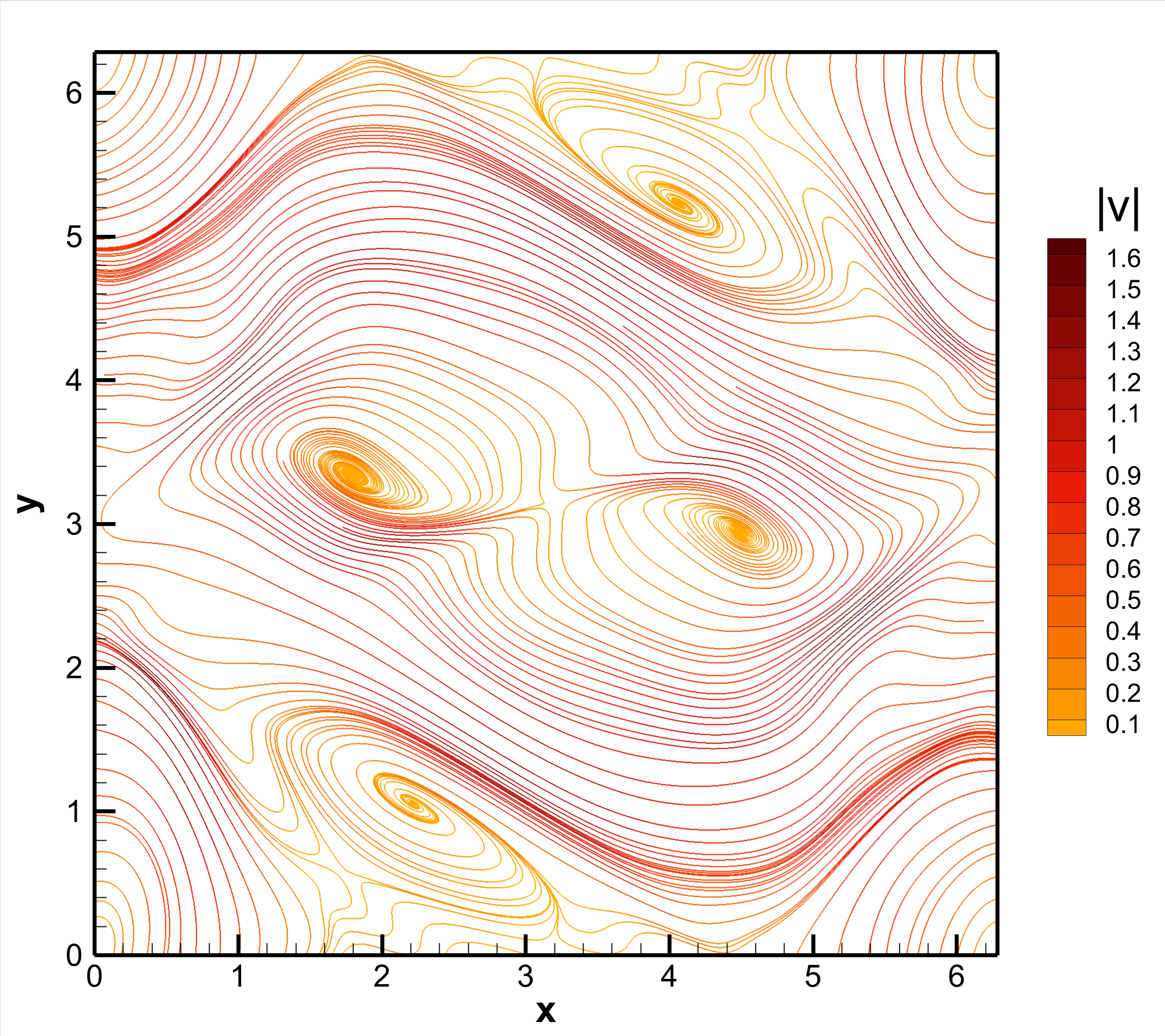}  
 	\end{tabular}
\caption{Numerical results for the 2d viscous and resistive Orszag-Tang vortex system at time $t=2$ t.u. obtained with our semi-implicit hybrid FV-FEEC method on a $\Delta x=\Delta y=1/200$ grid. The field lines of magnetic field and velocity are plotted} \label{fig:VROT}
\end{figure}

The viscous and resistive Orszag-Tang vortex is the same test of above (see Sec. \ref{sec:OT}) after considering also physical dissipation. In particular, the initial condition is
\begin{equation} \label{eq:VROrszagTang_ic}
\left\{
  \begin{array}{rl} 
\rho &= 1, \\
\v  &=  \left( - \sin\left(y\right), \sin \left(x \right), 0 \right),\\
p &=\frac{15}{4} + \frac{1}{4} \cos(4x) + \frac{4}{5} \cos(2x) \cos(y) - \cos(x) \cos(y) + \frac{1}{4} \cos(2y),\\
\B &=   \left( -  \sin\left(y\right), \sin \left(2x \right), 0\right)/\sqrt{4 \pi} .  
  \end{array}
	\right.
\end{equation} 
 with parameters $\gamma = \frac{5}{3}$, $\mu = \eta = 10^{-2}$, $c_v=1$ and a Prandtl number of $Pr=1$, see \cite{WarburtonVRMHD,ADERVRMHD}.
The domain is a periodic two-dimensional box $(x,y)\in[0,2\pi]^2$, and the solution space $V_0$ is discretized with a $\Delta x=\Delta y=1/200$ grid. Second order reconstruction is used for the finite-volume fluxes and additional artificial stabilization is needed, i.e. $c_h=c_\eta = 0$. Regarding the time-discretization, the security parameter is set to $\CFL=0.9$, combined with a time-step size related to the convective time-scale $\Delta t(\lambda^v)$, and implicit weights $\theta_b=\theta_p=0.55$, with magnetic -energy and -helicity \emph{stable} discretization. The final time is $t_f=2$. 
For this test, velocity and magnetic field streamlines are   plotted in Fig. \ref{fig:VROT}. Velocity streamlines cross each other because of the fluid compressibility, while the magnetic streamlines remains almost all closed because of the divergence-preserving property.

\subsection{3d Viscous and resistive MHD Orszag-Tang vortex}
\begin{figure} 
\centering 
			 \includegraphics[width=0.45\textwidth]{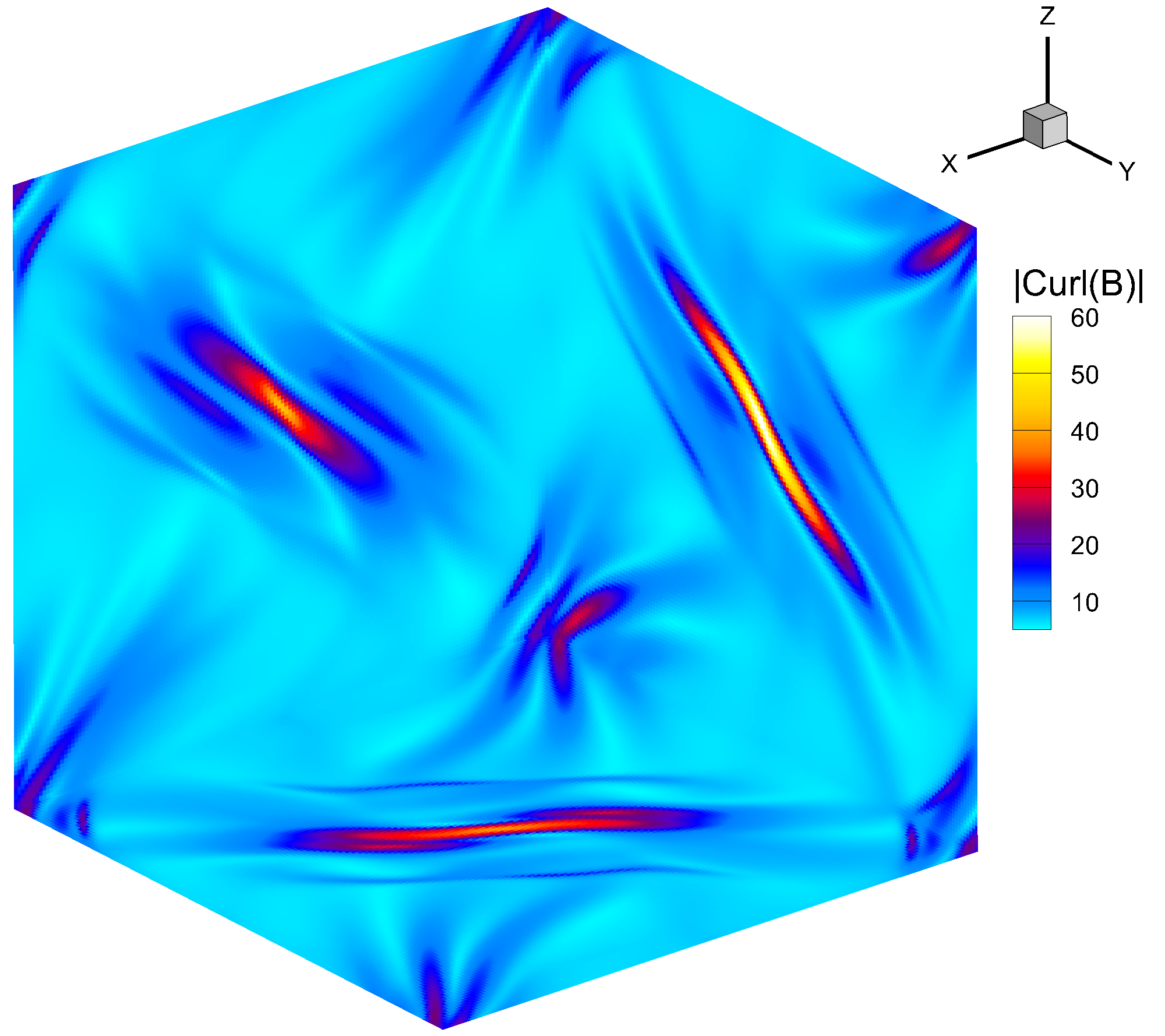}  
			 \includegraphics[width=0.45\textwidth]{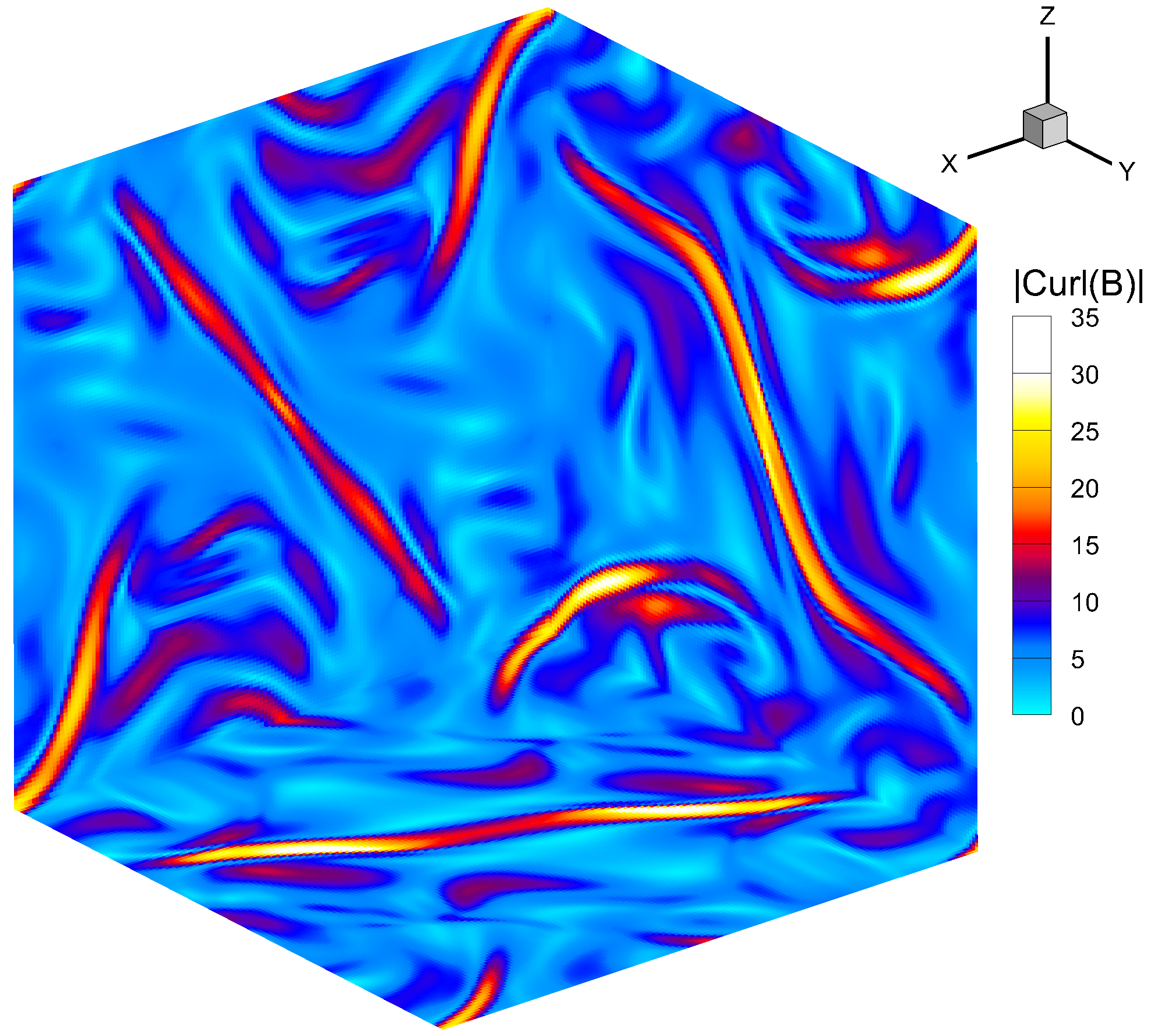}	 \\
			 \includegraphics[width=0.45\textwidth]{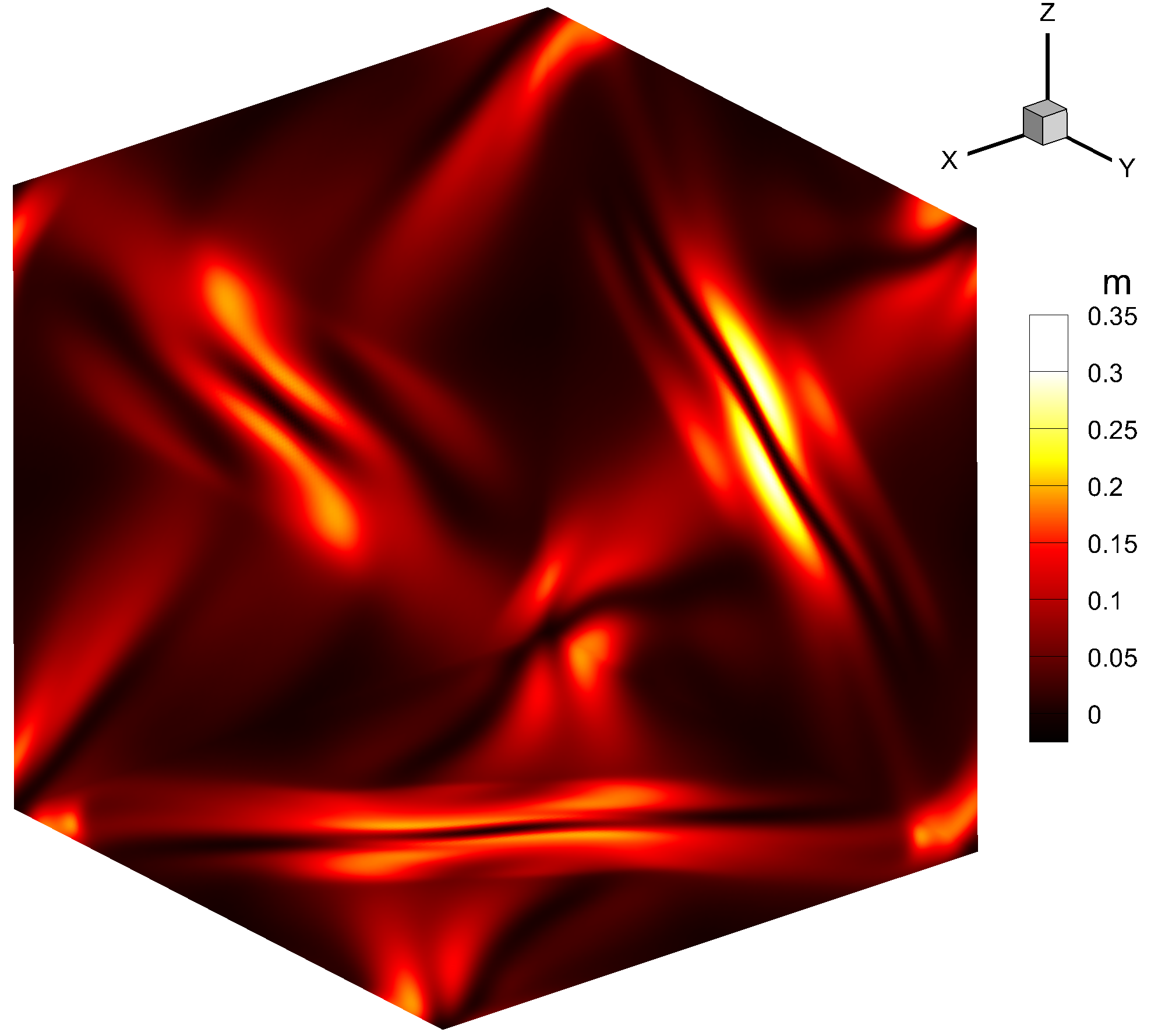}  
			 \includegraphics[width=0.45\textwidth]{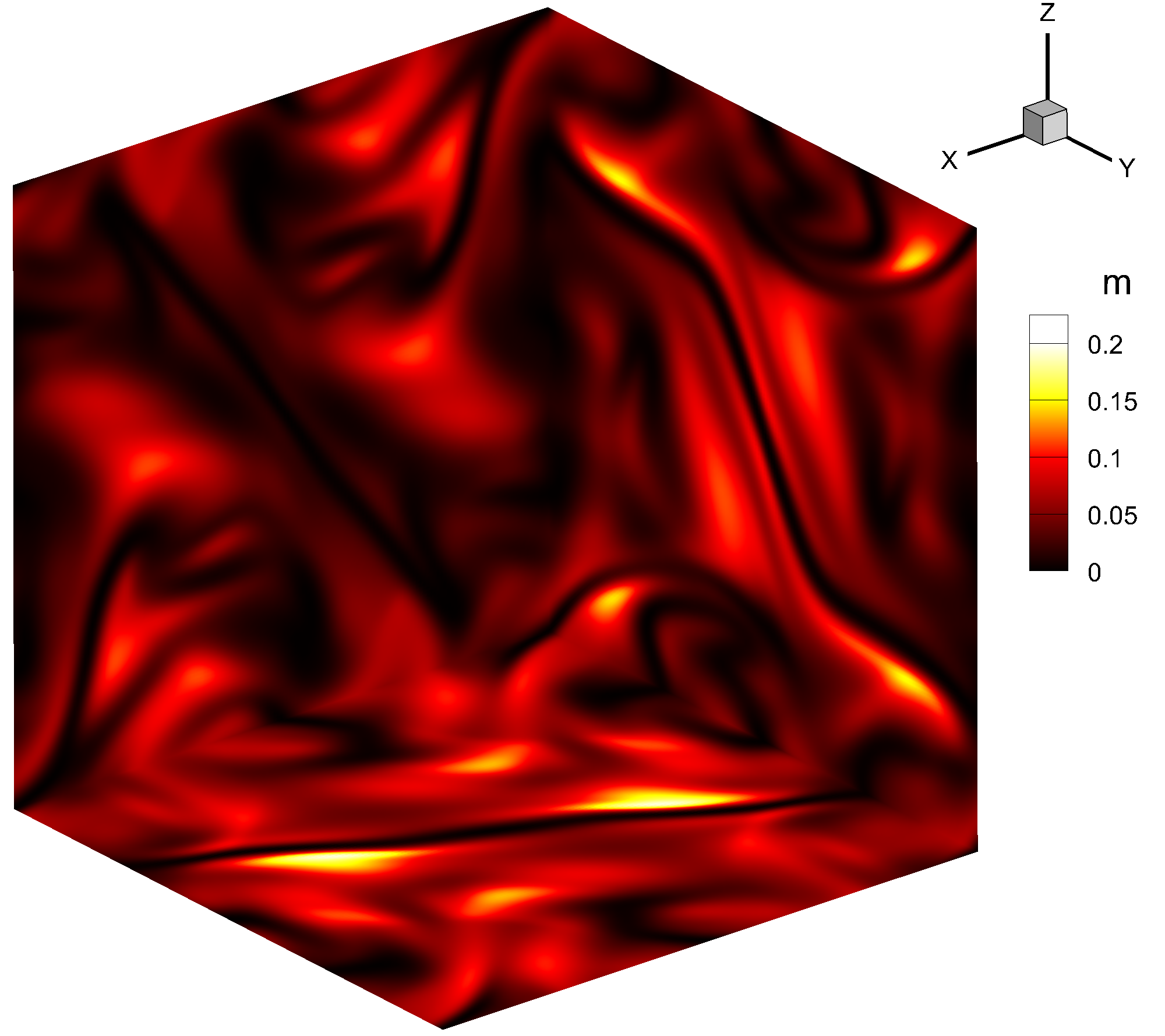}	 \\
			 \includegraphics[width=0.45\textwidth]{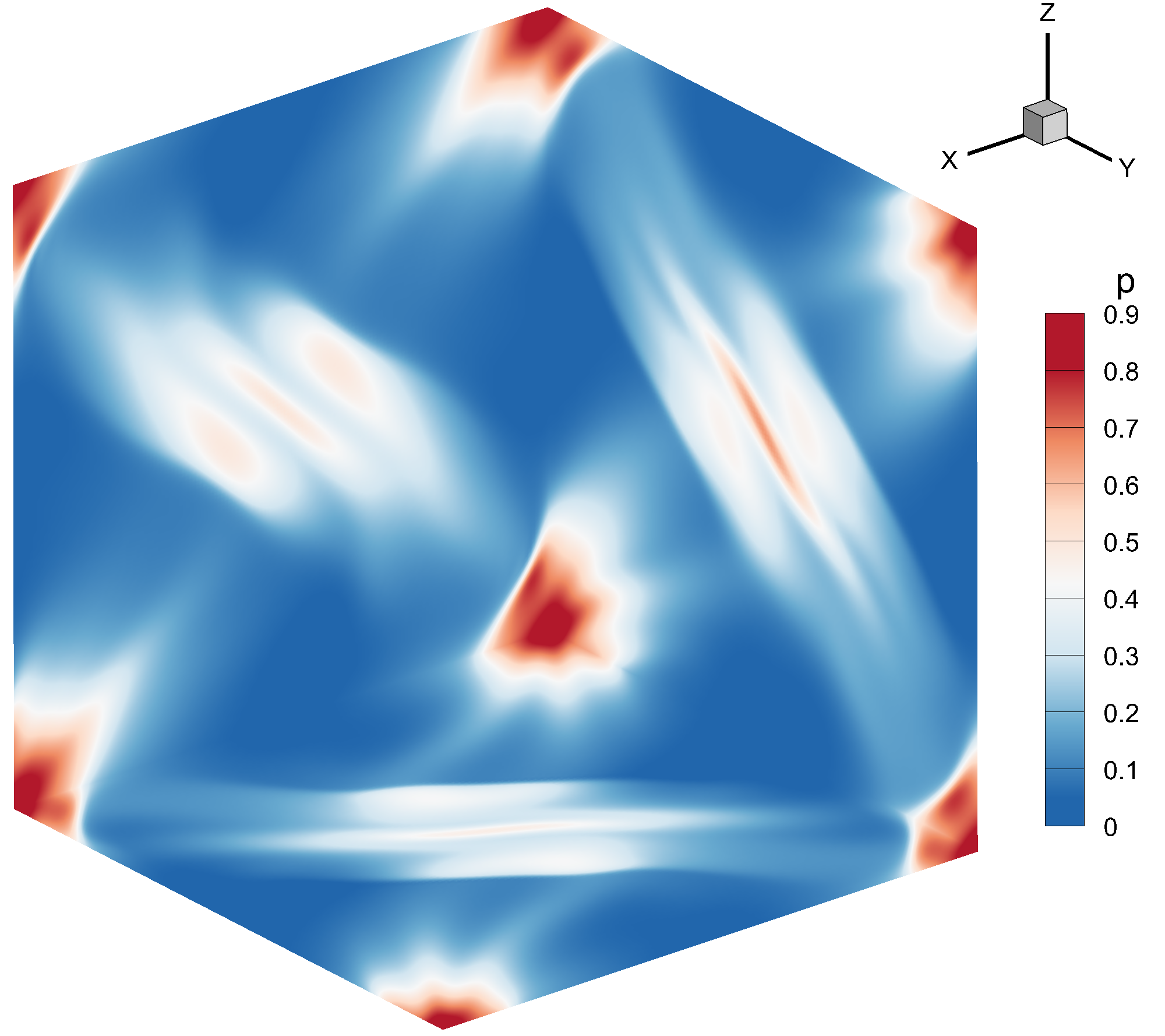}  
			 \includegraphics[width=0.45\textwidth]{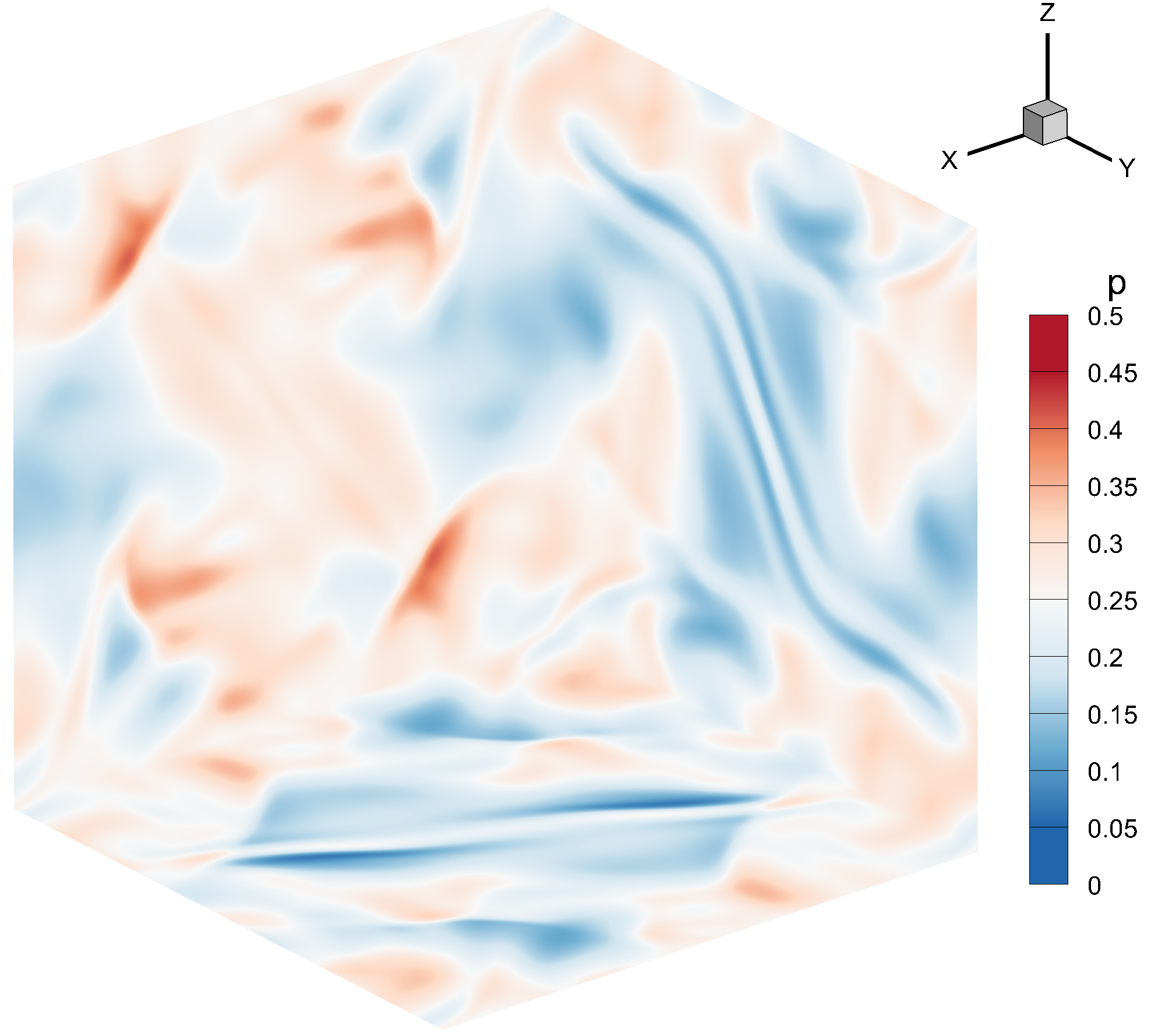}	 
\caption{Numerical results for the 3d viscous and resistive Orszag-Tang vortex system at time $t=0.25$ (left) and $0.5$ (right) t.u. obtained with our semi-implicit hybrid FV-FEEC method on a $\Delta x= \Delta y = \Delta z = 1/150$ grid. Contour plots of the numerical solution interpolated along the three orthogonal 2D planes $z=0$, $y=0$ and $x=0$ are shown: the magnitude of the magnetic current $curl(\B)$ (first row), magnetic energy $m$ (second row)  and fluid pressure (third row).} \label{fig:VROT3d}
\end{figure}

\begin{figure} 
\centering 
			 \includegraphics[width=0.45\textwidth]{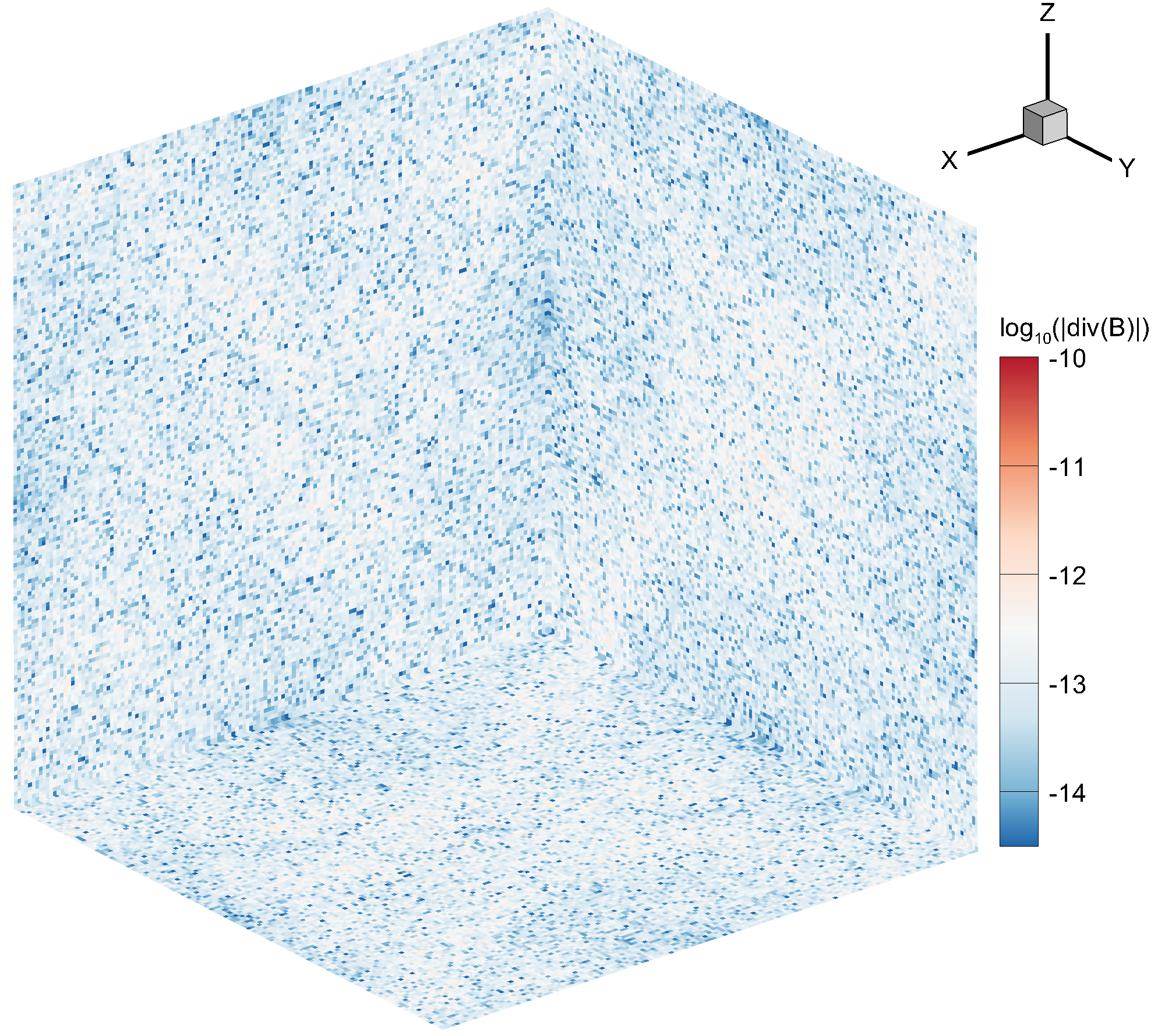}\includegraphics[width=0.45\textwidth]{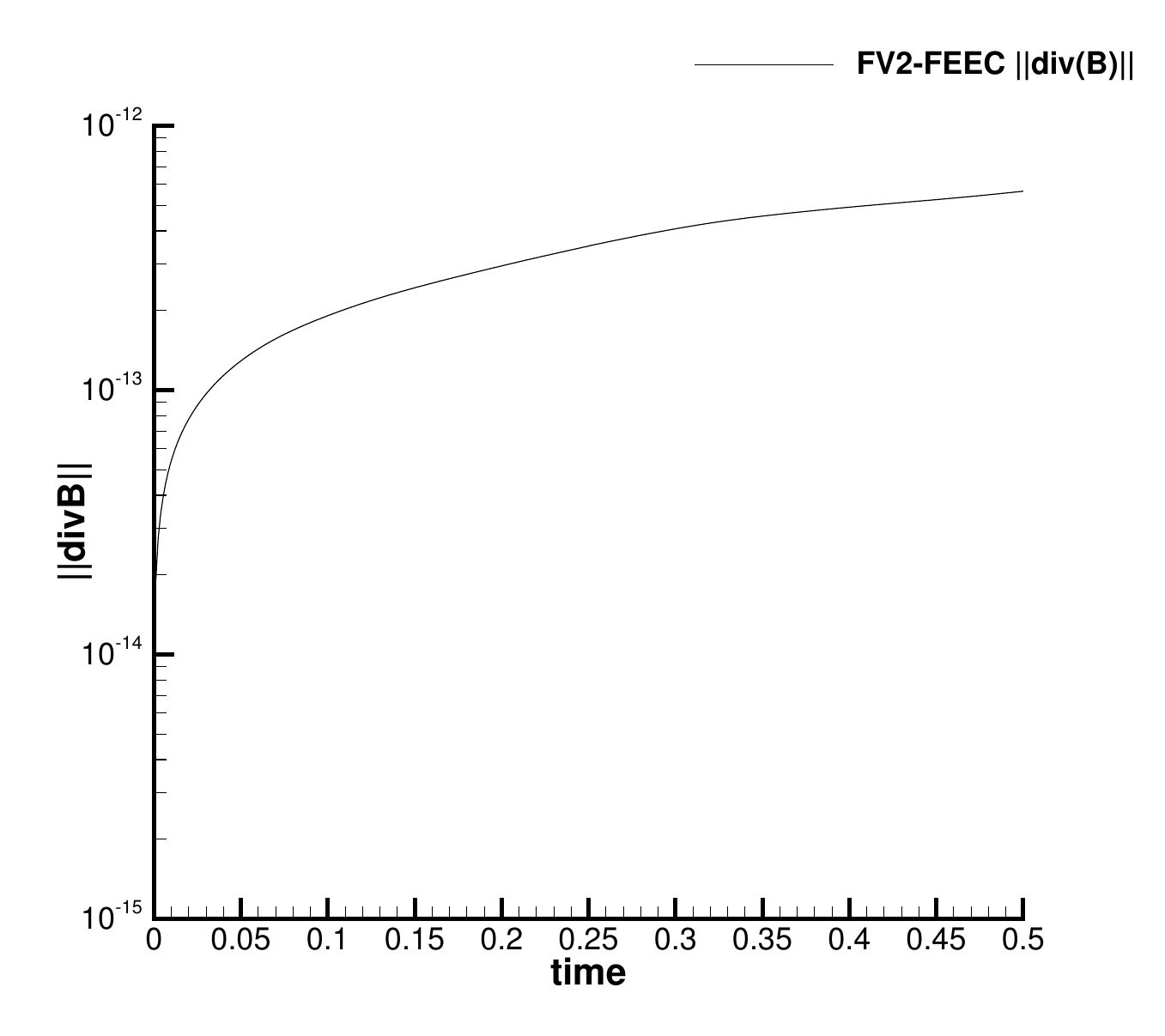}
\caption{Numerical divergence errors for the 3d viscous and resistive Orszag-Tang vortex system at time $t=0.5$ (right) t.u. obtained with our semi-implicit hybrid FV-FEEC method on a $\Delta x= \Delta y = \Delta z = 1/150$ grid. Contour plots of the numerical solution interpolated along the three orthogonal 2D planes $z=0$, $y=0$ and $x=0$ are shown at the left, the time evolution of the norm of divergence errors at the right.} \label{fig:VROT3d_divB}
\end{figure}

In this section, the numerical results for the three-dimensional viscous and resistive MHD Orszag Tang vortex are presented. This test is the three-dimensional extension of the previous one, see Sec. \ref{sec:VROT}. 
The initial condition comes from  \cite{BOHM2018}, and \cite{PopovElizarova2015} for ideal MHD, i.e.
\begin{equation} \label{eq:VROrszagTang3d_ic}
\left\{
  \begin{array}{rl} 
\rho &= \frac{25}{36 \pi}, \\
\v  &=   \left( - \sin\left(2\pi z\right), \sin \left(2\pi x \right), \sin \left(2\pi y \right)\right),\\
p &=\frac{5}{12 \pi },\\
\B &=   \left( -  \sin\left(2\pi z\right), \sin \left(4\pi x \right), \sin \left(4\pi y \right)\right).  
  \end{array}
	\right.
\end{equation} 
This time, the flow parameters are chosen to be  $\gamma = \frac{5}{3}$, $\mu = 6\times 10^{-6}$, $\eta = 10^{-3}$, $c_v=1$ and a Prandtl number of $Pr=0.72$. The physical domain is a periodic 3d box $(x,y,z)\in[0,1]^3$, and the solution space $V_0$ is discretized with a $\Delta x= \Delta y = \Delta z = 1/150$ grid.
The time-discretization is set up with a security parameter $\CFL=0.9$, implicit weights $\theta_b=\theta_p=0.65$, and time-step sizes of the convective time-scale $\Delta t (\lambda^v)$. Note that the fluid viscosity is very low, and the flow may become highly nonlinear in the three-space dimensions. Contour plots interpolated along the three planes $(x,y,z=0)$, $(x,y=0,z)$ and $(x=0,y,z)$ the magnitude of magnetic current $|curl(\B)|$, magnetic energy $m$ and fluid pressure $p$ are reported in Fig. \ref{fig:VROT3d}, the obtained divergence error are of order $10^{-12}$, see Fig. \ref{fig:VROT3d_divB}.


\section{Conclusions}
\label{sec:conclusions}

We designed a novel semi-implicit hybrid FV-FE method for the compressible and nonlinear viscous and resistive MHD equations. The nonlinear convective terms are thought to be solved by means of any very efficient and robust FV scheme, while the magneto acoustic terms are discretized with compatible FE based on a discrete and continuous de Rham complexes (FEEC). Thanks to the FEEC framework and operator-splitting, the implicit steps consist in solving two decoupled linear, symmetric and positive definite algebraic systems: one for the scalar pressure (acoustic step); one for the fluid velocity (Alfv\'enic step).
The time-discretization is an adaptation of the \emph{alternating direction implicit} (ADI) method to a three-operator splitting of the type \emph{explicit-implicit-implicit}. 
The final scheme is shown to be robust and accurate even at high Courant numbers, with time-steps based only on the convective time-scale $\Delta t(\lambda^v)$.
To validate the method, a second-order accurate FV-FE implementation is tested against a non-trivial set of MHD problems: in particular, shock-dominated flows even with Mach $25.5$ shocks, are accurately solved; a stationary MHD vortex is preserved with low dissipation at long times, and two- and three-dimensional nonlinear tests are solved to prove the robustness, low-dissipation and stability properties of the final algorithm.

\noindent
The numerical formulation is quite general, and it may be easily generalized to higher-order accurate implementations, eventually by considering a robust Discontinous Galerkin discretization of the convective part, see e.g. \cite{Dumbser2008,DumbserNSE}. Moreover, thanks to the high flexibility of finite-volume and finite-element methods, the extension to adaptive grids is currently under investigation, eventually introducing high-order broken FE spaces, see e.g. \cite{campos_pinto_gauss-compatible_2016,brokenFEECNS,gucclu2022broken}. 
We believe that the presented algorithm can be easily adapted to be rigorously asymptotic-preserving, and is a good candidate to solve multi-scale problems in plasma-physics, e.g. for  applications to astrophysics or magnetic confinement fusion research. Indeed, the presented schemes performs very well against stationary or slowly varying nonlinear MHD configurations, as well as against  fast-modes dominated flows, also in presence of high Mach MHD shocks.

\appendix
\section{Computing the vector potential} \label{app:A}
For the simple purpose of verifying the conservation of magnetic helicity, we report  here an efficient solver for the vector potential that can be built by using exactly the same FEEC discrete space of solutions and operators that are presented in this manuscript. In particular, we computed the vector potential $\mathbf{A}_h$, imposing $\nabla \cdot \mathbf{A}_h=0$, as the equilibrium solution of the following parabolic-elliptic equation with source term
\begin{align} \label{eq:dAdt}
& \frac{d}{d t}  \int \mathbf{w}_h \cdot  \mathbf{A}_h + \int \nabla \times   \mathbf{w}_h  \cdot  \nabla \times  \mathbf{A}_h  + \int \mathbf{w}_h \cdot \nabla \psi_h = \int \nabla \times   \mathbf{w}_h  \cdot  \mathbf{B}_h, \quad \forall \mathbf{w}_h\in V_1, \\ \label{eq:divA}
&\int  \mathbf{A}_h \cdot \nabla q_h  = 0, \quad \forall q_h \in V_0
\end{align} 
where $\psi_h \in V_0$ is the GLM multiplier that corresponds to the divergence-free constraint  $\nabla \cdot \mathbf{A}_h = 0$. 
The equation is solved by using again an implicit splitting between the $curl-curl$ parabolic term and the $div-grad$ coupling with the divergence-free condition. 
One may note that,  $\mathbf{A}_h$ can be interpreted as the equilibrium solution of the incompressible Navier-Stokes equation with a stationary source term. 
Also in this case, we defined the numerical equilibrium as $\left\|\mathbf{A}^{n(\tilde\epsilon)+1}-\mathbf{A}^{n(\tilde\epsilon)} \right\| \leq \tilde\epsilon \ll 1$.
This equation has been solved only to show that the presented algorithm is helicity conserving within the prescribed tolerances. For the numerical test, we choose $\tilde\epsilon=10^{-12}$. Here we assumed periodic boundary conditions.

\section*{Author contributions} 
\noindent
\textbf{F. Fambri and E. Sonnendr\"ucker}:  
Conceptualization;Methodology;Writing - original draft; Writing - review \& editing; Software; Validation.

\section*{Acknowledgments}
We would like to acknowledge Omar Maj for the inspiring discussions on the topic. \\ 
F.F. is member of the Gruppo Nazionale Calcolo Scientifico-Istituto Nazionale di Alta Matematica (GNCS-INdAM).

\bibliographystyle{plain}
\bibliography{references}


 
\end{document}